\documentclass[a4paper,twoside]{article}

\usepackage{amsmath,amsthm,amsfonts,latexsym,amscd,amssymb,enumerate}

\usepackage[hypertex,backref]{hyperref}

\swapnumbers
\theoremstyle{plain}
\newtheorem{theorem}{Theorem}[section]
\newtheorem{lemma}[theorem]{Lemma}
\newtheorem{sublemma}[theorem]{Sublemma}
\newtheorem{corollary}[theorem]{Corollary}
\newtheorem{proposition}[theorem]{Proposition}

\theoremstyle{definition}
\newtheorem{definition}[theorem]{Definition}
\newtheorem{example}[theorem]{Example}
\newtheorem{notation}[theorem]{Notation}

\newtheorem{set-up}[theorem]{Geometric set-up}
\newtheorem{remark}[theorem]{Remark}

\newcommand{\reals}{\mathbb{R}}
\newcommand{\complexs}{\mathbb{C}}
\newcommand{\naturals}{\mathbb{N}}
\newcommand{\integers}{\mathbb{Z}}
\newcommand{\rationals}{\mathbb{Q}}

\newcommand{\NeumannN}{\mathcal{N}}
\DeclareMathOperator{\id}{id}
\newcommand{\boundary}[1]{\partial#1}
\newcommand{\boundedops}{\mathcal{B}}

\newcommand{\abs}[1]{\left\lvert#1\right\rvert} 
\newcommand{\norm}[1]{\left\lVert#1\right\rVert}

\newcommand{\tensor}{\otimes}
\newcommand{\into}{\hookrightarrow}

\newcommand{\iso}{\cong}

\newcommand{\disjointunion}{\amalg}

\DeclareMathOperator{\im}{im}      
\DeclareMathOperator{\vol}{vol}    
\DeclareMathOperator{\Hom}{Hom}    

 \DeclareMathOperator{\tr}{tr}
\DeclareMathOperator{\TR}{TR}
 \DeclareMathOperator{\pr}{pr}
\DeclareMathOperator{\coker}{coker}

\DeclareMathOperator{\ind}{ind}
\DeclareMathOperator{\Ind}{Ind}
\DeclareMathOperator{\sgn}{sgn}

\newcommand{\forget}[1]{}

\newcommand{\innerprod}[1]{\langle #1 \rangle}
\newcommand{\Commutator}[2]{[ #1,#2]}

\def  \nuint {\raise10pt\hbox{$\nu$}\kern-6pt\int}

\newcommand\STR{\operatorname{STR}}
\newcommand\bSTR{\operatorname{bSTR}}

\newcommand\Tr{\operatorname{Tr}}
\newcommand\str{\operatorname{str}}

\def\N{\mathcal N}
\def \L{\mathcal L}

\newcommand\E{\mathcal E}
\newcommand\Q{\mathcal Q}

\newcommand\V{\mathcal V}
\newcommand\C{\mathcal C}

\def \L {{\cal L}}
\def \Sp {{\cal S}}
\newcommand\B{\mathcal B}

\def \H {{\cal H}}

\def \Ch {{\rm Ch}}
\def\Id{{\rm Id}}

\newcommand\ch{\operatorname{ch}}

\renewcommand\Im{\operatorname{Im}}

\newcommand\D{\mathcal D}
\newcommand\Di{D\kern-6pt/}
\newcommand\cDi{{\mathcal D}\kern-6pt/}
\newcommand\spi{S\kern-6pt/}
\newcommand \cspi{\Sp\kern-6pt/}
\newcommand\quotient{(C^*\Gamma)_{{\rm ab}}/<[1]>}
\newcommand\CC{\mathbb C}

\def \cal {\mathcal}
\def \C {{\cal C}}

\newcommand\NN{\mathbb N}

\newcommand\RR{\mathbb R}
\newcommand\ZZ{\mathbb Z}

\newcommand\pa{\partial}
\newcommand\Ker{\operatorname{Ker}}

\def\Si{{\cal D}^{{\rm sign}}}

{\catcode`@=11\global\let\c@equation=\c@theorem}

\allowdisplaybreaks[2]

\begin{document}
\pagestyle{myheadings}
\markboth{Paolo Piazza and Thomas Schick}{Bordism, rho-invariants and the Baum-Connes conjecture}


\title{Bordism, rho-invariants and the Baum-Connes conjecture}
\author{Paolo Piazza and Thomas Schick}
\maketitle

\begin{abstract}
Let $\Gamma$ be a finitely generated discrete group.  In this
paper we establish vanishing results for rho-invariants associated to
\\
(i) the
  spin-Dirac operator of a spin manifold with positive scalar
  curvature and fundamental group $\Gamma$\\
  (ii)  the signature operator of the disjoint union of a pair of homotopy
  equivalent oriented manifolds with fundamental group $\Gamma$.

  The invariants we consider are
  more precisely
  \begin{itemize}
\item the Atiyah-Patodi-Singer ($\equiv$ APS) rho-invariant associated to a
pair of finite dimensional unitary representations $\lambda_1,
\lambda_2: \Gamma \to U(d)$
\item the $L^2$-rho invariant of Cheeger-Gromov
\item the delocalized eta invariant of Lott for a non-trivial
conjugacy class of $\Gamma$ which is finite. 
\end{itemize}
We prove that all these rho-invariants vanish if the
   group $\Gamma$ is {\it torsion-free} and
the Baum-Connes map
  for the maximal group $C^*$-algebra is bijective. This condition is satisfied, for
  example, by torsion-free
  amenable groups or by torsion-free discrete subgroups of $SO(n,1)$
  and $SU(n,1)$. For the delocalized invariant we only assume the
  validity of
the Baum-Connes conjecture for the {\it reduced} $C^*$-algebra. In
addition to the examples above, this
condition is satisfied e.g.~by Gromov hyperbolic groups
or by  cocompact discrete subgroups of $SL(3,\CC)$.

  In particular, the three  rho-invariants associated to the signature operator
  are,
  for such groups, {\it homotopy invariant}.
For the APS and the Cheeger-Gromov rho-invariants
the latter result had been established by Navin Keswani.
Our proof re-establishes this result and  also extends it to
  the delocalized eta-invariant of Lott. The proof
  exploits in a fundamental way results from bordism theory as well as
  various generalizations
  of the APS-index theorem; it also embeds these
  results in general vanishing phenomena for
  degree zero higher rho invariants (taking values in $A/\overline{[A,A]}$ for
  suitable $C^*$-algebras $A$).
We also obtain precise information about the eta-invariants in
question themselves, which are usually much more saddle objects than
the rho-invariants.

\end{abstract}

\tableofcontents

\section{Introduction and main results}\label{sec:intro}

Throughout this paper, $M$ is an oriented closed manifold of odd
dimension $2n-1$. Let $u\colon M\to B\Gamma$ be a continuous map,
classifying a normal $\Gamma$-covering $\tilde M$ of $M$. Let
$(M^\prime, u^\prime: M^\prime \to B\Gamma)$ be a different
oriented $\Gamma$-covering. We shall say that $(M^\prime,
u^\prime: M^\prime \to B\Gamma)$ and $(M, u: M \to B\Gamma)$ are
oriented $\Gamma$-{\it homotopy equivalent} if there exists an
oriented  homotopy equivalence $h:M\to M^\prime$ such that
$u^\prime \circ h$ is homotopic to $u$.

\subsection{Atiyah-Patodi-Singer rho invariant}
\label{sec:atiyah-patodi-singer}

Let $\lambda_1,\lambda_2\colon \Gamma\to U(d)$ be two finite
dimensional unitary representations of $\Gamma$ of the same
dimension. Equivalently, $\lambda_1-\lambda_2$ is a virtual
representation in the representation ring of $\Gamma$ of dimension
$0$. The representations $\lambda_1, \lambda_2$ induce flat
bundles $L_1$ and $L_2$ on $M$ endowed with flat unitary
connections:  if $V_j$ denotes the representation space of
$\lambda_j$, then $L_j:= \tilde M\times_\Gamma V_j$ and  the
unitary connection  is induced from the trivial connection on the
product $\tilde M\times L_j$.

Let $D$ be a Dirac-type operator acting on the sections of a
Clifford module $E$. We shall always assume that the Clifford module is
 unitary and endowed with a unitary Clifford connection
$\nabla^E$. The operator $D$ can be twisted with the flat bundles
$L_1$ and $L_2$. We use the notation $D_{L_j}$, but also
$D_{\lambda_j}$, for the twisted operator, acting on sections of
$E\tensor L_j$.

The usual integral defines the eta-invariant
of the twisted operator
\begin{equation}\label{eq:define_eta}
  \eta(D_{L_j}) = \frac{1}{\sqrt{\pi}} \int_0^\infty s^{-1/2} \Tr(D_{L_j}\exp(-(s D^2_{L_j}))\;ds.
\end{equation}
Convergence near $s=0$ is ensured by our assumptions on $E$ and
$\nabla^E$, see  \cite{BGV}. In particular, if $M$ is equipped
with a spin structure and $D:= \Di$ is the associated Dirac
operator, then $\eta(\Di_{L_j})$ is defined.
 Similarly, if $D:=D^{{\rm sign}} $ is the signature
operator of an orientation and Riemannian metric on $M$, the
eta-invariant $\eta(D^{{\rm sign}}_{L_j})$ is  defined by
\eqref{eq:define_eta}.

\begin{definition}
\label{def:aps-rho invariant} The Atiyah-Patodi-Singer
rho-invariant associated to the Dirac-type operator $D$ and the
virtual representation $\lambda_1 - \lambda_2$ is the difference
\begin{equation}\label{eq:aps-rho invariant}
\rho(D)_{\lambda_1-\lambda_2} \,:=\,\eta(D_{L_1})-\eta(D_{L_2})
\end{equation}
\end{definition}

The first important result we shall prove in this paper is the following.
\begin{theorem}\label{theo:APS_rho}
  Assume that $\Gamma$ is
  a torsion-free discrete group and that the Baum-Connes
  assembly map is an isomorphism
  \begin{equation*}
  \mu_{max}\colon  K_*(B\Gamma) \to
  K^*(C^*\Gamma),
  \end{equation*}
  where $C^*\Gamma$ is the maximal $C^*$-algebra of $\Gamma$.
\begin{enumerate}
\item  Let $D$ be the Dirac operator of a spin structure, $D=\Di$. If the
metric on
  $M$ has positive scalar curvature then the Atiyah-Patodi-Singer
  rho-invariant  is zero.
\item    Let $D$ be the signature operator, $D=D^{{\rm sign}}$,  then
the Atiyah-Patodi-Singer rho invariant
      $\rho (D^{{\rm sign}})_{\lambda_1-\lambda_2}$ only
    depends  on the oriented $\Gamma$-homotopy type of $(M,u:M\to B\Gamma)$.
\end{enumerate}
  \end{theorem}

\begin{remark}
  The result about the APS-$\rho$-invariant of the {\it signature operator} is due
  to Navin Keswani \cite{Kes1}. We give an independent proof of
  Keswani's result.
 The first special case of this result is due to Mathai
  \cite{MR93j:58125}, who proved it for Bieberbach groups.
 We shall also try to present a ``philosophical'' reason for
  the theorem. It should be added Keswani's argument is most likely to
be adapted so as to cover the result on the spin Dirac
operator as well. Finally, Shmuel Weinberger  \cite{Wei}
proves the homotopy invariance of the
  APS-rho invariants for the signature operator
  under the assumption
 that $\Gamma$ is torsion
  free and that the $L$-theory
  isomorphism conjecture holds for $\Gamma$. Chang \cite{Chang}  has used the same
  assumption to prove  the corresponding result for the  $L^2$-rho
  invariant .
\end{remark}

  \begin{remark}
    Note that for \emph{finite} groups the $\rho$-invariants of the spin
    Dirac operator are effectively used to distinguish different
    metrics with positive scalar curvature, compare
    e.g.~\cite{BoGi}. We see that this is not possible under
    our assumption on the group $\Gamma$. This might be particularly
    surprising in view of the genuinely non-local definition of the
    eta-invariants.
\end{remark}

We shall use different $C^*$-algebras associated to the group $\Gamma$.

  \begin{definition}\label{def:group_Cstar_algebras}
   Let $\Gamma$ be a discrete group.
    The \emph{maximal $C^*$-algebra $C^*\Gamma$} is the completion of
    the group ring $\complexs[\Gamma]$ with
    respect to the maximal possible $C^*$-norm on this ring
    \cite{MR98a:46066}.

    The maximal $C^*$-algebra has the universal property that for each
    unitary representation $\phi\colon \Gamma\to U(H)$ on a Hilbert
    space $H$, there is a unique $C^*$-algebra homomorphism
    $C^*\Gamma\to \boundedops(H)$ extending $\phi$.

    More frequently used is the \emph{reduced $C^*$-algebra} of
    $\Gamma$. It is by definition the closure of $\complexs[\Gamma]$,
    considered as subalgebra of $\boundedops(l^2(\Gamma))$ with
    respect to the right regular representation ($
    \left(\sum_{\gamma\in\Gamma} \lambda_\gamma \gamma\right) \cdot g:=
    \sum_{\gamma\in\Gamma} \lambda_\gamma (\gamma g)$). Note that by
    the universal property of $C^*\Gamma$, we have a canonical
    homomorphism $C^*\Gamma\to C^*_{red}\Gamma$ extending the identity
    on $\Gamma$.

    We will also make use of the \emph{group von Neumann algebra}
    $\NeumannN\Gamma$. This is by definition the weak closure of the
    right regular representation of
    $\complexs[\Gamma]$ in $\boundedops(l^2\Gamma)$. In particular it
    is also the weak closure of $C^*_{red}\Gamma$, i.e.~we have an
    inclusion $C^*_{red}\Gamma\into \NeumannN\Gamma$. By the
    bicommutant theorem, it is also equal to
    $\boundedops(l^2(\Gamma))^\Gamma$, i.e.~all operators which
    commute with the right regular representation of $\Gamma$ on
    $l^2(\Gamma)$. Note that $\NeumannN\Gamma$ is also a $C^*$-algebra.
\end{definition}

\begin{remark}
  The Baum-Connes  conjecture for a {\it torsion-free} group $\Gamma$ asserts
  that the assembly map
  \begin{equation*}
    K_*(B\Gamma)\to K^*(C^*_{red}\Gamma)
  \end{equation*}
  is an isomorphism, where the right hand side is the \emph{reduced}
  $C^*$-algebra of the group $\Gamma$.
  It is a fact that $K_*(C^*\Gamma)$ and $K_*(C^*_{red}\Gamma)$ are
  identical for large classes of groups, e.g.~amenable groups or
  discrete subgroups of $SO(n,1)$ or $SU(n,1)$. However, for groups with Kazdan's
  property T, they definitely differ, and $K_*(B\Gamma)\to
  K^*(C^*\Gamma)$ is known not to be surjective in this case.

\end{remark}



\subsection{$L^2$-rho invariant}
\label{sec:l2-rho-invariant}

If we look at the $\Gamma$-covering $\tilde M$ directly, we can
lift any \emph{differential} operator $D$ from $M$  to a
$\Gamma$-invariant differential operator $\tilde D$ on $\tilde M$.
In particular,  this is possible for a Dirac-type operator $D$.
Moreover, by using  Schwartz kernels, but integrating only over a
fundamental domain for the covering $\tilde M\to M$, we get an
$L^2$-trace $\Tr_{(2)}$, see \cite{Atiyah-covering}. This yields
the $L^2$-$\eta$-invariant defined using essentially
formula \eqref{eq:define_eta} by
\begin{equation*}
  \eta_{(2)}(\tilde D) = \frac{1}{\sqrt{\pi}} \int_0^\infty s^{-1/2}
  \Tr_{(2)}(\tilde D \exp(-(s \tilde D ^2 ))\;ds,
\end{equation*}
where the $L^2$-trace is defined by
\begin{equation}
  \Tr_{(2)}(\tilde D\exp(-t \tilde D ^2)) = \int_{\mathcal{F}} \tr_x k_t(x,x).
\end{equation}
Here $\tr_x$ is the fiberwise trace and $\mathcal{F}$ is a
fundamental domain.

In order to ensure the convergence of the integral near $t=0$ we
still require the operator $D$ to be associated to a unitary
Clifford module endowed with a unitary Clifford connection. The
convergence of the integral for large $t$ is discussed in
\cite{ChGr} and also in \cite{Ra}.

\begin{definition}
  We define the $L^2$-$\rho$-invariant as the difference
  \begin{equation*}
    \rho_{(2)}(D) := \eta_{(2)}(\tilde D) - \eta(D)
  \end{equation*}
\end{definition}

\begin{theorem}\label{theo:L2_rho}
  If $\Gamma$ is torsion free and the Baum-Connes assembly map for the
  maximal $C^*$-algebra of
  $\Gamma$, $K_*(B\Gamma)\to K_*(C^*\Gamma)$, is an isomorphism, then
  \begin{enumerate}
 \item the $L^2$-rho invariant of the Dirac operator of a spin
    manifold with positive scalar curvature vanishes.
\item  the $L^2$-rho invariant of the signature operator
    depends only on the oriented $\Gamma$-homotopy type of $(M,u:M\to B\Gamma)$.
  \end{enumerate}
\end{theorem}

\begin{remark}
  For the signature operator, this result is proved by Keswani in
  \cite{Kes2} using methods different from ours. 

  The result for the spin Dirac operator was originally obtained (with
  a similar, but slightly more complicated method as the one presented
  here) in unpublished work of Nigel Higson and the second author.
It should be possible to adapt Keswani's arguments so as to cover  the result on
the spin Dirac operator.

Our proof also provides some rather delicate information about the eta-invariants directly (and not
just of the less saddle rho-invariants).
\end{remark}

\subsection{Delocalized eta invariants}
\label{sec:deloc-rho-invar}

We continue with the geometric set-up of Section
\ref{sec:l2-rho-invariant}, with a Galois $\Gamma$-covering
$\Gamma\longrightarrow \tilde M \longrightarrow M$ classified by
a map $u:M\to B\Gamma$. To construct the $L^2$-eta invariant, we
used the integral kernel $k_t(x,y)$ of $\tilde D\exp(-t\tilde
D^2)$ on $\tilde M$.  Fix now a non-trivial conjugacy class $<g>$ of
$\Gamma$. Define the \emph{delocalized} trace
\begin{equation}\label{eq:def_delocalized_trace}
  \Tr_{<g>}(\tilde D \exp(-t\tilde D^2)) := \sum_{g\in<g>} \int_{\mathcal{F}} \tr_x k_t(x,gx).
  \end{equation}
 Note that the fibers at $x$ and $gx$ of the pull back vector bundle
 on which $\tilde D$ acts are canonically identified, so that
 $k_t(x,gx)$ can be considered as an endomorphism of this fiber, and
 its fiberwise trace $tr$ is defined.

Convergence of this (possibly infinite) sum follows from exponential
decay, compare \cite{MR2000k:58039}.

\medskip
\noindent
{\bf Finite conjugacy classes}
If the conjugacy class $<g>$ contains only finitely many elements,
Lott proves in \cite{MR2000k:58039} that the following integral,
defining the \emph{delocalized eta invariant}, converges.
\begin{equation}\label{eq:delocalized_eta}
  \eta_{<g>}(\tilde D):= \frac{1}{\sqrt{\pi}} \int_0^\infty t^{-1/2}
  \Tr_{<g>}(\tilde D\exp(-t\tilde D^2)\,dt.
\end{equation}

We prove:
\begin{theorem}\label{theo:finite_conjugacy}
  Assume that $\Gamma$ is torsion-free and that the assembly map
  $\mu_{red}\colon K_*(B\Gamma)\to
  K_*(C^*_{red} \Gamma)$ for the reduced group $C^*$-algebra is an
isomorphism. Let $<g>$ be a
  non-trivial finite conjugacy class in $\Gamma$.
\begin{enumerate}
\item   If $D=\Di$ is the spin-Dirac operator of a spin manifold with positive scalar
  curvature  then
  \begin{equation*}
    \eta_{<g>}(\tilde \Di) =0.
  \end{equation*}
\item
If $D=D^{{\rm sign}}$ is the signature operator of an oriented Riemannian
   manifold  then $ \eta_{<g>}(\tilde{D}^{{\rm sign}})$ depends only on the
   oriented $\Gamma$-homotopy type of $(M,u:M\to B\Gamma)$.
   \end{enumerate}
\end{theorem}

\begin{remark}
Notice that our assumption here involves the {\it reduced} group
$C^*$-algebra. There are substantially more groups which satisfy the
Baum-Connes conjecture for the reduced $C^*$-algebra, even some with
property T, e.g.~all  Gromov hyperbolic groups \cite{Mineyev-Yu} or cocompact
discrete subgroups of $Sl(3,\complexs)$ \cite{Lafforgue}.
\end{remark}

\begin{remark}
The stated result is only interesting if the group contains elements with
finite conjugacy class. Not that this does by no means imply that the
corresponding element has finite order.

Evidently, each central element has a finite conjugacy class
(consisting only of itself). This means that, starting with an
arbitrary group $\Gamma$, each central extension with kernel
$\integers$ will be a group with large center, and if
$H^2(\Gamma,\integers)\ne \{0\}$, there are non-trivial such central extensions.
\end{remark}

\noindent
{\bf Infinite conjugacy classes.}
Lott establishes the convergence of the delocalized eta invariant
under more general assumptions: it suffices to
assume that  the Dirac-type operator on the covering has a gap
at 0 (with 0  allowed in the spectrum)
and that the conjugacy class is of polynomial growth with respect
to a word-metric on $\Gamma$. Thus, for the spin Dirac operator
of a manifold with positive scalar curvature, the delocalized eta invariant $\eta_{<g>}(\tilde\Di)$
is well defined, provided $<g>$ is of polynomial growth. If,      in addition,
$\Gamma$ is torsion-free and the reduced Baum-Connes map is bijective,
then we shall prove that, as in the case of finite conjugacy classes stated above,
$\eta_{<g>}(\tilde\Di)=0$.

\subsection{General principle of the proofs of vanishing results}
\label{sec:gener-princ-proofs}

In order to simplify the exposition we shall concentrate on the
case where $\Gamma$ is the fundamental group of our manifold and
the $\Gamma$-covering is the universal covering.

 First observe that the {\it homotopy invariance} of the rho
invariants for the signature operator reduces to the {\it
vanishing} of the rho invariant for the disjoint union
$X\disjointunion - X'$, if $X$ and $X'$ are homotopy equivalent,
since all rho invariants are certainly additive under disjoint
union.

For this reason we shall be  only concerned with vanishing results. To establish
these results , we apply the following {\it general principle}. To
avoid undo repetitions, let $\rho$ for the moment stand for any of
the $\rho$-invariants we want to investigate.
\begin{enumerate}
\item We first define a \emph{stable} variant $\rho^s$ of $\rho$. This
  will be defined as the invariant of a perturbation
  of our generalized Dirac operator. Such perturbations do not always
  exist, we need the vanishing of the index class of the
  generalized Dirac operator. This very strong assumption is satisfied
  for  geometric reasons if one looks at the Dirac operator of a spin
  manifold with positive scalar curvature, as well as for the
  signature operator on the disjoint union $X\disjointunion -X'$ of
  two homotopy equivalent manifolds.

  We study the main properties of $\rho^s$. Most
  important is that it appears as the correction term in an index theorem
  for manifolds with boundary for suitably perturbed Dirac
  operators.
We use this fact and the assumed {\it  surjectivity} of the Baum-Connes
map in order to show that the stable rho-invariant is {\it well
defined}, independent of the chosen perturbation (we are always
under the assumption that the index class of our
operator is zero). Under these assumptions we also prove the
fundamental fact that
$\rho^s$ is a {\it bordism invariant}:  suitable
  index theorems on manifolds with boundary will again play a crucial role 
here.

\item  Then we use our injectivity assumption on the Baum-Connes map,
fundamental results in bordism theory and the bordism invariance
 established in (1)   in order to show that
  the stable invariant $\rho^s$, whenever it is defined, is equal to the stable
  invariant of a particularly nice manifold. 
  For this nice manifold
  we compute the stable invariant and show that {\it it vanishes}.

  To this point, we have therefore shown that in certain special
  situations one can define an invariant $\rho^s$ which turns out to be zero.
\item As a last step we show that in the two geometric situations we
  are studying, the stable invariant $\rho^s$ coincides with the unstable
  invariant $\rho$. This will be done by constructing very special
  perturbations (used in the definition of the stable invariant) which
  make the direct comparison of the stable and unstable invariant
  possible.   For the signature operator on $X\disjointunion -X'$ we use
  perturbations that are inspired by the
  work of Hilsum-Skandalis

  In fact, our results here are much more precise: they give
  information directly about the unstable $\eta$-invariants. This is
 quite remarkable because of the non-stable nature of the
  $\eta$-invariants under perturbations and might lead to future
  applications or calculations of $\eta$-invariants.

  In the case of a spin manifold with positive scalar
  curvature, we don't have to perturb at all, so the last step is
  trivial.

\end{enumerate}

\subsection{Examples}
\label{sec:non-vanish-results}

Let $\Gamma$ be a torsion-free discrete group satisfying our basic assumption,
the bijectivity of the Baum-Connes map. One might very well wonder whether
the signature rho-invariants considered in this paper are non-zero
and whether they can be effectively used in order to distinguish
manifolds that are not homotopy equivalent.
In the last part of the paper we give a careful treatment
of some
non-trivial examples.
In particular, we construct manifolds with the same cohomology
but that are not in the same homotopy class and we distinguish their
homotopy type through their rho-invariants. This can be done for all non-trivial
groups which satisfy our basic examples. The results might also be
obtainable using Blanchfield forms of classical algebraic topology;
however, it seems that one has to use some advanced version to cover
all the cases covered by our invariant. The advantage of our approach
is that it is very easy to carry out the calculations.

 We also show, along the way, that the
vanishing of the signature index class
is not sufficient for establishing the vanishing of the
signature-rho-invariants.
The mere vanishing of the signature index class
does not imply the vanishing of our rho-invariants.

\subsection{Plan of the paper.}\label{subsec:plan}

In {\bf Section \ref{sec:degree-zero-as}}
we gather the index theoretic results that will be needed
throughout the paper. We treat the general case
of a Dirac operator twisted by a bundle $\mathcal{L}$ of finitely generated
projective $A$-modules, with $A$ a $C^*$-algebra.
We state $A/\overline{\Commutator{A}{A}}$-valued index theorems on
closed manifolds and on manifolds
with boundary.
Proofs are given in {\bf Appendix \ref{sec:proof}}.

In {\bf Section \ref{sec:discrete groups}} we specialize to the $C^*$-algebras defined
by a discrete group $\Gamma$ and to the corresponding Mishchenko-Fomenko
bundle; we also discuss the corresponding notions when we
twist with the group von Neumann algebra $\mathcal{N}\Gamma$.

In {\bf Section \ref{sec:surjective+stable}} we use the {\it surjectivity} of the assembly
map in order to define the {\it stable} rho-invariants and establish
their bordism invariance. The stable rho-invariants are $C^*$-algebraic objects
and they are defined under the additional assumption that the index class
in $K_1 (C^* \Gamma)$ is equal to zero.

In {\bf Section \ref{sec:injectivity}} we employ the {\it injectivity} of the assembly
map in order to construct a suitable bordism between $d$ copies
of the manifold $X=M\cup (-M^\prime)$, with $M$ and $M^\prime$ homotopy equivalent,
and a manifold of a special type
A similar result is proved if
$X$ is a manifold with positive scalar curvature.
These bordisms take the classifying maps
into account.

In {\bf Section \ref{sec:vanishing-stable-rho}}
we show that the stable rho-invariants of these
special manifolds are zero. This fact and the bordism invariance
of Section \ref{sec:surjective+stable}
are then  used in {\bf Section \ref{sec:vanishing-stable-rho-1}} in order to prove that under
our assumptions the stable rho-invariants for $D=D^{{\rm sign}}$ are zero if $X=M\cup (-M^\prime)$,
with $M$ and $M^\prime$ homotopy equivalent.
The corresponding result for  $X$  a manifold with positive scalar
curvature and $D=\Di$ is established in {\bf Section \ref{sec:vanishing-spin}}.

In {\bf Section \ref{sec:unst-rho-invar}} we introduce {\it unstable} rho-invariants; these
are von Neumann objects producing the three rho-invariants defined in Subsections
\ref{sec:atiyah-patodi-singer},
\ref{sec:l2-rho-invariant}, \ref{sec:deloc-rho-invar}
as special cases. The unstable rho-invariants are always defined.

In {\bf Section \ref{sec:homot-invar-sign}}, {\bf Section
\ref{sec:stable-=-unstable}}  and {\bf Section
\ref{sec:homot-invar-unst}} we show that if $X=M\cup (-M^\prime)$,
with $M$ and $M^\prime$ homotopy equivalent, then the unstable
rho-invariants are suitable limits of stable
rho-invariants. 
This
step completes the proof of our theorems for the signature
operator, since we know that stable rho-invariants are zero under
our assumptions. The proof for the spin Dirac operator, which is
much easier since there is no need for this limit-argument,
 is given in Section \ref{sec:vanishing-spin}


In {\bf Section \ref{sec:remarks}} we gather several remarks on
delocalized eta-invariant for
infinite conjugacy classes. We also prove vanishing results for {\it
  higher}
rho-invariants.

In {\bf Section \ref{sec:center-valued-l2}}
we state a general von Neumann signature formula for
manifolds with boundary.

In {\bf Section \ref{sec:examples-non-trivial}}
we give examples of closed manifolds with torsion free fundamental
group satisfying the Baum-Connes assumption and for which the relevant
rho-invariants are non-zero. We  also   construct manifolds with isomorphic
homology but with different homotopy type  and we distinguish them
by using the $L^2$-rho invariant or the delocalized eta invariant.

In four additional Appendices we recall signatures and the signature
operator of Hilbert $A$-module chain complexes with symmetry, and we give a detailed account of the
relationship
between spectral invariants on  $\Gamma$-coverings and certain algebra-valued
invariants; we also discuss  naturality properties of these
algebra-valued invariants.
All appendices either recall known results or immediate extensions of
known results; they are included for the reader's convenience.

\subsection{Acknowledgments}
\label{sec:acknowledgements}

We thank Paul Kirk, Eric Leichtnam, Victor Nistor, George
Skandalis and Shmuel Weinberger for helpful discussions and
remarks. {Special thanks go to Nigel Higson and John Roe for pointing out a
 gap  in an
  earlier version of the paper.}
Part of this work was carried out during visits of the authors to
G{\"o}ttingen, Paris and  Rome funded by Ministero Istruzione
Universit{\`a} Ricerca, Italy (Cofin {\it Spazi di Moduli e Teorie
di Lie}), Institut de Math{\'e}matiques de Jussieu , CNRS,
Graduiertenkolleg ``Gruppen und Geometrie'' (G{\"o}ttingen).


\section{Index theory: statement of results}
\label{sec:degree-zero-as}

 Many of
the results of this paper are based on suitable generalizations of
the Atiyah-Singer and Atiyah-Patodi-Singer index theorem. In this
section we gather the index theoretic results that we shall need
in the rest of the paper.
Most of these index theorems  are
described (sometimes implicitly)  in the literature; they are due to Mishchenko-Fomenko
in the closed case and Leichtnam and the first author in the
boundary case.
Nevertheless, we
shall give a direct and  simplified account
of these results, showing in particular that once the index class
in $K_* (C^*_r \Gamma)$ is carefully described, it is possible
to give heat-kernel-proofs of all needed results following
the same steps as for the numeric case (see, for example,
\cite[Introduction]{Melrose}).

This section is devoted to the
statements; Appendix \ref{sec:proof}  will contain the proofs. We start with the
closed case.

\subsection{Atiyah-Singer index theory for arbitrary $C^*$-algebras}
\label{sec:as-index-theory}

\subsubsection{Geometric set-up.}\label{set-up:A_bundle_without_boundary}
  Let $M$ be a closed manifold. We assume that $D$ is a Dirac-type
  operator acting on a
  finite dimensional  Clifford module $E$ on $M$,
see \cite{BGV} for the definitions. Our main interest
  will be the
  spin-Dirac operator $\Di$ of a spin structure, and the signature operator
  $D^{{\rm sign}}$ of an orientation and Riemannian metric on $M$.
Let $A$ be a $C^*$-algebra, and $\L$ a bundle of finitely
generated
  projective Hilbert $A$-modules\footnote{{\bf A comment on notations:}
   objects in the $C^*$-algebraic context will always be
 denoted by calligraphic  letters}, with an $A$-connection
  $\nabla_{\L}$. We define $D_{\L}$
to be the operator $D$ twisted with the
  bundle with connection $\L$, acting on sections of $E\tensor \L$.
  Then $D_{\L}$ is an $A$-linear differential operator,
  $$D_{\L} \in {\rm Diff}^1_{A}(M;E\otimes \L,E\otimes \L)\subset \Psi^*_{A}(M;E\otimes \L,E\otimes
  \L).$$
  On the right hand side the  Mishchenko-Fomenko pseudodifferential
  calculus appears.
$D_{\L}$ is elliptic in the sense of Mishchenko-Fomenko with a
well defined index class $\Ind (D_{\L})\in K_{\dim (M)} (A)$. We
shall recall below the definition of $\Ind (D_{\L})\in K_{\dim
(M)} (A)$ when $\dim(M)$ is even.

\medskip

In many situations, instead of working with the differential
operator $D_{\L}$ itself, we will have to perturb the operator
slightly (thereby leaving the world of differential operators).
We do this because on the one hand  we will have to improve some technical properties
of the operator, such as the large time behavior of the
associated heat-kernel; on the other hand, when the operator
arises as a boundary operator, we shall need invertibility in
order to define a Fredholm problem on the manifold with boundary.
Of course, each time such a perturbation is introduced,
 we shall have then
 to control how
invariants such as the index or the eta-invariant behave.

 We  start
with some generalities  about perturbations; this material will be
needed throughout the paper.

\begin{definition}
Let $\C\in \Psi^{-\infty}_{A}$ be a smoothing operator in the
  Mishchenko-Fomenko calculus, acting on sections of $E\tensor \L$.
Set
  \begin{equation*}
  D_{\L,\C}:= D_{\L}+\C.
\end{equation*}
  This is a (smoothing) \emph{perturbation} of $D_{\L}$.
\end{definition}

 The following lemma is a direct consequence of the
Mishchenko-Fomenko ($\equiv$ MF) pseudodifferential calculus:

\begin{lemma}
 $D_{\L,\C}$ is
an elliptic element in the MF-calculus
\begin{equation}
  D_{\L,\C}\in  \Psi^1_A(M; E\otimes \L,E\otimes \L)\subset
  \Psi^*_{A}(M;E\otimes \L,E\otimes \L)
\end{equation}
Moreover, the operators $D_{\L}$ and $D_{\L,\C}$ have the same
index $$\Ind(D_{\L})=\Ind(D_{\L,\C})\in K_{\dim(M)}(A)\,.$$
\end{lemma}

Fundamental in what follows is the following Mishchenko-Fomenko
decomposition theorem.

\begin{theorem}\label{lem:Mishchenko_decomposition_closed}
 Let $M$ be  even dimensional, so that $E=E^+\oplus E^-$.
 Let $(E\otimes \L)^\pm := E^\pm \otimes \L$.
  There is a Mishchenko-Fomenko decomposition of
the space of sections
  of $E\tensor \L$ with respect to $D_{\L}$, i.e.
  \begin{equation}\label{Mishchenko_decomposition_closed}
    C^\infty(M,(E\tensor \L)^+) = \mathcal{I}_+ \oplus
    \mathcal{I}_+^{\perp},\quad C^\infty(M,(E\tensor
    \L)^-)=\mathcal{I}_- \oplus D_{\L}(\mathcal{I}_+^{\perp}).
  \end{equation}
 By completion, we obtain a decomposition of the Sobolev Hilbert
  $A$-modules $H^m(M,(E\tensor \L)^{\pm})$ for any $m\in \NN$.

  The second decomposition is not, a priori, orthogonal.
  However, $D_{\L}$ induces an isomorphism (in the Fr{\'e}chet
  topology) between $\mathcal{I}_+^{\perp}$ and
  $D_{\L}(\mathcal{I}_+^{\perp})$, and $\mathcal{I}_+$ and $\mathcal{I}_-$
  are finitely generated projective Hilbert $A$-modules (consisting of
  smooth section of $L^2(M,E\tensor \L)$). The projections
  $\Pi_{\mathcal{I}_+}$ onto
  $\mathcal{I}_+$ (orthogonal) and $\Pi_{\mathcal{I}_-}$ onto $\mathcal{I}_-$ (along
  $D_{\L}(\mathcal{I}_+^{\perp})$) are smoothing operators in the
  Mishchenko-Fomenko calculus. Because $\mathcal{I}_{\pm}$ are
  already complete finitely generated projective Hilbert $A$-modules
  they are unchanged under the completions.
The Mishchenko-Fomenko index class is given
by
\begin{equation}\label{index class closed}
\Ind (D_L)= [\mathcal{I}_+] - [\mathcal{I}_-]\;\in\;K_0 (A)\,.
\end{equation}
\end{theorem}

\begin{remark}
The theorem is  ultimately a consequence of the ellipticity of
$D_{\L}$ and more precisely of the existence of an inverse for
$D_{\L}$ modulo smoothing operators. Thus, an identical statement
remains true for the more
    general (elliptic) operator $D_{\L,\C}$.
    The existence of the decomposition  is implicitly proved in
    \cite{Mish-Fom}; the structure of the two projections is
    analyzed in \cite[Appendix A]{LPMEMOIRS}.
    \end{remark}

\subsubsection{The $A/\overline{\Commutator{A}{A}}$-valued index}
\label{sec:perturbations}

We want to give a heat kernel proof of a suitable Atiyah-Singer
($\equiv$AS) index theorem for the operator $D_{\L}$ of
\ref{set-up:A_bundle_without_boundary}. First of all, we introduce
the index we want to compute. Consider the index class $\Ind
(D_{\L})$ expressed through the MF-decomposition theorem (see
\ref{index class closed}):
 $\Ind
(D_{\L})= [\mathcal{I}_+] - [\mathcal{I}_-]\;\in\;K_0 (A)\,.$
Expressing the finitely generated projective modules
$\mathcal{I}_\pm$ as the images  of idempotent matrices and taking
the difference of the traces, we obtain a well defined  element in
$A/\overline{\Commutator{A}{A}}$, with $\Commutator{A}{A}$ equal to the $\complexs$-subspace
generated by the
  commutators $\Commutator{a}{b}:=ab-ba$
  for $a,b\in A$.

\begin{definition}\label{def:of_commutator}
Consider the closure $\overline{\Commutator{A}{A}}$ of  the subspace
$\Commutator{A}{A}$.
   We set
\begin{equation}\label{definition-of-ab}
A_{{\rm ab}}:= A/\overline{\Commutator{A}{A}}.
\end{equation}
where the subscript {\it ab} stands for {\it abelianization}.
It should be noticed that $A_{{\rm ab}}$ is a commutative Banach algebra.
\end{definition}

In this way we define a homomorphism of abelian
  groups
\begin{equation}\label{algebraic-trace}
\tr^{{\rm alg}}:K_0 (A)\longrightarrow A_{{\rm ab}}
\end{equation} Since this is nothing but the zero-degree part of
the Karoubi-Chern character \cite{Karoubi1},  we  denote  $\tr^{{\rm
alg}}([\mathcal{I}_+] - [\mathcal{I}_-])$ by
\begin{equation}\label{def:algebraic_trace}
[\mathcal{I}_+]_{[0]} -
[\mathcal{I}_-]_{[0]}\;\in\;A_{{\rm ab}}\,,
\end{equation}
or, equivalently, by $\Ind_{[0]}
(D_{\L})\;\in\;A_{{\rm ab}}$.

Our aim is to give a formula for $\Ind_{[0]}
(D_{\L})\;\in\;A_{{\rm ab}}$ and to  prove it via heat
kernel techniques. First we need the existence of the heat-kernel
for the Dirac Laplacian $D_{\L}^2$ and its perturbation
$D_{\L,C}^2$. This is not completely obvious, given that we are in
the $C^*$-algebraic context. We shall give the precise statement
in Lemma \ref{lem:heat semigroup closed} in Appendix \ref{subsec:proof-closed}
For the time being we shall content ourselves with the statement
that the heat-semigroup $\{e^{-tD_{\L,C}^2}, t>0\}$ exists,
provides a fundamental solution of the heat-equation
and is such that $e^{-tD_{\L,C}^2}\, \in \Psi^{-\infty}_A$.

Since the heat operator
 is a smoothing operator in the
  Mishchenko-Fomenko calculus, we can consider its supertrace
 $$\STR(e^{-tD_{\L}^2})\,\in\, A_{{\rm ab}} \,.$$
This is defined as follows:
we restrict the heat-kernel to the
  diagonal $\Delta$ in $M\times M$, $\Delta\leftrightarrow M$, thus obtaining
an element in $C^\infty (M, {\rm End}(E\otimes \L))$, we take
 for  each fiber over $x$  the $A_{{\rm ab}} $-valued
  supertrace on  ${\rm End} (E\tensor \L)_x$ and integrate
\begin{equation}
\STR(e^{-tD_{\L} ^2}) :=\int_M {\rm str}^{\rm alg}_x (e^{-tD_{\L} ^2}(x,x))
{\rm vol}_M\,\in\,A_{{\rm ab}}\,.
\end{equation}

Our first result concerns the large time behavior of the heat-kernel
in the Mishchenko-Fomenko calculus.

\begin{proposition}\label{prop:large time closed no alpha}
  \begin{equation*}
  \lim_{t\to\infty}\STR (e^{-tD^2_{\L}}) = [I_+]_{[0]} - [I_-]_{[0]}\equiv
  \Ind_{[0]} (D_{\L}) \in A_{{\rm ab}} ,
\end{equation*}
where part of the assertion is that the limit exists.
\end{proposition}

The proof of this proposition can be found in Appendix
\ref{subsec:proof-closed}. Next we want to connect $\Ind_{[0]}
(D_{\L})\in A_{\rm ab}$ defined using the index class and the
algebraic trace  $\tr^{{\rm alg}}:K_0 (A) \to A_{\rm ab}$ to the
integral-kernel-trace, $\TR$, of the projection operators
$P_+\,:=\,  \Pi_{\mathcal{I}_+}$, $P_- \,:=\, \Pi_{\mathcal{I}_-}$
given in the Mishchenko-Fomenko decomposition. We recall that
these are smoothing operators; the trace $\TR$ is thus well
defined. We state the result here and defer the proof to Appendix
\ref{subsec:proof-closed}.

\begin{proposition}\label{prop:alg=int}
  The algebraic trace \begin{equation*}
\Ind_{[0]}(D_{\L}):=
  [I_+]_{[0]}-[I_-]_{[0]}\;\;\;\text{of}\;\;\; Ind (D_{\L}),
\end{equation*}
  i.e.~the image under the map
$\tr^{{\rm alg}}\colon K_0(A)\to A_{\rm{ab}},$
of $\Ind (D_L)$,can be calculated as
  \begin{align*}
  [I_+]_{[0]}-[I_-]_{[0]} &= \TR (P_+) - \TR (P_-)\\&\equiv  \int_M
  \tr^{\rm alg}_x P_+(x,x) - \int_M
  \tr_x^{\rm alg}
  P_-(x,x) \in
  A_{\rm{ab}},
  \end{align*}
  where $P_+$ and $P_-$ are the projections onto $\mathcal{I}_+$ and $\mathcal{I}_-$ as
  given by the Mishchenko-Fomenko decomposition.
\end{proposition}

Proceeding now as in the classical case, one proves the following
$A_{\rm{ab}}$-valued Atiyah-Singer index theorem:

\begin{theorem}\label{theo:general_AS}
  We have
  \begin{equation*}
      \Ind_{[0]} (D_{\L})      =
      \TR(P_+)-\TR(P_-)
      =
      \int_M AS(D)(x)\wedge \ch{\L}(x)_{[\dim M]} \in
      A_{\rm{ab}}
  \end{equation*}
  where $AS(D)(x)$ is the local integrand in the Atiyah-Singer
  formula for $D$.
\end{theorem}

In the above  formula the differential form
 \begin{equation*}
(x\mapsto AS(D)(x)\wedge \ch(E)(x)\wedge \ch{\L}(x)_{[\dim
    M]}) \in \Omega^{\dim M} (M,A_{ab}),
\end{equation*}
can be calculated  as usual using Chern-Weyl theory and the curvature of the
  connections, see \cite{Schick03}.

\subsection{APS index theory for arbitrary $C^*$-algebras}
\label{sec:aps-index-theory}


\subsubsection{Geometric set-up.}\label{sit:A_bundle_with_boundary}

Let $W$ a compact manifold of even dimension with
  boundary  $\boundary W=M$. We assume that $D$ is a Dirac-type
  operator acting on a
  finite dimensional Dirac-bundle $E$ on $W$. Our main interest
  will be the
  spin-Dirac operator of a spin structure, and the signature operator
  of an orientation and Riemannian metric on $W$. We always assume
  that all these structures are of product type near the boundary.

  Let $A$ be a $C^*$-algebra, and $\L$ a bundle of finitely generated
  projective Hilbert $A$-modules, with an $A$-connection
  $\nabla_{\L}$ (again everything has a product structure near the
  boundary). We define $D_{\L}$ to be the operator $D$ twisted with the
  bundle with connection $\L$. Associated to this are the boundary
  operators $D_M\equiv D_{\partial W}$ of $D$ and $D_{M,\L}\equiv D_{\partial W,\L}$ of $D_{\L}$.
We can attach an infinite cylinder $(-\infty,0]\times M$ to $W$
along its boundary $\partial W=M$, thus obtaining a manifold with
cylindrical ends $\widehat{W}$ and with product metric $dt^2 +
g_{\partial W}$ along the cylinder. The operators $D_{\L}$ extends
in a natural way to the manifold $\widehat{W}$. The change of
coordinates $t=\log x$ compactifies $\widehat{W}$ to a manifold
with boundary and with product $b$-metric $dx^2/x^2 + g_{\partial
M}$ near the boundary. The operator on $\widehat{W}$ then defines
in a natural way a $b$-differential operator on the compactified
manifold. We refer the reader to the book \cite{Melrose} of Melrose
for basics about the $b$-calculus. As we shall
mainly work in the framework of the $b$-calculus, we keep denoting
the compactified manifold by $W$ and the resulting
$b$-differential operators by $D_{\L}$.

\subsubsection{Trivializing perturbations}\label{subsec:trivializing}

Let  $M$ be odd  dimensional and without boundary, and let us
assume that $\Ind(D_{M,\L})=0\in K_{1}(A)$. This will be the case
in the following examples:
\begin{itemize}
\item  $D_M$ is the signature operator on $M=X\sqcup (-X^\prime)$ with
  $X$ and $X^\prime$ homotopy equivalent
\item  $D_M$ is the signature
  operator and $(M,\L)$ bounds,  i.e.~there is a manifold $W$ with
  $\boundary W=M$ such that $\L$ extends to $W$
\item
 $M$ is spin with positive scalar curvature and $D_M$ is the spin Dirac
  operator.
\end{itemize}
For more on these  examples we refer the reader to
\cite{HiSka}, \cite{RosI} \cite{LPGAFA} and \cite{RosIII}. According to
 \cite{WuI} \cite[Theorem 3]{LP03},
  we can find a
non-commutative spectral section $\mathcal{P}$ for
$D_{M,\L}$.\footnote{ Thus, by definition, $\mathcal{P}$ is a
  self-adjoint
projection, $\mathcal{P}\in \Psi^0_{A}$ and there exists
functions $\chi_1, \chi_2 \in C^\infty (\RR,[0,1])$ such that
$\chi_i (t)=0$ for $t<<0$, $\chi_i (t)=1$ for $t>>0$,
$\chi_2\equiv 1$ on a neighborhood of the support of $\chi_1$ and
$\Im \chi_1 ({D}_{M,\L})\subset \Im\mathcal{P}\subset \Im \chi_2
({D}_{M,\L}).$} Using the projection $\mathcal{P}$ one can construct
\cite[Proposition 2.10]{LPGAFA}
a smoothing
operator $\mathcal{C}_{\mathcal{P}}\in \Psi^{-\infty}_{A}(M;
E\tensor \L)$
 such that $${D}_{M,\L}+\mathcal{C}_{\mathcal{P}}\in \Psi^1_{A}(M; E\tensor
 \L)
\;\;\text{is invertible  in}\;\; \Psi^*_{A}(M;E\tensor \L).$$
Moreover, $\mathcal{C}_{\mathcal{P}}$ is symmetric with respect to
the $A$-valued scalar product on $C^\infty (M,E\tensor \L)$.
We
shall consider in what follows the set of allowable perturbations:
\begin{multline}\label{allowables}{\mathfrak P}:=
{\mathfrak P}_{{D}_{M,\L}}:=\{\mathcal{C}\in
\Psi^{-\infty}_{A}(M;E\tensor \L)\;\;
|\;\;\mathcal{C}\;\;\text{is symmetric}\\ \text{and}\;\;
{D}_{M,\L}+\mathcal{C}\;\;\text{is invertible in}\;\;
\Psi^*_{A}(M;E\tensor \L)\}.
\end{multline}


Summarizing, we have just seen that
\begin{equation}\label{eq:non_empty_perturb}
 \Ind({D}_{M,\L})=0 \;\;\implies \;\; {\mathfrak P}\ne
 \emptyset\,.
\end{equation}

Let us go back to the case where $M=\partial W$, with $W$ even
dimensional. As mentioned above, by cobordism invariance, the
index class of $D_{M,\L}$ is zero in $K_{1} (A)$; thus ${\mathfrak
  P}\ne \emptyset$. Let
$\mathcal{C}\in {\mathfrak P}$ and consider the invertible operator
$D_{M,\L} +\mathcal{C}$. According to  \cite[Lemma 6.1]{LPGAFA}
there exists a smoothing operator in the Mishchenko-Fomenko
$b$-calculus $$\mathcal{C}_W^+ \in \Psi^{-\infty}_{b,A}
(W;E^+\otimes \L, E^-\tensor \L)$$ such that ${D}_{\L}^+
+\mathcal{C}_W^+ $ has ${D}_{M,\L} +\mathcal{C}$ as boundary
operator. Using the invertibility of  ${D}_{M,\L} +\mathcal{C}$
one proves the invertibility of the {\it indicial family}
associated to $ {D}^+_{\L} +\mathcal{C}^+_W$; this property can in
turn be used to prove that $ {D}^+_{\L} +\mathcal{C}^+_W$ is
$A$-Fredholm as an operator between suitable Sobolev Hilbert
$A$-modules. Thus there is a well defined index class in $K_0
(A)$, denoted $\Ind_b ({D}_{\L} +\mathcal{C}_W)$. More precisely
\begin{equation}
\label{bindexclass} \Ind_b
(D_{\L}+\mathcal{C}_W)=[\mathcal{I}_+]-[\mathcal{I}_-]\,\in K_0
(A)
\end{equation}
with $\mathcal{I}_\pm $ finitely generated projective $A$-modules and
$$\mathcal{I}_\pm\subset x^\epsilon H^\infty_b (W, (E\tensor
L)^{\pm}),\quad \epsilon>0.$$ The construction of $\mathcal{C}_W$
from $\mathcal{C}$ involves choices; the index class, on the other
hand, does not depend on these choices. We refer the reader to
\cite{LPGAFA} for the details.

\begin{remark}
The results in
\cite{LPGAFA}, \cite{LP03} are noncommutative generalizations of
the results first proved in \cite{MPI}, \cite{MPII}  for
families of Dirac operators (i.e.~for $C(B)$-linear operators,
with $B$ a compact manifold).
\end{remark}

\begin{notation}
We shall often denote by $\Ind_b (D_L,\mathcal{C})$ the index class
in (\ref{bindexclass}).
\end{notation}

\begin{remark}
Notice that the index class $ \Ind_b ({D}_L,
\mathcal{C})$ {\it does depend} on $\mathcal{C}$ \footnote{For example if
$\mathcal{C}_{\mathcal{P}}\in {\mathfrak P}$ 
and $\mathcal{C}_{\mathcal{Q}}\in {\mathfrak P}$ 
are defined by spectral sections $\mathcal{P}$ and $\mathcal{Q}$
respectively, then (see \cite{LPAGAG}, \cite{LPFOLIATED})
$$\Ind_b ({D}_{\L}, \mathcal{C}_{\mathcal{P}})-
\Ind_b ({D}_{\L},
\mathcal{C}_{\mathcal{Q}})=[\mathcal{Q}-\mathcal{P}] \;\;\;
\text{in}\;\;\; K_0 (A)\,,$$
 where on the right hand side the
difference class defined by the two projections $\mathcal{Q}$,
$\mathcal{P}$ appears. This is known to be non-zero in general.}

\end{remark}

\begin{definition}
Let $\tr^{{\rm alg}}\colon K_0 (A)\rightarrow A_{{\rm ab}}$ be
the trace introduced in Subsection
\ref{sec:perturbations}. The $A_{\rm ab}$-valued b-index
is by definition $\tr^{{\rm alg}}(\Ind_b (D_{\L},\mathcal{C}))\in
A_{{\rm ab}}$. As in the closed case we use the notation
\begin{equation*}
\tr^{{\rm alg}}(\Ind_b (D_{\L},\mathcal{C}))=:\Ind_{b,[0]}
(D_{\L},\mathcal{C})= [\mathcal{I}_{+}]_{[0]} -
[\mathcal{I}_{-}]_{[0]} \in A_{{\rm ab}} \,.
\end{equation*}
\end{definition}

Our goal is to prove an index formula for $\Ind_{b,[0]}
(D_{\L},\mathcal{C})$. First of all, we
introduce the boundary correction term that will appear in the
formula.

\subsubsection{Eta invariants for perturbed operators.}
\label{subsec:eta-perturbed}

Let $\C\in {\mathfrak P}_{{D}_{M,\L}}$. Consider the
pseudo-differential operator
\begin{equation*}
D_{M,\L}+\C\in\Psi^1_A(M;E\tensor \L).
\end{equation*}
We define the $A_{\rm ab}$-valued eta invariant
$\eta_{[0]}(D_{M,\L}+\C) \in A_{{\rm ab}}$
by
\begin{equation}\label{perturbed-eta_general_A}
\eta_{[0]}(D_{M,\L}+\C):=\frac{1}{\sqrt{\pi}}\int_0^\infty
\TR \left ( (D_{M,\L}+\C)e^{-t (D_{M,\L}+\C)^2}
\right)\;\frac{dt}{\sqrt{t}} .
\end{equation}
The integral converges for large $t$ because of the
invertibility of $D_{M,\L}+\C$. The convergence for $t\downarrow
0$ follows from the local $A_{{\rm ab}}$-valued index
theorem on closed manifold and the observation that $$ e^{-t
(D_{M,\L}+C)^2}= e^{-t D_{M,\L}^2}+ t
C^{\infty}([0,\infty),\Psi^{-\infty}_A(M,E\otimes \L))\,,$$ a consequence of Duhamel's formula and
fact that $\C\in\Psi_A^{-\infty}(M,E\tensor \L)$.

\subsubsection{The $A/\overline{\Commutator{A}{A}}$-valued
 Atiyah-Patodi-Singer index formula}
\label{sec:general-degree-zero}

In Appendix \ref{sec:aps-index-theory-proof} we shall recall the precise form of the
Mishchenko-Fomenko decomposition theorem in the $b$-context.
From this result we get the two finitely generated
projective modules $\mathcal{I}_\pm$ entering into the definition
of the $b$-index class $\Ind_{b,[0]}(D_{\L},\mathcal{C})$.
$b$-elliptic regularity implies that the projections $P_\pm$ onto
these two modules have well defined traces $\TR (P_\pm)\in
A_{{\rm ab}}$. Analyzing the large-time behavior of the
$b$-supertrace of the heat-kernel for $(D_{\L}+\C_W)^2$, computing
the short-time limit and using the commutator formula for the
$b$-trace, one finally proves the following theorem (see
Appendix \ref{sec:aps-index-theory-proof} for proofs and relevant definitions)

\begin{theorem}\label{theo:general_APS}
  \begin{equation*}
    \begin{split}
      \Ind_{b,[0]}(D_{\L},\mathcal{C}) \equiv &
      [\mathcal{I}_+]_{[0]}-[\mathcal{I}_-]_{[0]}\\=& \lim_{t\to +\infty} b\STR (e^{-t
      (D_{\L}+\mathcal{C}_W)^2})=\TR (P_+) - \Tr (P_-) \\
      =& \int_W{AS(D)}(x)\wedge \ch(E)(x)\wedge \ch{\L}(x)\;dx -
      \frac{1}{2}\eta_{[0]}(D_{M,\L}+\C)  \in A_{{\rm ab}} .
    \end{split}
  \end{equation*}
  where $AS(D)$ is the
local integrand in the APS
  index theorem for $D$.
  \end{theorem}

\section{Discrete groups, twisting bundles and index theorems}
\label{sec:discrete groups}

 In this section we show how we can derive the
particular index theorems needed for the results stated in the
introduction from the general ones presented above. We shall
specialize to the maximal or the reduced $C^*$-algebra of the
group, and also to the group von Neumann algebra.

\subsection{The group $C^*$-algebra and the classical Mishchenko-Fomenko twisting bundle}
\label{sec:class-mishch-fomenko}

We now specialize the $C^*$-algebra $A$ and the Hilbert $A$-module
bundle $\L$ in
\ref{sit:A_bundle_with_boundary} and
\ref{set-up:A_bundle_without_boundary} to the cases which lead to
essentially all our applications.

Specifically, let $\Gamma$ be a discrete finitely generated
 group. Let $\Gamma\to \tilde M \to M$
be a Galois covering and let
  $u\colon M\to B\Gamma$ be the associated classifying map,
  i.e.~$\tilde M$ is the pull back of the universal $\Gamma$-covering
  $E\Gamma\to B\Gamma$ under $u$.
We twist our Dirac type operator $D$ with the flat
$C^*\Gamma$-bundle $\L$ associated to the covering, i.e.~with
$$\L:= \tilde M \times_\Gamma C^* \Gamma.$$

\begin{theorem}\label{theo:MS_AS_index}
  For this operator
  \begin{equation}\label{equation: MS AS index}
    \Ind_{[0]} (D_{M,\L}) = \left(\int_M AS(D)\right)\cdot 1 \in
    (C^*\Gamma)_{{\rm ab}}.
  \end{equation}
  Here, AS is the integrand in the classical
  Atiyah-Singer index theorem and we recall that $(C^*\Gamma)_{{\rm
      ab}}:= C^*\Gamma/\overline{\Commutator{C^*\Gamma}{C^*\Gamma}}$
\end{theorem}
\begin{proof}
  We apply Theorem \ref{theo:general_AS}. Since $\L$ is flat, and each
  fiber is isomorphic to the free $C^*\Gamma$-module of rank $1$,
  $\ch(\L)(x)=1\in (C^*\Gamma)_{{\rm ab}}$, so
  that we get the constant $1$ in the formula.
\end{proof}

We now look at a compact manifold $W$ with boundary $M$, again with a
$\Gamma$-covering $\Gamma\to \tilde W\to W$, classified by a map (also
called $u$) $u\colon W\to B\Gamma$. Note that this induces, by
restriction, a $\Gamma$-covering $\tilde M$ of the boundary.

Let $D$ be a Dirac type operator on $W$. Then let $\mathcal{C}\in {\mathfrak P}_{{D}_{M,L}}$ give rise to an allowable
perturbation of the boundary operator $D_{M,\L}$.
\begin{theorem}\label{theo:Cstart_Gamma_index_theorem}
  \begin{equation}\label{eq:aPS for Mishchenko twist}
\Ind_{b,[0]}(D_{\L},\C)= \left( \int_W AS(D) \right) \cdot 1 - \frac{1}{2}
\eta_{[0]}(D_{M,\L}+\C)\quad \in (C^*\Gamma)_{{\rm ab}}.
\end{equation}
 \end{theorem}

 \begin{proof}
This is a special
case of Theorem \ref{theo:general_APS}. Note that,
since the twisting
bundle $\L$ is flat and fiberwise isomorphic to the free Hilbert
$C^*\Gamma$-module of rank $1$, again $\ch{\L}(x)=1\in
(C^*\Gamma)_{\rm ab}$, so that it does only contribute the $1$ in the
formula.
 \end{proof}

\begin{remark}\label{remark:maximal-reduced}
Consider the reduced group $C^*$-algebra  $C^*_{r}\Gamma$ of
Definition \ref{def:group_Cstar_algebras}. One can form a
Mishchenko-Fomenko twisting bundle $\widetilde M\times_\Gamma
C^*_r \Gamma$ and prove $(C^*_r
\Gamma)_{\rm{ab}}$-valued index formulas
completely analogous to formulas (\ref{equation: MS AS index}) and
(\ref{eq:aPS for Mishchenko twist}).
\end{remark}

\section{Surjectivity of the Baum-Connes map and stable
rho-invariants}\label{sec:surjective+stable}

\subsection{The $\rho$-index and its
vanishing}\label{subsect:rho-index}

 Let $\Gamma$ be a discrete
group. Let $\Gamma\to \tilde W\to W$ be an even-dimensional Galois
covering of a compact manifold with boundary $W$, classified by $u\colon W\to
  B\Gamma$. Let $\Gamma\to
\tilde M\to M$ be the boundary covering.
 As before, we consider $\mathcal{D}:=D_{\L}$, a Dirac-type operator on $W$ twisted with
  the Mishchenko-Fomenko line bundle
  $\mathcal{L}$. Since the boundary operator $D_{M,\L}$ has trivial index in
  $K_{1}(C^*\Gamma)$ we can choose a trivializing perturbation $\mathcal{C}\in {\frak
    P}_{{D}_{M,\L}}$.
With respect to this perturbation,
 we
$\Ind_b (\mathcal{D},\C) \in K_0(C^*\Gamma)$. By applying $
\tr^{{\rm alg}}\colon K_0(C^*\Gamma)\to
(C^*\Gamma)_{{\rm ab}}$ we have finally
defined the $(C^*\Gamma)_{{\rm ab}}$-valued
$b$-index   $\Ind_{b,[0]} (\mathcal{D},\C) \in
(C^*\Gamma)_{{\rm ab}}$.

\smallskip
  Since
 $(C^*\Gamma)_{{\rm ab}}$ is a vector
  space with one dimensional subspace generated by
  $[1]=1+\overline{\Commutator{C^*\Gamma}{C^*\Gamma}}$, we can project the degree
  zero part of the index onto the quotient
  $(C^*\Gamma)_{{\rm ab}}/<[1]>$.

\begin{definition}\label{def:rho-index}
The  \emph{$\rho$-index} associated to $\D$ is the image of
$\Ind_{b,[0]} (\mathcal{D},C) $ in the quotient $\quotient$. We
shall denote the $\rho$-index by
\begin{equation}\label{display-rho-index}
\ind^{\rho}_{b,[0]}
(\mathcal{D},C)\in \quotient\,.
\end{equation}
\end{definition}

\begin{remark}
Suppose now that $W$ is {\it closed} and let
$\Ind_{[0]}(\mathcal{D})\in (C^*\Gamma)_{{\rm ab}}$ be
the associated $(C^*\Gamma)_{{\rm ab}}$-valued index.
Using the degree-zero Atiyah-Singer index formula
\ref{theo:MS_AS_index}, we see that
the $\rho$-index  of $\mathcal{D}$ vanishes:
\begin{equation}
\ind^{\rho}_{[0]}
(\mathcal{D})\,=\,0\in \quotient.
\end{equation}
This simple observation plays a fundamental role.
\end{remark}

\begin{lemma}\label{lemma:Atiyahs_index_theorem}
If $\Gamma$ is torsion-free and the Baum-Connes map $\mu_{max}:K_0
(B\Gamma)\to K_0 (C^*\Gamma)$ is surjective, then the $\rho$-index
vanishes: \begin{equation*}\ind^{\rho}_b (\mathcal{D},C)_{[0]} = 0
\;\;\text{in}\;\; \quotient\,.\end{equation*}
\end{lemma}
\begin{proof}
Observe that we can define a homomorphisms
  \begin{equation*}
\Phi^\rho \colon K_0(C^*\Gamma)\to \quotient\,,
\end{equation*}
by simply composing $ \tr^{{\rm alg}}: K_0(C^*\Gamma)\to
(C^*\Gamma)_{{\rm ab}}$ with the quotient map  $(C^*\Gamma)_{{\rm
ab}}\to \quotient$.
 Using the above remark, we see that this
  homomorphism
  is zero on the image of $\mu_{max}$, since the image of $\mu_{max}$
  consists precisely of index classes associated to {\it closed} manifolds.
  We use the Baum-Douglas \cite{Baum-Douglas-1}
  description of the K-homology of $B\Gamma$ here.
  By assumption, $\mu_{max}$ is surjective, so that $\Phi^\rho =0$ on
  all of $K_0 (C^*\Gamma)$. Thus
$0=\Phi^\rho (\Ind_b (\mathcal{D},C))=\ind^{\rho}_b
(\mathcal{D},C)_{[0]} $
   and the assertion follows.
\end{proof}

\subsection{The stable rho-invariant}\label{subsect:stable-rho}

 Assume that $\mathcal{D}$ is a Dirac type operator on a
 closed
  odd-dimensional  manifold $M$, twisted with
  the Mishchenko-Fomenko line bundle
  $\mathcal{L}$ associated to some classifying map $u\colon M\to
  B\Gamma$.
 Assume that  $\Ind(\mathcal{D})=0\in
K_{1}(C^*\Gamma)$. Pick a trivializing  perturbation ${C}\in
{\frak P}_{\D}$.

\begin{lemma}\label{lemma:independence-outlook}
If $\Gamma$ is  torsion-free and the
  Baum-Connes map
$\mu_{max}\colon K_0(B\Gamma)\to K_0(C^* \Gamma)$
  is surjective, then
\begin{equation}
  \label{eq:degree_zero_rho}
[\eta_{[0]}(\mathcal{D}+C)] \quad\in
  \quotient
\end{equation}
i.e. the image of $\eta_{[0]}(\mathcal{D}+C)\in (C^*\Gamma)_{{\rm ab}}$ in the
quotient $\quotient$,
 does
{\bf not} depend on the particular perturbation $\C$ chosen.
\end{lemma}

\begin{proof}
Let $\mathcal{C}^\prime\in {\frak P}$ be a different perturbation.
Consider the cylinder $[-1,1]\times M$. On the boundary of the
cylinder we have an invertible operator, obtained by considering
the operator $\mathcal{D}+\mathcal{C}$ on $\{-1\}\times M$ and the
operator $\mathcal{D}+\mathcal{C}^\prime$ on $\{1\}\times M$. As
already explained, there is a well defined  $b$-index class in
$K_0 (C^* \Gamma)$, obtained by lifting the two perturbations to
the cylinder and defining a $b$-pseudodifferential operator
$\mathcal{D}_{[-1,1]\times M }+\mathcal{C}_{[-1,1]\times M}$
with {\bf invertible} indicial family. We denote this index class
by
$$\Ind_b (\mathcal{D}_{[-1,1]\times M }+\mathcal{C}_{[-1,1]\times M})\,.$$
Using the  APS-index formula
\ref{theo:Cstart_Gamma_index_theorem} we obtain
  \begin{multline*}
\ind^\rho_b(\mathcal{D}_{[-1,1]\times M }+\mathcal{C}_{[-1,1]\times
  M})_{[0]}= \\
  -\frac{1}{2}\left(
  \rho_{[0]}(\mathcal{D}+\mathcal{C})-
  \rho_{[0]}(\mathcal{D}+\mathcal{C}')\right)\quad \in
  \quotient
\end{multline*}
On the other hand, by the assumed surjectivity of the Baum-Connes
map we know that $\ind^\rho_b (\mathcal{D}_{[-1,1]\times
M}+\mathcal{C}_{[-1,1]\times
  M})=0$
and the Lemma is proved.
\end{proof}

\begin{definition}\label{def:stable_rho}
Assume
 that $\Ind(\mathcal{D})=0\in K_1(C^*\Gamma)$
and that the max-Baum-Connes map $K_0 (B\Gamma)\to K_0 (C^*
\Gamma)$ is surjective. Then the \emph{stable degree-zero
$\rho$-invariant} is the class
\begin{equation}\label{l2rhostable}
\rho_{[0]}^s (\mathcal{D}):= [\eta_{[0]}(\mathcal{D}+{C})] \in
\quotient
\end{equation}
for any perturbation ${C}\in {\frak
    P}_{\mathcal{D}}$, and it is well defined because of Lemmas
  \ref{lemma:independence-outlook}.
\end{definition}




\begin{notation} Let $r\colon M\to B\Gamma$ be a classifying map.
We shall frequently use the notation $\rho_{[0]}^s (M,r)$ with the understanding that
this rho-invariant is associated either to the signature operator on
the oriented Riemannian manifold $M$, or to the spin Dirac operator on
the spin manifold $M$ (both twisted by the Mishchenko-Fomenko bundle),
depending on the context.
Sometime, we shall omit the classifying map from the notation thus
writing $\rho_{[0]}^s (M)$.
\end{notation}



\subsection{Bordism invariance}\label{subsec:bordism-invariance}

Let $(M,g)$ be odd dimensional endowed with a a bundle of Clifford modules
defining a Dirac-type operator $D_M$. Let  $(M,r\colon M\to B\Gamma)$ be
a Galois covering and let $\D:=D_{M,\L}$ be the twisted Dirac
operator in the Mishchenko-Fomenko calculus. Let
$(M^\prime,g^\prime)$ another odd dimensional Riemannian manifold
endowed with a Clifford module and let
$(M^\prime,r^\prime\colon M^\prime\to B\Gamma)$ be a Galois covering on
$M^\prime$; we let $\D^\prime:=D_{M^\prime,\L^\prime}$. Assume now
the existence of a bordism $(W, R:W\to B\Gamma)$ between
$(M,r:M\to B\Gamma)$ and $(M^\prime,r^\prime:M^\prime\to
B\Gamma)$; we also assume the existence of a Riemannian metric and
of a bundle of Clifford modules on $W$ restricting to the given data on $M$
and $M^\prime$.

\begin{proposition}\label{prop:bordsim-invariance-stable}
Assume $\Gamma$ is torsion free and such that the Baum-Connes
map $\mu_{max}$ is surjective. Assume that $\Ind (\D)=0$ in $K_1
(C^*\Gamma)$ so that the stable rho-invariant $\rho^s_{[0]}(\D)$
is well defined. By the bordism invariance of the index class we
know that $\Ind (\D^\prime)=0$ so that $\rho^s_{[0]}(\D^\prime)$
is also well defined. We have
\begin{equation*}
\rho^s_{[0]}(\D)=\rho^s_{[0]}(\D^\prime)
\end{equation*}
\end{proposition}

\begin{proof}
The argument establishing Lemma \ref{lemma:independence-outlook}
can be repeated here, provided we substitute the cylinder
$[-1,1]\times M$ there with the bordism $W$ here.
\end{proof}

\section{Injectivity of the Baum-Connes map and bordism}
\label{sec:injectivity}

\subsection{Statement of results}\label{subsec:statement-injectivity}

In contrast with the previous section we now consider only the spin
Dirac operator $\Di$ and the signature operator $D^{{\rm sign}}$ and  not arbitrary
generalized Dirac type operators.

\begin{definition}
  Consider $(M,u\colon M\to B\Gamma)$. If $d\in \NN\setminus \{0\}$, we denote by
  $d (M,u:M\to B\Gamma)$ the disjoint union of $d$ copies of $M$ with
  obvious induced map $du\colon dM\to B\Gamma$, and with obvious
  induced structure (orientation, etc.).
\end{definition}

\begin{proposition}\label{prop:bordism}
  Assume that $M$ is a closed oriented manifold, $u\colon M\to B\Gamma$ a
  classifying map. Let $\mathcal{D}^{{\rm sign}}$ be the associated
  Mishchenko-Fomenko signature operator. Assume that $\mu_{max}\colon
  K_*(B\Gamma)\tensor\rationals \to K_*(C^*\Gamma)\tensor\rationals$
  is injective, and assume that
$\Ind(\mathcal{D}^{{\rm sign}} )=0$ in $K_{\dim M} (C^* \Gamma)\tensor\rationals$.
Then there exists $d\in\NN\setminus \{0\}$ such that  $d\cdot [M,u\colon M\to B\Gamma]$
is bordant to
$$\cup_{j=1}^k (A_j\times B_j,r_j\times 1:A_j\times
B_j\to B\Gamma)$$
in $\Omega^{SO}_{\dim M}(B\Gamma)$,
with $\dim B_j= 4b_j$,
$\pi_1(B_j)=1$ and  $<L(B_j),[B_j]>={\rm sgn}(B_j)=0$. Here $r_j:A_j\to B\Gamma$
is a continuous map and $(r_j\times 1) (a,b):=r_j(a)$.
\end{proposition}

\begin{proposition}\label{prop:spin_bordism}
  Assume that $M$ is a closed spin manifold, $u\colon M\to B\Gamma$ a
  classifying map. Let $\cDi:=\Di_{\L}$ be the associated Mishchenko-Fomenko
  spin Dirac operator. Assume that $\mu_{max}\colon
  K_*(B\Gamma)\tensor\rationals \to K_*(C^*\Gamma)\tensor\rationals$
  is injective and that $\Ind(\cDi)=0 \in K_*(C^*\Gamma)\tensor \rationals$.
  Then there is
  $d\in\naturals\setminus\{0\}$ such that $d [M,u\colon M\to B\Gamma]$ is
  bordant to
  \begin{equation*}
    \bigcup_{j=1}^k (A_j\times B_j, r_j\times 1\colon A_j\times B_j\to
    B\Gamma)
  \end{equation*}
  in $\Omega^{spin}_{\dim M}(B\Gamma)$, with $\dim B_j= 4b_j$,
$\pi_1(B_j)=1$ and  $<\hat A(B_j),[B_j]>=0$. Here $r_j:A_j\to B\Gamma$
is a continuous map and $(r_j\times 1) (a,b):=r_j(a)$.
\end{proposition}

\begin{remark}
  In Propositions \ref{prop:bordism} and \ref{prop:spin_bordism}, the
  condition that the index of the Mishchenko-Fomenko operator in the
  K-theory of the maximal $C^*$-algebra vanishes can be replaced by
  the more familiar condition that the index in the K-theory of the
  reduced $C^*$-algebra vanishes, if we replace the assumption on the
  rational injectivity of the maximal Baum-Connes map
by the analogous one for the reduced $C^*$-algebra.
Moreover, it should
  be remarked that the rational injectivity of the  {\it complex}
  Baum-Connes map
is equivalent
to the rational injectivity of the {\it real} Baum-Connes
  map. This is a well known result, a proof can be found
  in \cite{math.KT/0311295}.

\end{remark}

We prove Propositions \ref{prop:bordism} and \ref{prop:spin_bordism}
with the same method from algebraic topology, which characterizes rationalized
homology
theories.

\subsection{Rational homology theories}
\label{sec:rati-homol-theor}

Let $h_*$ and $k_*$ be two generalized homology theories. We will be interested in
the examples $\Omega^{spin}_*$ of spin bordism, $\Omega^{SO}_*$ of
oriented bordism, $K_*$ and $KO_*$ of complex or real K-homology (the
homology theory dual to K-theory).

Let $\check{h}_*:= h_*(pt)$ and $\check{k}_*:=k_*(pt)$ be the
coefficients.

It is easy to see that $h_*\tensor \rationals$ and
$H_*(X;\rationals)\tensor_\rationals (\check{h}_*\tensor\rationals)$
(graded tensor product, the
$n$-th group is $\oplus_{p+q=n} H_p(X)\tensor \check{h}_q$) are again
homology theories, compare e.g.~\cite[3.18 ff]{Hilton}\footnote{One remark on
the reference \cite{Hilton}: there, only cohomology is considered, but
everything works exactly in the same way if cohomology is replaced by
homology.}.

\begin{proposition}\label{prop:general_trafos}
  \begin{enumerate}
  \item Every natural transformation $T\colon h_*\tensor\rationals\to
    k_* \tensor\rationals$
    is determined by $\check{T}\colon \check{h}_*\tensor\rationals\to
    \check{k}_*\tensor\rationals$.
  \item Every homomorphism $\check{T}\colon \check{h}_*\tensor\rationals\to
    \check{k}_*\tensor\rationals$ has a unique extension to a natural
    transformation $T\colon h_*\tensor\rationals\to k_*\tensor\rationals$.
  \item There is a unique natural transformation
    $H_*(\cdot;\rationals)\tensor (\check{h}_*\tensor\rationals) \to
    h_*(\cdot)\tensor\rationals$ which is the identity on the
    point. It is a natural isomorphism
of $\check{h}_*\tensor\rationals$-modules if $h_*$ is a
multiplicative homology theory.
\end{enumerate}
\end{proposition}
\begin{proof}
  Compare \cite[(3.20),(3.21), (3.22)]{Hilton}.
\end{proof}

Now, we apply this to the examples we want to study. By definition,
elements in $K_*(X)$ are represented by triples $(M,E,\phi)$, where
$M$ is a closed oriented Riemannian manifold of dimension congruent to $* \pmod{2}$, $E$
is a bundle of Clifford modules on $M$ and $\phi\colon M\to
X$ is continuous. Similarly, $KO_*(X)$ is defined by triples
$(M,E,\phi)$ as above, with the exception that $E$ is a
bundle of real Clifford modules. Of course, a suitable equivalence relation has to
be factored out. See \cite{Baum-Douglas-1} and \cite{keswani3}.

In particular, we get natural transformations
\begin{align*}
  &T_{SO}\colon \Omega^{SO}_*(X) \to K_*(X); & [M\xrightarrow{\phi} X] \mapsto
  [M,\Lambda^{{\rm sign}},\phi]\\
  & T_{spin}\colon \Omega^{spin}_*(X) \to KO_*(X); & [M\xrightarrow{\phi} X] \mapsto
  [M, \spi,\phi].
\end{align*}
Here, $\Lambda^{{\rm sign}}$ stands for the complex Clifford module
defining the signature operator on the oriented Riemannian manifold $M$, whereas
$\spi$ is the spin Clifford module defining the real Dirac operator of the given spin
structure\footnote{Topologists
refer to these transformations  as the
{\it natural K-orientation associated to an (ordinary) orientation} and
{\it natural KO-orientation associated to a spin structure} on $M$,
respectively.} .
If
$X$ is a point, $T_{SO}$ is just given by the signature (using the
isomorphism $K_{2n}(pt)=\integers$). Also if
$X=\{pt\}$, $T_{spin}$ is the $\alpha$-invariant, which is (up to a
multiple) equal to $\hat{A}(M)$ for dimensions divisible by $4$ (again
use the isomorphism $KO_{4n}(pt)=\integers$). Note that, in other
dimensions, $KO_k(pt)$ is finite, and $K_{2n+1}(pt)=0$.

Using Proposition \ref{prop:general_trafos} we get commutative
diagrams
\begin{equation}
  \label{eq:special_trafos_SO}
  \begin{CD}
    H_*(X;\rationals)\tensor (\Omega^{SO}_*(pt)\tensor\rationals) @>{{\id}\tensor
      T^{\rationals}_{SO}}>> H_*(X;\rationals)\tensor (K_*(pt)\tensor\rationals)\\
    @VV{\iso}V  @VV{\iso}V\\
    \Omega^{SO}_*(X)\tensor\rationals @>{T^{\rationals}_{SO}}>> K_*(X)\tensor\rationals.
  \end{CD}
\end{equation}
\begin{equation}
  \label{eq:special_trafos_spin}
  \begin{CD}
    H_*(X;\rationals)\tensor (\Omega^{spin}_*(pt)\tensor\rationals) @>{{\id}\tensor
      T^{\rationals}_{spin}}>> H_*(X;\rationals)\tensor (KO_*(pt)\tensor\rationals)\\
    @VV{\iso}V  @VV{\iso}V\\
    \Omega^{spin}_*(X)\tensor\rationals @>{T^{\rationals}_{spin}}>> KO_*(X)\tensor \rationals.
  \end{CD}
\end{equation}

It is not really necessary for us to understand completely the
vertical isomorphisms. The important point is that
$\Omega^{spin}_*(X)\tensor\rationals$ is a \emph{free} module over
$\Omega^{spin}_*(pt)\tensor\rationals$ (free generators given by a basis of
$H_*(X;\rationals)$ under the vertical isomorphism). The corresponding
statements hold for $\Omega^{SO}_*(X)\tensor\rationals$ and $K_*(X)\tensor\rationals$,
$KO_*(X)\tensor\rationals$. Moreover, we also see that the basis can
be chosen in a natural way and such that the transformations $T^{\rationals}_{SO}$
and $T^{\rationals}_{spin}$ are diagonal with respect to them. In particular,
$\ker(T^{\rationals}_{SO}(X))$ is a free module over $\ker(T^{\rationals}_{SO}(pt))$ and
$\ker(T^{\rationals}_{spin}(X))$ is a free module over $\ker(T^{\rationals}_{spin}(pt))$.

This means that every element in the kernel of $T^{\rationals}_{SO}$
or $T^{\rationals}_{spin}$,
respectively, will after multiplication with a suitable non-zero
integer (to clear the denominators) have the form described in Proposition
\ref{prop:bordism} or Proposition \ref{prop:spin_bordism},
respectively.
Here we use the fact that the $\check{\Omega}_*$-module structure of
$\Omega_* (X)$ is given by the  Cartesian product.

Finally, recall that
\begin{equation*}
\Ind(\mathcal{D}^{{\rm sign}})=
\mu_{max}(T_{SO}([M\xrightarrow{u}B\Gamma]))\in K_{\dim M}(C^*\Gamma))
\end{equation*}
in the situation of Proposition \ref{prop:bordism},
and that
\begin{equation*}
\Ind(\cDi_{\reals})=
\mu_{\reals,max}(T_{spin}([M\xrightarrow{u}B\Gamma]))\in KO_{\dim M}(C_{\reals}^*\Gamma)
\end{equation*}
in the situation of Proposition \ref{prop:spin_bordism}. Because we
assume that $\mu_{max}\tensor \rationals$
is injective, Proposition \ref{prop:bordism} on the signature operator
follows  from the
descriptions of $\ker(T^{\rationals}_{SO})$ we have just
given.
Concerning
Proposition  \ref{prop:spin_bordism},
we point out that the rational vanishing
of the complex index class associated to $\cDi$
is equivalent to the rational vanishing of the corresponding real
index
class for the real spin Dirac operator $\cDi_{\reals}$.
Recalling that the rationally injectivity of the real
Baum-Connes map $\mu_{\RR,max}$ is equivalent to the rationally
injectivity
of the complex Baum-Connes map $\mu_{max}$ \cite{math.KT/0311295}, we see
that the result stated in Proposition \ref{prop:spin_bordism}
is again  a  consequence of the
description of  $\ker(T^{\rationals}_{spin})$ we have just
given.

\section{Vanishing of stable $\rho$-invariants for certain product manifolds}
\label{sec:vanishing-stable-rho}

\begin{theorem}\label{theo:general_atimesb}
Assume that $M=U\times V$, where $V$ is a simply connected oriented
Riemannian manifold of
dimension divisible by $4$ with vanishing signature, and $U$ is any
oriented Riemannian manifold of odd dimension.
Let $A$ be a $C^*$-algebra and $\L$ a \emph{flat} finitely generated
projective Hilbert $A$-module bundle on $U$ (defining a bundle, also
called $\L$, on $M$ by pullback).
Let $D_{\L}$ be the corresponding signature operator.

Then $\Ind D_{\L}=0 \in K_{\dim M}(A)$
and there exists a perturbation $\mathcal{C}\in {\frak
    P}_{{D}_{\L}}$ such that
  \begin{equation*}
    \eta_{[0]}(D_{\L}+\C) = 0 \in A_{{\rm ab}}.
  \end{equation*}
\end{theorem}
\begin{proof}
For the differential forms on $U\times V$ we
have the following decomposition:
$$ L^2\Omega(U\times V)= L^2\Omega(U)\otimes L^2\Omega^+(V) \oplus
L^2\Omega(U)\tensor L^2\Omega^-(V)$$
where $\Omega^\pm (V)$ refers to the grading associated to the
signature operator on $V$.
Since the bundle $\L$ pulls back from $V$, we get a corresponding
decomposition of the Hilbert $A$-module of forms with coefficients in $\L$:
$$ L^2\Omega(U\times V;\L)= L^2\Omega(U;\L)\otimes L^2\Omega^+(V) \oplus
L^2\Omega(U;\L)\tensor L^2\Omega^-(V).$$
With respect to this decomposition, the signature operator splits as
$$ D_{U\times V}= \begin{pmatrix} D_U\tensor 1_V & 1_U\tensor D_V^-\\
                    1_U\tensor D_V^+ & -D_U \tensor 1_V
\end{pmatrix} = D_U\tensor \tau + 1_U\tensor D_V,$$
where $\tau$ is the grading operator on $L^2\Omega(V)$. Since $\L$ is
flat, we get a corresponding splitting of the signature operator
twisted with $\L$:
\begin{equation}
\label{eq:split_fo_D_LUteimsV}
 D_{\L,U\times V}= \begin{pmatrix} D_{\L,U}\tensor 1_V & 1_U\tensor D_V^-\\
                    1_U\tensor D_V^+ & -D_{\L,U} \tensor 1_V
\end{pmatrix} = D_{L,U}\tensor \tau + 1_U\tensor
D_V.
\end{equation}

Let us now define the perturbed operator. Choose an isometry
$\Psi\colon \ker(D_V^+)\to \ker(D_V^-)$. Since $V$ is compact with
vanishing signature, such an isomorphism of finite dimensional vector
spaces exists. Then
\begin{equation*}
C:=  \begin{pmatrix} 0 & \Psi^* \\ \Psi & 0 \end{pmatrix} \colon
L^2\Omega^+(V)\oplus L^2\Omega^-(V)\to L^2\Omega^+(V)\oplus
L^2\Omega^-(V)
\end{equation*}
is a smoothing operator on $V$ (here ``$\Psi$'' stands for the
projection onto $\ker(D_V^+)$ follows by the operator $\Psi$).

The operator $D_{\L, U\times V} + 1_U \tensor C$
is invertible.
Let $f( D_{\L,U})$ be any smoothing function of $D_{\L,U}$ which is
the identity on a sufficiently large neighborhood of zero in the
spectrum of $D_{\L,U}$ so that
\begin{equation*}
\D_C:= D_{\L, U\times V} + f(
D_{\L,U})\tensor C
\end{equation*}
becomes invertible.
This will be our perturbation; notice that
since we have found a smoothing perturbation of $D_{\L,U\times V}$
which is invertible, we conclude that $\Ind (D_{\L,U\times V})=0$ in
$K_{\dim M} (A)$, as required.

We want to show that the degree $0$ eta invariant of this perturbation
$\D_C$
vanishes. Since $L^2\Omega(V)$ is an ordinary Hilbert space where
orthogonal complements exist, we get
another orthogonal decomposition into projective Hilbert $A$-modules
\begin{equation*}
  L^2\Omega(U\times V;\L) = \left( L^2\Omega(U;\L)\tensor \ker(D_V)
  \right) \oplus \left(
  L^2\Omega(U;\L)\tensor \ker(D_V)^\perp\right).
\end{equation*}

The operator $D_{\L,U\times V}$ as well as  $f( D_{\L,U})\tensor C$
preserves the splitting $$L^2\Omega(U;\L)\tensor \ker(D_V)\oplus
L^2\Omega(U;\L)\tensor (\ker(D_V)^\perp)\,.$$ Moreover, the second
operator is zero on
the second summand. Therefore, the integrand in the definition of
$\eta_{[0]}(\D_C)$ splits into two summands. We continue by
investigating the two summand separately. We will use freely from
Equation \eqref{eq:split_fo_D_LUteimsV}: $D_{\L,U\times V} = D_{\L,U}\tensor \tau + 1_U\tensor D_V$.
We restrict to $L^2\Omega(U;\L)\tensor \ker(D_V)^\perp$.
We denote the restriction of $\D_C$
to this space
by $(\D_C)_r$ and the restrictions of $\tau$ and $D_V$ to
$\ker(D_V)^\perp$ by $\tau_{r}$ and $(D_V)_{r}$.
  Then
\begin{equation}
\begin{split} (\D_{C})_r \, e^{-t (\D_{C})_r^2} &=
   D_{\L,U\times V}
  e^{-t  D_{\L,U}^2}\tensor e^{-t (D_V)_r^2} \\
  &=  D_{\L,U}e^{-t D_{\L,U}^2} \tensor (\tau)_r e^{-t(D_V)_r^2} +
  e^{-t D_{\L,U}^2}\tensor (D_V)_r e^{-t(D_V)_r^2} .
\end{split}
\end{equation}
We claim that the degree zero trace with values in
$A_{{\rm ab}}$ of both summands vanishes. Note that $\TR$ is
multiplicative for tensor products. Here, we get $\TR$ of the first
factor multiplied with the ordinary complex valued trace of the second
factor.

For the first summand
$$\Tr((\tau)_r e^{-t(D_V)_r^2})=0$$
by the McKean-Singer formula (we
restrict to the orthogonal complement of the kernel, so the index is
zero); therefore
$$\TR(D_{\L,U}e^{-t D_{\L,U}^2} \tensor (\tau)_r e^{-t(D_V)_r^2})= \TR(
D_{\L,U} e^{-t D_{\L,U}^2}) \Tr((\tau)_r e^{-t(D_V)_r^2}) =0. $$
For the
second summand, the same argument applies since $\TR(e^{-t
  D_{\L,U}^2})=0$ on the odd dimensional manifold $U$ by the usual
symmetry argument.

 Now, we restrict $\D_C$ to $L^2\Omega(U;\L)\tensor\ker(D_V)$
(and we again denote the relevant restriction with $(\;\;)_r$).

We will finish by showing that this summand also is identically
equal to zero, which then implies immediately that
$\eta_{[0]}(\D_C)=0$.  To see this, consider the operator
$$B:=1_U\tensor (\tau\circ C)\colon L^2\Omega(U;\L)\tensor \ker(D_V)\to
L^2\Omega(U;\L)\tensor \ker(D_V).$$
This is an isometry which
anticommutes with $D_{\L,U}\tensor \tau$ (since $C$ anticommutes and
$\tau$ commutes with $\tau$) and which anticommutes with $1_U\tensor
C$ (since $C$ commutes and $\tau$ anticommutes with $C$).
Similarly, $B$ anticommutes with $f(D_{\L,U})\tensor C$.
Consequently, $B$ anticommutes with the restriction of $\D_C$ to
$L^2\Omega(U;\L)\tensor \ker(D_V)$, which is just given by
$D_{\L,U}\tensor \tau+ f(D_{\L,U})\tensor C$. Of course, this implies
also that $B$ commutes with $e^{-t (D_{\L,U\times V})_r^2}$. Therefore
(everything restricted to $L^2\Omega(U;\L)\tensor \ker(D_V)$)
\begin{equation*}
\begin{split}  \TR((\D_C)_r e^{-t (\D_{C})_r^2}) =&
  \TR(B^{-1}B (\D_{C})_r e^{-t
    (\D_{C})_r^2})\\
  =& - \TR (B^{-1}(\D_{C})_r e^{-t (\D_{C})_r^2}B)\\
  =& - \TR ((\D_{C})_r e^{-t(\D_{C})_r^2}) \in A_{{\rm ab}}.
\end{split}
\end{equation*}
The last equality follows from the trace property.
Consequently, the integrand also vanishes identically.


Putting everything together, we get
\begin{equation*}
\begin{split}
  \eta_{[0]}(\D_C) =& \frac{1}{\sqrt{\pi}}\int_0^\infty t^{-1/2}
\TR(\D_{C} e^{-t \D_{C}^2}) = \frac{1}{\sqrt{\pi}}
\int_0^\infty 0 =0.
\end{split}
\end{equation*}

\end{proof}

As a special case of the theorem, we obtain the following Corollary.

\begin{corollary}\label{theo:atimesb}
Assume that $M=U\times V$, where $V$ is a simply connected oriented
4k-dimensional Riemannian manifold such that $\langle L(V),[V]\rangle$=0, and $U$ is any
oriented Riemannian manifold.
Let $r\colon U\to B\Gamma$ be a continuous map.
Then
there exists a smoothing  perturbation $\C$ such that
$\mathcal{D}+\C$ is invertible and
$\eta_{[0]}(\mathcal{D}+\C)$ vanishes, where $\mathcal{D}$ is
the Mishchenko-Fomenko signature operator associated to $r\times 1$.
In particular, if the Baum-Connes map is surjective, then the stable
rho-invariant
is defined in $C^*\Gamma_{\rm{ab}}/<[1]>$ and equal to zero:
$\rho^s_{[0]}({U}\times V,r\times 1)=0 $.

\end{corollary}

There are similar results in the spin context:

\begin{theorem}\label{theo:atimesb_Dirac}
Assume that $M=U\times V$, where $V$ is a simply connected spin manifold of
dimension divisible by $4$ with $\hat{A}(V)=0$, and $U$ is any
spin manifold of odd dimension.
Let $A$ be a real $C^*$-algebra and $\L$ a \emph{flat} finitely generated
projective Hilbert $A$-module bundle on $U$ (defining a bundle, also
called $\L$, on $M$ by pullback).
There is a product metric on $U\times V$ such that for the resulting
twisted spin Dirac operator
$D_\L$ we get
  \begin{equation*}
    \eta_{[0]}(D_\L) = 0 \in A_{{\rm ab}}.
  \end{equation*}
\end{theorem}
\begin{proof}
  By Stephan Stolz' solution \cite{Sto3} of the Gromov-Lawson-Rosenberg conjecture
  for simply connected manifolds, $V$ admits a metric with positive
  scalar curvature. If we shrink this metric appropriately, the
  corresponding product metric on $U\times V$ will still have positive
  scalar curvature. The Lichnerowicz formula then implies that its
  spin Dirac operator $D_L$ is invertible (here we use that the bundle
  is flat). In particular, we don't need a perturbation and can define
  $\eta_{[0]}(D_{\L})$.

  Secondly, we can use the product structure (of the Dirac operator)
  and argue in a way similar to the proof of Theorem
  \ref{theo:atimesb} (but with the simplification that the operator on
  $V$ has no kernel) to conclude that
  \begin{equation*}
    \eta_{[0]}(D_L) =0 \in (C^*\Gamma)_{{\rm ab}}.
  \end{equation*}
\end{proof}

\section{Vanishing of stable rho-invariants}
\label{sec:vanishing-stable-rho-1}

\begin{theorem}\label{theo:vanish-stable}
Let $M$ be an odd dimensional oriented closed compact manifold
and let
$u\colon M\to B\Gamma$ be a continuous map, classifying a Galois
$\Gamma$-covering $\tilde{M}\to M$.
If the max-Baum-Connes map
$  \mu_{max}\colon K_*(B\Gamma) \to K_*(C^*\Gamma)$
is bijective and the signature
index class of $(M,u\colon M\to B\Gamma)$ vanishes in $K_{\dim M}(C^*\Gamma)$ then
\begin{equation}\label{form:vanish-stable-sign}
\rho^s_{[0]} (M,u)=0.
\end{equation}
\end{theorem}

\begin{proof}
We use the injectivity of the max-Baum-Connes map, the assumption
$\Ind (\mathcal{D})=0$ and Proposition
\ref{prop:bordism} in order to conclude that there exists a bordism
$(W,F:W\to B\Gamma)$ between $(dM,du\colon M\to B\Gamma)$ and
$$\cup_{j=1}^k (A_j\times B_j,r_j\times 1:A_j\times B_j\to B\Gamma)$$
with $\pi_1(B_j)=1$ and  $<L(B_j),[B_j]>=0$.
Denote briefly
\begin{equation*}
\cup_{j=1}^k (A_j\times B_j,r_j\times 1:A_j\times B_j\to B\Gamma)
\end{equation*}
by $(N,v\colon N\to B\Gamma)$. By cobordism invariance, the index of
the signature operator associated to $(N,v\colon N\to B\Gamma)$
is zero in $K_{\dim M}(C^*\Gamma)$. From the surjectivity
of the max-Baum-Connes map we know that there are well
defined stable rho-invariants $\rho_{[0]}^s (M,u)$,
$\rho_{[0]}^s (N,v)$. Fix  allowable perturbations
$\mathcal{C}_M$, $\mathcal{C}_N$. Then there exists a well
defined signature $b$-index class,
in $K_{\dim W}(C^*\Gamma)$,
associated to the bordism $(W,F:W\to B\Gamma)$
and to the perturbations
$d\mathcal{C}_M$, $\mathcal{C}_N$, which give rise to
a Mishchenko-Fomenko $b$-smoothing operator $\mathcal{A}_W$ on $W$.
We denote this $b$-index class
by
$\Ind_b (\mathcal{D}_W+\mathcal{A}_W)\in K_{\dim W}(C^*\Gamma)\,.$
Now we proceed as in the proof of Lemma \ref{lemma:independence-outlook}.
On the one hand, by the surjectivity of $\mu_{{max}}$
we have
$$ \Ind^{\rho}_b (\mathcal{D}_W+\mathcal{A}_W)=0 \in
(C^*\Gamma)_{{\rm ab}}/<1> \,;$$
on the other hand by applying the APS index theorem
\ref{theo:Cstart_Gamma_index_theorem}
 we get
$$ \Ind^{\rho}_b (\mathcal{D}_W+\mathcal{A}_W))
 = -\frac{1}{2} \left(d\rho^s_{[0]} (M,u)-
\rho^s_{[0]} (N,v)\right)
\in (C^*\Gamma)_{{\rm ab}}/<1> \,;$$
from which we deduce  that
$$ \rho^s_{[0]} (M,u)=\frac{1}{d} \rho^s_{[0]} (N,v) = \frac{1}{d} \sum_{j=1}^k
\rho^s_{[0]} (A_j \times B_j , r_j\times 1)$$
We finish the proof by applying Corollary \ref{theo:atimesb}.
\end{proof}

\begin{theorem}\label{theo:vanish_Dirac_stable}
Let $M$ be an odd dimensional closed compact spin manifold
and let
$u\colon M\to B\Gamma$ be a continuous classifying map. Let $\cDi:=\Di_{\L}$ be
the associated Mishchenko-Fomenko spin Dirac operator.
If the max-Baum-Connes map
 $ \mu_{max}\colon K_*(B\Gamma)
\to K_*(C^*\Gamma)$ 
is bijective and the
index class $\Ind(\cDi)=0 \in K_{\dim M}(C^*\Gamma)$ vanishes, then
\begin{equation}\label{form:vanish-stable-Dirac}
\rho^s_{[0]} (M,u)=0.
\end{equation}
\end{theorem}
\begin{proof}
  The proof follows exactly the same pattern of the proof of Theorem
  \ref{theo:vanish-stable}. More details in Subsection
  \ref{subsect:higher}
(see in particular the proof of Theorem \ref{theo:vanishing-higher}).
\end{proof}

\section{Unstable rho invariants}
\label{sec:unst-rho-invar}

In this section, let $A$ be a von Neumann algebra. Let $Z$ be a
commutative von Neumann algebra, and $\tau\colon A\to Z$ a positive
and normal
trace on $A$ with values in $Z$. 

In Section \ref{sec:surjective+stable} we managed to define, if the
index class is zero and the Baum-Connes map is surjective, the stable
rho-invariant $\rho_{[0]}$.
In Section
\ref{sec:vanishing-stable-rho-1} we showed  that, under the additional
assumption of injectivity of the Baum-Connes map, it is zero.

We now introduce unstable rho-invariants which are potentially more
interesting, since they are defined under much more general
hypothesis.


Assume that $M$ is a
closed manifold of odd dimension and let $D$ be  a Dirac type operator on $M$, acting on
sections of a bundle $E$. Let $\L$ be a
bundle of finitely generated projective Hilbert $A$-modules on $M$,
with a connection preserving all the structure. Let $D_\L$ be the
corresponding twisted Dirac operator. Each fiber of $\L$ is a finitely
generated projective module; $\L_x = \Im (p_x)$, with $p_x$ a
projection in $M_{k\times k}(A)$ for some $k$;
 let
  $\tau ({\L}_x):= \tau (p_x)$
where we extend the trace $\tau$ to $M_{k\times k}(A)$
in the obvious way.
We assume that  $\tau ({\L}_x)$
is constant in $x$.

\begin{definition}\label{def:tau-eta}
We define the $\tau$-eta invariant as
\begin{equation}
  \label{eq:def_eta_tau}
  \eta_\tau(D_\L):= \frac{1}{\sqrt{\pi}}\int_0^\infty \tau(\TR(D_\L
  \exp(-t D_\L^2))) \; \frac{dt}{\sqrt{t}}\;\in\;Z\,.
\end{equation}
\end{definition}

We have to check that this integral converges. For $t\to 0$, this follows from the
usual local heat expansion. Since $\tau$ is positive and
normal, the estimates of Cheeger and Gromov \cite{ChGr} can be used to
obtain convergence for large times.

\begin{definition}\label{def:of_general_rho}
Let $A_1 $ and $A_2$ be von Neumann algebras;
let $\L_1, \L_2$ be two Hilbert $A_j$-module bundles over $M$;
let $\tau_j\colon A_j
\to Z_j$ be positive normal traces with values in commutative
von Neumann algebras $Z_j$.
Let $\beta_j \colon Z_j\to V$ be homomorphisms to a fixed
target space $V$ (a vector space). Assume that
$\beta_1 (\tau_1 ((\L_{1})_x)) + \beta_2 (\tau_2 ((\L_2)_x))=0 \in V$.
The \emph{unstable rho invariant}
 of $D_\L$ with respect to $
    \beta_j$, $\tau_j$ and $\L_j$ is defined as
  \begin{equation*}
    \rho_{(\beta_j,\tau_j, \L_j)}(D) := \beta_1 (\eta_{\tau_1}(D_{\L_1}))
+ \beta_2 (\eta_{\tau_2}(D_{\L_2}))
\in V  \,.
\end{equation*}
The definition can be extended to a finite number of summands in the
evident way.
\end{definition}

We give examples showing the interest of such a definition.

\smallskip
\noindent
{\bf General example.}
Let $M$ be a closed oriented Riemannian manifold of odd
  dimension. Assume that $\Gamma$ is a
  discrete group; let $u\colon M\to B\Gamma$ be a map classifying a
  covering $\tilde M= u^* E\Gamma$. Let
  $\alpha_j\colon C^*\Gamma\to A_j$ be homomorphisms to
  unital von Neumann
  algebras $A_j$ (with $j=1,2$), and $\tau_j \colon
  A_j\to Z_j$ positive normal traces with
  values in  abelian von Neumann algebras $Z_j$.
Let $\L_j:= \tilde M\times_\Gamma A_j$ be the associated Hilbert
  $A_j$-module bundle, where $\Gamma$ acts on $A_j$ via $\Gamma\to
  C^*\Gamma\xrightarrow{\alpha_j} A_j$. Then $\tau_j ((\L_j)_x) =\tau_j
  (\alpha_j (1))$. If
  $\beta_j\colon Z_j\to V$ satisfy
  \begin{equation*}
\beta_1 (\tau_1(\alpha_1 (1)))  + \beta_2 (\tau_2 (\alpha_2 (1)))  = 0 \in V
\end{equation*}
then the  unstable rho-invariant $\rho_{(\beta_j,\tau_j,\L_j)}(D) \in
V$ is well defined and equal to
 \begin{equation}\label{definition-of-urho}
\beta_1 (\eta_{\tau_1}(D_{\L_1}))  + \beta_2 (\eta_{\tau_2} (D_{\L_2}))  \in V.
\end{equation}

\begin{example}\label{example:finite-repres} (Atiyah-Patodi-Singer
  rho-invariant.)\\
We refer to the {\it general example}.
  The relevant von Neumann algebras here are $A_1=M_d(\complexs)=A_2$, with two
  homomorphisms $\alpha_j \colon C^*\Gamma\to
  M_d(\complexs)$
induced  by two representations $\lambda_1,
  \lambda_2\colon
\Gamma\to U(d)$. The
  relevant trace is (in both cases) the usual trace $\tau\colon
  M_d(\complexs)\to\complexs$ on
  matrices. Then $\tau(\alpha_1(1))=d=\tau(\alpha_2(1))$, so that,
  with $V=\CC$, we can choose
  $\beta_1=\id$ and $\beta_2=-\id$. By Equation \eqref{eq:aps_equals_tr}, the eta-invariants
appearing in the definition of the APS-rho invariant
 and those appearing in formula
  \eqref{definition-of-urho}
  coincide. Thus with these choices
$\rho_{(\beta_j,\tau_j,\L_j)}(D) =\rho_{\lambda_1 - \lambda_2} (D)$.
\end{example}

\begin{example}\label{example:center_valued} (center-valued rho-invariant.)\\
  We refer to the {\it general example}.
Let  $A_1:=\NeumannN\Gamma$ be the group
  von Neumann algebra of $\Gamma$ and $Z$ its center.
Let $A_2=\CC$.
Let $\alpha_1$ be induced by  the natural map $C^*\Gamma \to \NeumannN\Gamma$
and let $\alpha_2\colon C^*\Gamma\to\complexs$ be
  induced by the trivial
  representation. Then (by definition) $\mathcal{L}_1=\NeumannN:=\tilde
  M\times_{\Gamma}\NeumannN\Gamma$ and $\mathcal{L}_2=M\times
  \complexs$.

We take $V=Z$, $\tau_1=\tau\colon\NeumannN\Gamma\to Z$
  equal to  the canonical center
  valued trace (compare \cite[Chapter 8]{MR98f:46001b}); $\tau_2\colon\complexs\to Z$ given by
$\tau_2 (z):=z\cdot 1$. Note that both are positive and normal
traces. Here $\tau_1(\alpha_1(1))=1 =\tau_2(\alpha_2(1))$, so
  that we can choose again  $\beta_1=\id$, $\beta_2=-\id$.
By Lemma \ref{lem:stable_eta_preserved},
  \begin{equation}\label{eq:center_valued_rho}
  \rho_{(\beta_j,\tau_j,\L_j)}(D) =  \eta_\tau(D_\NeumannN) - \eta(D)\cdot 1 \in Z
  \end{equation}
We define this element in $Z$ to be the {\it center-valued rho
  invariant} of $M$.
\end{example}

\begin{example}\label{example:delocalized} ($L^2$-rho and delocalized
  eta invariants.)\\
Let $<g>$ be a finite conjugacy class in $\Gamma$. This defines a
trace
\begin{equation}\label{eq:defin_of_deloc_trace}
  \tau_{<g>}\colon\NeumannN\Gamma\to\complexs;\, \sum_{h\in\Gamma}
  \lambda_h h\mapsto \sum_{h\in <g>} \lambda_h.
\end{equation}
By the universal property of the central valued trace $\tau$,
$\tau_{<g>}=\tau_{<g>}\circ \tau$, where we use the restriction of
$\tau_{<g>}$ to the center $Z$ on the right hand side.
We now apply $\tau_{<g>}$ to both sides of Equation
\eqref{eq:center_valued_rho}. If $g\ne 1$, $\tau_{<g>}(1)=0$, and we
obtain (using Proposition \ref{prop:compute_L2_eta})
\begin{equation*}
  \tau_{<g>}(\rho_{(\beta_j,\tau_j,\mathcal{L}_j )}(D)) =
  \eta_{<g>}(\tilde D).
\end{equation*}
If $g=1$, then $\tau_{<g>}=: \tau_{\Gamma}$ is the canonical trace,
and we obtain by Proposition \ref{prop:compute_L2_eta}, using
$\tau_\Gamma(1)=1$,
\begin{equation*}
  \tau_{\Gamma}(\rho_{(\beta_j\tau_j\mathcal{L}_j)}(D)) =
  \rho_{(2)}(\tilde D).
\end{equation*}

\end{example}

Our goal is to prove vanishing results for such generalized
unstable
rho-invariants and, thanks to the examples just given,
 derive the assertions of the Introduction as
special cases.




\section{Special perturbations of the signature operator}
\label{sec:homot-invar-sign}

Let $f:M\to M^\prime$ be a smooth orientation preserving homotopy
equivalence between two closed manifolds.
Let $A$ be a $C^*$-algebra and $\mathcal{V}'$ a flat bundle of
finitely generated Hilbert $A$-modules on $M'$.
We consider the signature operator on $M^\prime$ with values in $\mathcal{V}'$.
We denote
this operator by $D^\prime_{\mathcal{V}'}$. Next, we consider the signature
operator on $M$ with values in the flat bundle $\mathcal{V}:=f^* \mathcal{V}^\prime$
and we denote it by $D_{\mathcal V}$.
These two operators come from the de Rham complexes and the Hodge star
operator on $M^\prime$
and $M$ with
values in the flat bundles $\mathcal{V}^\prime$ and $f^*
\mathcal{V}^\prime$,
respectively. We  denote these de Rham differentials
by $d^\prime_{\mathcal{V}^\prime}$ and $d_{\V}$.

In this section we shall construct an
{\it explicit} trivializing perturbation for the Mishchenko-Fomenko
signature operator
\begin{equation}\label{signature-on-union}
\D:= \begin{pmatrix} D_{\mathcal{V}} & 0
\cr 0  & -D^\prime_{\mathcal{V}'} \cr
\end{pmatrix}
\end{equation}
on $M\sqcup (-M^\prime)$. 

Recall that, since  $f$ is an
orientation preserving homotopy equivalence,  
 the index in
$K_*(A)$ of $\D$ is zero\footnote{This result,
the homotopy
invariance of the signature index class,
has been established by several people and with different techniques.
Mishchenko and Kasparov prove it by
showing its equality with the $C^*$-algebraic Mishchenko symmetric signature,
an a-priori homotopy invariant.
Kaminker-Miller give a more analytical treatment,
adapting to  the noncommutative context
the proof that Lustzig gives  in the case $\Gamma=\ZZ^k$;
Hilsum -Skandalis prove the homotopy invariance in a purely 
analytical fashion; we shall follow their approach.
For a thorough treatment of the homotopy invariance of the
signature index class and its connections
with surgery theory we refer the reader to 
the recent  papers by Higson and Roe \cite{Hig-Roe}}.
Thus, we
already  know that there exists
a smoothing perturbation of the operator which is invertible,
see  \cite{WuI} \cite[Theorem 3]{LP03}.
In this section we shall sharpen this result, showing
that we can construct certain {\it special}  perturbations which are
{\it spectrally concentrated near zero} (definition below). 
We will use these special perturbations in Section
\ref{sec:stable-=-unstable} in order to show that, in this situation, the stable
rho invariant of $\D$ in (\ref{signature-on-union}) coincides with the unstable one, getting a
more precise result about their eta-invariants in that case.

\begin{theorem}\label{theo:odd_good_signature_perturbation}
  Let $f:M\to M^\prime$ be a smooth orientation preserving homotopy
equivalence between two closed manifolds, with
B  $\dim(M)=\dim(M')$  odd. Then, for each $\epsilon>0$ we can
  find a special self-adjoint smoothing perturbation $\B_\epsilon$ such that
  $\D+\B_\epsilon$ is invertible and with the additional
property that $\B_\epsilon$ is
$\epsilon$-spectrally concentrated, that is 
$$\B_\epsilon\circ
  \phi(\D)=0=\phi(\D)\circ \B_\epsilon$$ for each function
  $\phi\colon\reals\to\reals$ with $\phi(t)=0$ for $\abs{t}<\epsilon$
  \end{theorem}

\begin{remark}
  We  expect that a corresponding result holds if $\dim(M)$ is
  even. In this case we will also need that
  $\B_\epsilon$ is an {\it odd} operator with respect to the signature
  grading (where we use on $M'$ the reverse orientation). Since we
  don't need the even dimensional case in this paper, we leave this for
  further investigation. 
\end{remark}



\begin{proof}
In order to prove Theorem \ref{theo:odd_good_signature_perturbation}, we must
carefully recall the definition of the signature operator for  odd
dimensional manifolds. This is done in Appendix 
\ref{sec:appendix-odd-sign}.


  We denote the Hilbert $A$-module
  $L^2(M,\Lambda^* M\otimes f^* \mathcal{V}^\prime)$ by $\mathcal{L}^2
  (M)^*$.  Similarly we set $L^2(M^\prime,\Lambda^* M^\prime\otimes
  \mathcal{V}^\prime)=:\mathcal{L}^2 (M^\prime)^*$. 
We would like to compare $\mathcal{L}^2 (M)^*$
 with
$\mathcal{L}^2 (M^\prime)^*$ using the pull-back map $f^*$; but this
is in general not $L^2$-bounded.. As in Hilsum-Skandalis \cite[p. 90]{HiSka},
we modify $f^*$ in order to obtain a $L^2$-bounded
cochain map between
 $\mathcal{L}^2 (M^\prime)^*$
and
$\mathcal{L}^2 (M)^*$
as follows.
From \cite[p. 90]{HiSka}, for suitably large $N$,
there is a submersion
$ F\colon D^N \times M  \rightarrow M^\prime$ such that
$F(0, m) = f (m)$. Here $D^N$ is an open ball
in an Euclidean space of dimension $N$.
Fix $v \in \Omega^N_c(D^N)$ with $\int_{D^N} v = 1$.
Define a {\it bounded} cochain map
$$T:  \mathcal{L}^2 (M^\prime)^* \longrightarrow \mathcal{L}^2 (M)^*\,,\quad
\omega\rightarrow  \int_{D^N} v \wedge F^*(\omega).$$

\smallskip
Let $\D_{\mathcal{V}}$ and $\D_{\mathcal{V}'}^\prime$ be the 
signature operators on $M$ and $M^\prime$, as introduced above.
There is a
  well-defined functional calculus associated to these regular operators
 on
  Hilbert modules. In particular, if
  $\phi\colon \RR\to \RR$ is a rapidly decreasing (i.e.~Schwartz)
  function, then $\phi
  (\D_{\mathcal{V}})\in\Psi^{-\infty}_{A}$ (one can prove this e.g.~by
  using the heat kernel and
\begin{definition}\label{def:of_perturbation}
 We define
  \begin{equation*}
\C_{\phi,f}:=\phi(D_{\V}) \circ T \circ \phi(\D_{\mathcal{V}'}')\colon \L^2(M')^*\to
  \L^2(M)^*.
\end{equation*}
If  $\phi$ is an even function, $\phi(D_{\V})$  is a function of
the  Laplacian $\D_{\V}^2=\Delta$ and therefore preserves the degree of the
differential forms. Consequently, the same is true for $\C_{\phi,f}$.
\end{definition}

\begin{lemma}\label{lem:properties_of_T_eps}
  The operator  $\C_{\phi,f}$ is an integral $A$-linear operator
  with smooth kernel on $M\times  M^\prime $. Moreover, 
if $\phi$ is even then $d_{\V}\circ
  \C_{\phi,f} = \C_{\phi,f}\circ d'_{\V^\prime}$.
\end{lemma}

\begin{proof}
Since $T$ is bounded and $\phi
  (\D_{\mathcal{V}}),
\phi  (\D_{\V^\prime}^\prime)\in\Psi^{-\infty}_{A}$
we immediately get that  $\C_{\phi,f}$ is a $A$-linear smoothing operator. 
We know that $d_{\V} g^*\omega=g^*(d'_{\V^\prime}\omega)$ for any
smooth map $g$. Moreover, $d_{\V}$ commutes with
$(d_{\V}+d_{\V}^*)^2=\D_{\V}^2$ and therefore also with $\phi(\D_{\V})$
since $\phi$ is even; similarly
$d'_{\V'}$ commutes with $\phi(\D^\prime_{\V^\prime})$.
Since 
$$ \C_{\phi,f}(\omega)=\phi(\D_{\V}) \circ\left(
 \int_{D^N} v \wedge 
 F^* 
 \phi(\D_{\mathcal{V}'}')(\omega)\right) 
$$
we also get 
$d_{\V}\circ
  \C_{\phi,f} = \C_{\phi,f}\circ d'_{\V^\prime}$ as required, as
  explained in \cite[Proof of Theorem 3.3 on p.~90]{HiSka}, where we
  use that the form $v$ is closed. 
\end{proof}

\begin{definition}
  We use the usual inner product (coming from a fixed Riemannian
  metric) on the space of differential forms. We use $^*$ for the
  adjoint with respect to this inner product.

  The involution $\tau$ coming from the Hodge-$*$ operator (compare
  Appendix \ref{sec:appendix-odd-sign}) and the above inner product
  define the signature quadratic form on the space of differential
  forms, we use $^\dagger$ for the adjoint with respect to this
  quadratic form. Recall that 
  \begin{equation*}
    A^* = \tau A^\dagger\tau,\qquad A^\dagger = \tau A^*\tau.
  \end{equation*}
\end{definition}

Let now $\phi_\epsilon$ be a smooth even
function $\phi_\epsilon\colon\reals\to\reals$ such that $\phi_\epsilon(t)= 1$ if
$\abs{t}\le\epsilon/4$, and $\phi_\epsilon(t)=0$ if $\abs{t}\ge
\epsilon/2$. Recall
that $f\colon M\to M'$ was an orientation preserving homotopy
equivalence and write  
\begin{equation*}
T_\epsilon:=\C_{\phi_\epsilon,f}\colon \Omega^* (M',\V^\prime)\to \Omega^*(M,\V).
\end{equation*}
Clearly $T_\epsilon$ extends to an $L^2$-bounded $A$-linear operator.
Let $T^*_\epsilon$ be the adjoint of $T_\epsilon$ and define $T^\dag_\epsilon:= \tau^\prime T^*_\epsilon
\tau$, with $\tau$ and $\tau^\prime$ denoting the involutions defined
by the Hodge $\star$-operators on $M$ and $M^\prime$ (see Appendix
\ref{sec:appendix-odd-sign}).
Introduce a  new differential
$d$ by setting $d\alpha:= i^{|\alpha|} d_{\V} \alpha$ and similarly
for $d'$.

\begin{lemma}\label{compressed-chain}
The bounded operator $T_\epsilon$ verifies the following properties:\\
(a) $T_\epsilon ({\rm dom} (d))\subset {\rm dom} (d^\prime)$;
$T_\epsilon d^\prime = d T_\epsilon$\\
(b) $T_\epsilon$ induces an isomorphism in cohomology, with inverse
induced by $T_\epsilon^\dagger$\\
(c) there exists a bounded operator $y_\epsilon$ of degree $-1$ on $\L^2(M^\prime)^*$,
with $ y_\epsilon^\dagger =-y_\epsilon$ and 
such 
that $y_\epsilon({\rm dom} (d^\prime))\subset {\rm dom} (d^\prime)$ and
$1-T^\dag_\epsilon T_\epsilon= d^\prime y_\epsilon + y_\epsilon d^\prime $\\
(d) $y_\epsilon=y_1+y_2$ where $y_1$ is $\epsilon$-spectrally concentrated and $y_2$
commutes with $\Delta'$.\\ 
(e) $T_\epsilon$ is $\epsilon$-spectrally concentrated.
\end{lemma}
\begin{proof}
  The  results (a) to (c) for the operator $T$ are proved by
  Hilsum-Skandalis. 
We need first of all to extend them to $T_\epsilon$. Observe that
Lemma \ref{lem:properties_of_T_eps} is precisely assertion (a).
  Next, we observe that the $T$ of Hilsum-Skandalis and our
  $T_\epsilon$ are chain homotopic. Indeed, recall that $T_\epsilon =
  \phi_\epsilon(d_{\V}+d_{\V}^*)\circ T\circ
  \phi_\epsilon(d'_{\V'}+(d'_{\V'})^* )$. Here, we can replace $\D_{\V}$
  by $d_{\V}+d^*_{\V}$, since 
  $\D_{\V}^2=(d_{\V}+d^*_{\V})^2$, and $\phi$ is even. Choose a
  homotopy $\phi^t_\epsilon$ of even functions with $\phi^0_\epsilon=\phi_\epsilon$,
  $\phi^1_\epsilon=1$ and such that $\phi^t_\epsilon(x)=1$ for all $t\in[0,1]$ and
  $\abs{x}<\epsilon/4$. Then $\frac{d}{dt}\phi_t(x)= xg_t(x)=g_t(x)x$ with a
  suitable smooth odd family of functions $g_t\colon\reals\to\reals$. Observe that 
  \begin{equation}\label{eq:T_T_eps_homotopy}
    \begin{split}
      T-T_\epsilon = & \int_0^1 \frac{d}{dt} \left(\phi_t(d_{\V}+d^*_{\V})\circ
        T\circ \phi_t(d'_{\V}+(d'_{\V'})^*)\right)\,dt\\
      =&  (d_{\V}+d_{\V}^*)\int_0^1\left( g_t(d_{\V}+d^*_{\V}) \circ T\circ
        \phi_t(d'_{\V'} +(d'_{\V'})^*)\right)\,dt \\
      &+ \int_0^1 \left(\phi_t(d_{\V}+d^*_{\V})\circ T\circ
        g_t(d'_{\V'} +(d'_{\V'})^*)\right)\,dt\circ (d'_{\V'}
      +(d'_{\V'})^*)\\
      = & d_{\V} z+zd_{\V'}'        = dw+wd'
\end{split}
\end{equation}
\begin{equation*}
  \begin{split}
\text{with}\qquad\qquad    z= &\int_0^1 g_t(d_{\V}+d^*_{\V}) \circ T\circ
      \phi_t(d'_{\V'} +(d'_{\V'})^*)dt \\
    &+ \int_0^1 \phi_t(d_{\V}+d^*_{\V})\circ T\circ
        g_t(d'_{\V'} +(d'_{\V'})^*)dt ,\\
     w(\omega) =& i^{1-\abs{\omega}} z(\omega).
      \end{split}
\end{equation*}
The last but one equality in \eqref{eq:T_T_eps_homotopy} is true since
$d^*_{\V}g_t(d_{\V}+d^*_{\V})=g_t(d_{\V}+d^*_{\V})d$ and
$d_{\V}T=td'_{\V}$.

Since by \cite{HiSka} $T$ induces an isomorphism in homology, and
$T_\epsilon$ is chain homotopic to $T$, so does $T_\epsilon$, proving
(b).

Apply now $\cdot^\dag$ to Equation \eqref{eq:T_T_eps_homotopy} and use
that $d^\dag=-d$ to get
\begin{equation*}
  T^\dag-T_\epsilon^\dag = -w^\dag d' - dw^\dag.
\end{equation*}

Let now  $y$ be an operator such that $1-T^\dag T=dy+yd'$, as
constructed in \cite{HiSka}. Then
\begin{equation}\label{eq:t_eps_hom}
  \begin{split}
    1- T_\epsilon^\dag T_\epsilon =& 1-T^\dag T +
    (T^\dag-T_\epsilon^\dag)T + T_\epsilon^\dag(T-T_\epsilon) \\
    =& dy+yd' + (-dw^\dag-w^\dag d') T + T_\epsilon^\dag(dw+wd')\\
    =& d( y-w^\dag T+T_\epsilon^\dag w) + (y-w^\dag T+T_\epsilon^\dag
    w) d',
  \end{split}
\end{equation}
using that $T$ and $T_\epsilon$ are chain maps, and therefore in
particular $T^\dag_\epsilon d'=-T^\dag_\epsilon (d')^\dag=-d^\dag
T^\dag_\epsilon=dT^\dag_\epsilon$. Equation \eqref{eq:t_eps_hom} means
that we get (c) with $\tilde y_\epsilon= y-w^\dag T + T^\dag_\epsilon w$.

We will now modify $\tilde y_\epsilon$ in such a way that it splits as
required in (d). To do this, we make the following general
observation. Assume that $\psi_1\colon\reals\to\reals$ vanishes in a
neighborhood of $0$ and write $\psi_1(x)=x\psi_2(x)$ with a smooth
function $\psi_2$. Define 
\begin{equation}\label{eq:def_of_good_u}
u:=(d'_{\V'})^* \psi_2(\Delta') =\psi_2(\Delta')(d'_{\V'})^*.
\end{equation}
Then
\begin{equation}\label{eq:good_u}
  d'_{\V'}u + ud'_{\V'} = d'_{\V}(d'_{\V'})^* \psi_2(\Delta') +
  (d_{\V'}')^*d'_{\V'}\psi_2(\Delta') =\psi_1(\Delta').
\end{equation}

Choose $\psi\colon \reals\to\reals$ with support contained in
$[\epsilon/2,\infty)$ and such that $\psi(x)=1$ for $x\ge \epsilon$. 
Define now 
\begin{equation*}
y_1 := (1-\psi)(\Delta') \circ \tilde y_\epsilon \circ (1-\psi)(\Delta').
\end{equation*}
By construction, this operator is $\epsilon$-spectrally concentrated. Set $y_2:= \psi(\Delta')\tilde y_\epsilon (1-\psi)(\Delta')$,
$y_3:=\tilde y_\epsilon \psi(\Delta')$. We compute (since $T_\epsilon^\dag$ and $T_\epsilon$ are $\epsilon/2$-spectrally
concentrated)
\begin{equation*}
  \begin{split} d' y_2+y_2d' &=
    d' \circ (\psi(\Delta')\tilde y_\epsilon(1-\psi)(\Delta'))
     +(\psi(\Delta')\tilde y_\epsilon(1-\psi)(\Delta')) \circ d'\\
    = & \psi(\Delta') (d' \tilde y_\epsilon + \tilde y_\epsilon d')
    (1-\psi)\Delta')\\
    = &\psi(\Delta')\circ (1-\psi)(\Delta')\\
    =& d' \circ u_1 + u_1 d'
 \end{split}
\end{equation*}
 with a suitable operator $u_1$ defined as in \eqref{eq:def_of_good_u}
 and \eqref{eq:good_u} which commutes with $\Delta'$.
 Similarly,
 \begin{equation*}
   d'y_3+y_3d' = d'\circ \tilde y_\epsilon \psi(\Delta') + \tilde y_\epsilon
   \psi(\Delta')\circ d' = d' u_2 +u_2 d'
 \end{equation*}
 with a $u_2$ which also commutes with $\Delta'$.

Consequently,
\begin{equation*}
  \begin{split}
    1- T_\epsilon^\dag T_\epsilon =& d'\tilde y_\epsilon +\tilde
    y_\epsilon d'\\
    =&  d' (y_1 + y_2 + y_3) + (y_1+y_2 + y_3) d'\\
    = & d'(y_1+u_1+u_2) + (y_1+u_1+u_2)d'.
  \end{split}
\end{equation*}
We now set $y_\epsilon:=y_1+u_1+u_2$ and observe that $(u_1+u_2)$
commutes with $\Delta'$. This proves (d).\\
(e) is an immediate consequence of the construction of $T_\epsilon$.

\end{proof}

Following closely \cite{HiSka},
consider now the operators on $\L^2 (M')^*\oplus \L^2 (M)^*$:
\begin{equation}\label{key-oper}
R_\epsilon=\begin{pmatrix} 1 & 0
\cr - T_\epsilon\gamma   & 1
\end{pmatrix} \,,\quad L_{\epsilon,\alpha} = \begin{pmatrix}
  1-T_\epsilon^\dag T_\epsilon & (\gamma+\alpha y_\epsilon)T^\dag_\epsilon
\cr  T_\epsilon (-\gamma-\alpha y_\epsilon)   & 1
\end{pmatrix}\,,\quad \delta_{\epsilon,\alpha}=\begin{pmatrix}
  d^\prime & \alpha T^\dag_\epsilon
\cr 0   & -d
\end{pmatrix}
\end{equation}
with $\gamma\omega:= (-1)^{|\omega|}\omega$ and $\alpha$
a real number \footnote{The
  differences
in signs between these operators and the ones  appearing in \cite{HiSka}
are due to the fact that we take $\tau^\prime$ for the grading on
$M^\prime$,
whereas Hilsum-Skandalis take $-\tau^\prime$} . Note that $\gamma$ anticommutes
 with $d$ and $\tau$ and
commutes with $T$ and $T^\dagger$, and that $\gamma^\dagger=-\gamma$. 
It is
clear that $R_\epsilon$ and $L_{\epsilon,\alpha}$
are bounded. The crucial relation is
$L_{\epsilon,\alpha}\delta_{\epsilon,\alpha}= -
\delta_{\epsilon,\alpha}^\dagger L_{\epsilon,\alpha}$.

We notice that $R_\epsilon$ is invertible and that
$R_\epsilon^\dag R_\epsilon=L_{\epsilon,0}$.
Thus $L_{\epsilon,\alpha}$ is invertible for $\abs{\alpha}$ small enough.
Let $S_{\epsilon,\alpha}:=\frac{\tau \circ L_{\epsilon,\alpha}}{ |\tau \circ L_{\epsilon,\alpha}|}$,
with $\tau=\begin{pmatrix} \tau' & 0
\cr 0   & \tau
\end{pmatrix}$.
Then $S_{\epsilon,\alpha}$ is an involution. 
We now endow $\Omega^*(M',\V')\oplus \Omega^*(M,\V)$ with the new inner  product
$$<\omega_1,\omega_2>_{\epsilon,\alpha}:=<\omega_1,|\tau\circ L_{\epsilon,\alpha}|\omega_2>.$$
Notice that $|\tau \circ L_{\epsilon,\alpha}|$ is positive and self-adjoint with
respect to both  scalar-products.
Let $\widehat{\D}_{\epsilon,\alpha}:=-i (\delta_{\epsilon,\alpha}
S_{\epsilon,\alpha} +S_{\epsilon,\alpha}\delta_{\epsilon,\alpha})$ 
be the signature operator associated to
$\delta_{\epsilon,\alpha}$
and to the grading $S_{\epsilon,\alpha}$.
Using \cite[Lemme 2.1]{HiSka}, $\widehat{\D}_{\epsilon,\alpha}$
 is {\it invertible} and self-adjoint with respect to
 $<\cdot,\cdot>_{\epsilon,\alpha}$.
For the adjoint with respect to the original inner
product we therefore get:
 $$(\widehat{\D})_{\epsilon,\alpha}^*= 
=( |\tau\circ L_{\epsilon,\alpha}|)^{-1}\circ
\widehat{\D}_{\epsilon,\alpha} \circ ( |\tau\circ L_{\epsilon,\alpha}| )\,.$$

\begin{definition}\label{def:of_perturbed_D}
  We define the special perturbed signature operator
  \begin{equation*}
    \D_{\epsilon,\alpha} := -i \left( \delta_{\epsilon,\alpha}\tau
      L_{\epsilon,\alpha} + \frac{\tau
        L_{\epsilon,\alpha}}{\abs{\tau L_{\epsilon,\alpha}}}
      \delta_{\epsilon,\alpha} \frac{\tau
        L_{\epsilon,\alpha}}{\abs{\tau L_{\epsilon,\alpha}}} \tau L_{\epsilon,\alpha}\right)
  \end{equation*}
  with 
  \begin{equation*}
    \delta_{\epsilon,\alpha} = 
    \begin{pmatrix}
       d & \alpha T^\dagger_\epsilon\\ 0 & - d'
    \end{pmatrix}, \quad L_{\epsilon,\alpha} = 
    \begin{pmatrix}
      1- T_\epsilon^\dagger T_\epsilon & (\gamma+\alpha y_\epsilon)
      T^\dagger_\epsilon\\
      T_\epsilon (-\gamma -\alpha y_\epsilon) & 1 
    \end{pmatrix}
  \end{equation*}
\end{definition}

Then  $ {\D_{\epsilon,\alpha}}:= \widehat{\D_{\epsilon,\alpha}}\circ |\tau\circ L_{\epsilon,\alpha}|$ is self-adjoint with respect
the {\it original} inner product and is invertible for $\alpha\ne 0$
sufficiently small (depending on $\epsilon$). 
We complete the proof of the Theorem with the following

\begin{lemma}\label{lem:analysis_of_perturbations}
Consider the operators of Definition \ref{def:of_perturbed_D}. Then
$\tilde\delta_{\epsilon,\alpha}:=\delta_{\epsilon,\alpha} - 
\left(\begin{smallmatrix}
  d & 0 \\ 0 & -d'
\end{smallmatrix}\right)$ and $\tau L_{\epsilon,\alpha}- \tau$ are $\epsilon$-spectrally
concentrated. Consequently, $\frac{\tau L_{\epsilon,\alpha}}{\abs{\tau
    L_{\epsilon,\alpha}}}- \tau$ and $\D_{\epsilon,\alpha}-\D$ are
$\epsilon$-spectrally concentrated (recall that $\abs{\tau}=1$). 

Moreover, all these differences belong to $\Psi_A^{-\infty}$, in
particular $\D_{\epsilon,\alpha}-\D\in \Psi_A^{-\infty}$. 
\end{lemma}

\begin{proof}
First of all, using the definition of $L_{\epsilon,\alpha}$ and Lemma \ref{compressed-chain} we
claim  that $L_{\epsilon,\alpha}=1+\Theta_{\epsilon,\alpha}$ with 
$\Theta_{\epsilon,\alpha}\in\Psi^{-\infty}_A$
and $\epsilon$-spectrally concentrated. Indeed
\begin{equation*}
\begin{split}
L_{\epsilon,\alpha} &=1+
\begin{pmatrix} -T_\epsilon^\dag T_\epsilon & i\gamma T^\dag_\epsilon
\cr  -T_\epsilon \gamma   & 0
\end{pmatrix}+ \begin{pmatrix} 0 & \alpha y T^\dag_\epsilon
\cr  -T_\epsilon \alpha y   & 0
\end{pmatrix}\\&=
1+\Theta_{\epsilon, \alpha}^1+\begin{pmatrix} 0 & \alpha y_1 T^\dag_\epsilon
\cr - T_\epsilon \alpha y_1   & 0
\end{pmatrix} + \begin{pmatrix} 0 & \alpha y_2 T^\dag_\epsilon
\cr  -T_\epsilon \alpha y_2   & 0
\end{pmatrix}\\&=1+\Theta_{\epsilon,\alpha}^1+\Theta_{\epsilon,\alpha}^2+
+ \begin{pmatrix} 0 & \alpha y_2 T^\dag_\epsilon
\cr  -T_\epsilon \alpha y_2   & 0
\end{pmatrix}
\end{split}
\end{equation*}
with $\Theta_{\epsilon,\alpha}^j$, $j=1,2$ smoothing and $\epsilon$-spectrally concentrated
because of lemma \ref{compressed-chain}.
It remains to be checked that   the last term appearing in
the above formula is also smoothing and $\epsilon$-spectrally concentrated: but
since $y_2$ commutes with $\Delta'$, it also commutes with
$\phi(\D'_{\V'})$ if $\phi$ is even; using the very definition of
$T_\epsilon$ we then get that 
\begin{multline*}
  \Theta^3_{\epsilon,\alpha}=
\begin{pmatrix}
  0 & \alpha \phi_\epsilon(D_\V)\circ y_{2,\epsilon}\circ T\circ
  \phi_\epsilon(D'_{\V'})\\
 - \alpha \phi_\epsilon(D_\V) T\circ \circ y_{2,\epsilon}\circ
  \phi_\epsilon(D'_{\V'}) & 0
\end{pmatrix}
\end{multline*}
is smoothing and $\epsilon$-spectrally concentrated. Our claim now
follows with $\Theta_{\epsilon,\alpha}=\Theta_{\epsilon,\alpha}^1+\Theta_{\epsilon,\alpha}^2+\Theta_{\epsilon,\alpha}^3$.\\
 Thus $\tau \circ
L_{\epsilon,\alpha}=\tau+\Phi_{\epsilon,\alpha}$, with $\Phi_{\epsilon,\alpha}\in\Psi^{-\infty}_A$
and $\epsilon$-spectrally concentrated. 
Let us now consider
$|\tau \circ
L_{\epsilon,\alpha}|=\sqrt{(\tau \circ
L_{\epsilon,\alpha})^2}$.
We can write 
$$|\tau \circ
L_{\epsilon,\alpha}|=\frac{i}{2\pi} \int_{\C} \lambda^{\frac{1}{ 2}}((\tau \circ
L_{\epsilon,\alpha})^2-\lambda)^{-1} d\lambda
$$
with $\C$ equal to a circle of radius larger than the norm of $(\tau \circ
L_{\epsilon,\alpha})^2$. From this integral representation it follows that
$|\tau \circ
L_{\epsilon,\alpha}|=1+\Lambda_{\epsilon,\alpha}$ with 
$\Lambda_{\epsilon,\alpha}\in\Psi^{-\infty}_A$
and $\epsilon$-spectrally concentrated; indeed, 
$(\tau \circ
L_{\epsilon,\alpha})^2-\lambda$ is equal to $(1-\lambda)$  plus a $\epsilon$-spectrally concentrated
smoothing operator. Thus by Lemma \ref{lem:inverse_of_one_plus_smoothing},   its inverse will be of  type
$(1-\lambda)^{-1} (1+\Xi_{\epsilon,\alpha}(\lambda))$ with $\Xi_{\epsilon,\alpha}
(\lambda)\in \Psi^{-\infty}_A$ and $\epsilon$-spectrally concentrated; carrying
out the integral we obtain immediately that $|\tau \circ
L_{\epsilon,\alpha}|=1+\mathcal{F}_{\epsilon,\alpha}$ with 
$\mathcal{F}_{\epsilon,\alpha}\in\Psi^{-\infty}_A$
and $\epsilon$-spectrally concentrated.
It then follows that the same is true for $|\tau \circ L_{\epsilon,\alpha}|^{-1}$.
Thus for $S_{\epsilon,\alpha}$ we obtain: $S_{\epsilon,\alpha}=\tau+ \mathcal{G}_{\epsilon,\alpha}$
with $\mathcal{G}_{\epsilon,\alpha}$ smoothing and $\epsilon$-spectrally concentrated.
Consider now $\widetilde{\D}$; this equals
$-i (\delta_{\epsilon,\alpha} \circ S_{\epsilon,\alpha} + S_{\epsilon,\alpha} \circ \delta_{\epsilon,\alpha})
|\tau \circ L_{\epsilon,\alpha}|$.
By its very definition
and Lemma \ref{compressed-chain}, $\delta_{\epsilon,\alpha}=d+ \E_{\epsilon,\alpha}$
with $\E_{\epsilon,\alpha}$ smoothing and $\epsilon$-spectrally concentrated.
Thus
$${\D}_{\epsilon,\alpha}= -i \left( (d+\E_{\epsilon,\alpha})\circ
  (\tau+\mathcal{G}_{\epsilon,\alpha})+(\tau+\mathcal{G}_{\epsilon,\alpha})
\circ (d+\E_{\epsilon,\alpha}) \right) \circ (1+ \mathcal{F}_{\epsilon,\alpha}),$$
which is equal to $\D$ plus a smoothing operator $\epsilon$-spectrally
concentrated.
\end{proof}
Theorem \ref{theo:odd_good_signature_perturbation} is proved.
\end{proof}

\section{From stable to unstable eta-invariants}\label{sec:stable-=-unstable}

We continue with the setting of Section \ref{sec:homot-invar-sign}.
In this section, we analyze the $\tau$-$\eta$-invariant of the special
perturbations $\D_{\epsilon,\alpha}$ of $\D$, where we use the
notation and conventions of Section \ref{sec:unst-rho-invar}. In
particular, $A$ is a von Neumann algebra and $\tau\colon A\to Z$ a
positive normal trace with values in the commutative von Neumann
algebra $Z$.

\subsection{Limits of eta invariants}\label{subsec:limits}

\begin{theorem}\label{theo:perturbed_to_unperturbed_eta}
  Let $f\colon M'\to M$ be a smooth oriented homotopy equivalence
  between closed Riemannian oriented manifolds as in Section \ref{sec:homot-invar-sign}.
  This gives rise to twisted signature operators $D_{\L}$ on $M$ and
  $D'_{\L'}$ on $M'$, with $\L':=f^* \L$.
Consider the (unperturbed) $\tau$-eta invariants
$ \eta_\tau(D_\L)$, $\eta_\tau(D'_{\L'}) \in Z$, see
\ref{def:tau-eta}.

Assume, in addition to the above, that the von Neumann algebra $A$ admits a positive faithful
normal trace $\tau_A\colon A\to\complexs$.

  There are sequences $\epsilon_k>0$  and  $\alpha_k>0$ such that
  $\alpha_k$ is small enough for $\D_{\epsilon_k,\pm\alpha_k}$ to be
  an invertible perturbation of $\D$ and such that
  \begin{equation*}
    \lim_{k\to\infty} \frac{\eta_{\tau}(\D_{\epsilon_k,\alpha_k})
      +\eta_\tau (\D_{\epsilon_k,-\alpha_k})}{2} = \eta_\tau(\D).
  \end{equation*}
  Here, $\D_{\epsilon,\alpha}$ is the special smoothing perturbation of the signature
  operator $\D$ on the manifold $M\disjointunion(-M')$ as defined in
  Definition \ref{def:of_perturbed_D}. Note that on the left hand side we
  have a sequences of averaged perturbed eta-invariants, converging to the
  unperturbed eta-invariant on the right hand side.
\end{theorem}

\begin{proof}
We shall prove that
 there are sequences $\epsilon_k>0$  and  $\alpha_k>0$ such that
  \begin{equation*}
    \lim_{k\to\infty} \eta_{\tau}(\D_{\epsilon_k,\alpha_k})=\eta_\tau
    (\D)+\sigma\;,\quad
     \lim_{k\to\infty}  \eta_\tau (\D_{\epsilon_k,-\alpha_k})=
     \eta_\tau(\D) -\sigma\;,\quad \sigma\in\RR
  \end{equation*}
Let $p_\epsilon:= \chi_{[-\epsilon,\epsilon]} (\Delta)$, where
$\chi_Y$ is the characteristic function of a set $Y$.

Since $\D_{\epsilon,\alpha}$ is an $\epsilon$-perturbation of $\D$, we
get 
\begin{equation*}
\D_{\epsilon,\alpha}=p_\epsilon\D_{\epsilon,\alpha}p_\epsilon +
(1-p_\epsilon)\D_{\epsilon,\alpha}(1-p_\epsilon),\qquad
(1-p_\epsilon)\D_{\epsilon,\alpha}(1-p_\epsilon) = (1-p_\epsilon)\D(1-p_{\epsilon}).
\end{equation*}

Consequently, by the spectral theorem and the definition of
$\eta_\tau$
\begin{equation}\label{eq:split_eta_spectrally}
  \eta_\tau(\D_{\epsilon,\alpha})=\eta_\tau(p_\epsilon
  \D_{\epsilon,\alpha}p_\epsilon) + \eta_\tau((1-p_\epsilon)\D(1-p_\epsilon)).
\end{equation}

To analyze $\eta_\tau$ of the compressed operators, we will make use
of the following translation into simple functional calculus.
\begin{lemma}\label{lem:translate_eta}
  Let $p$ and $B$ be self-adjoint bounded Hilbert $A$-module morphisms on a
  Hilbert $A$-module $H$, $p$ a projection, such that
  $\tau(p)<\infty$ and such that $B=pBp$. Let $q\colon \reals\to \reals$ be defined by
  $q(x)=1$ for $x>0$, $q(0)=0$ and $q(x)=-1$ for $x<0$.
  We define
  \begin{equation}\label{eq:def_of_sgn}
    \sgn_\tau(B):= \tau(q(B)).\qquad\text{Then }    \eta_\tau(B) = \sgn_\tau(B).
  \end{equation}
\end{lemma}
\begin{proof}
  We simply observe that 
  \begin{equation*}
q(x) = \frac{1}{\sqrt{\pi}} \int_0^\infty
x\exp(-tx^2)\,\frac{dt}{\sqrt{t}} = \frac{2}{\pi}\int_0^\infty x \exp(-t^2x^2)\,dt,
\end{equation*}
and $\abs{q(x)}\le 1$, therefore $q(B)= \frac{1}{\pi}
\int_0^\infty B\exp(-tB^2)\,dt/\sqrt{t}$. Here we observe that
$\int_0^T B\exp(-t^2B^2)\,dt$ is defined for each $T\ge 0$ since $B$
is bounded and is equal to the functional calculus of $B$ applied to
$\int_0^T x\exp(-t^2x^2)\,dt$. The latter family of functions,
depending on the parameter $T$, converges
pointwise to $q$, therefore the limit $\int_0^\infty B\exp(-t^2B^2)\,dt$
exists in the von Neumann algebra of Hilbert $A$-module morphisms on $H$ and
equals $q(B)$. Finally observe that we can throughout write $pBp$
instead of $B$. By assumption $p$ is of $\tau$-trace class. Since
those operators form an ideal, it follows that both  $pBp$ and
$pBp\exp(-t^2B^2)$ are $\tau$-trace class. Normality implies that for a strongly convergent
sequence $X_n$ of operators, $\tau(pX_n)\xrightarrow{n\to\infty}
\tau(pX)$. We can now interchange in the
above argument integration, passage to limit and application of $\tau$,
and the statement follows.
\end{proof}

\begin{lemma}\label{lem:large_t_eta_contrib}
  \begin{equation*}
  \lim_{\epsilon\to 0} \eta_\tau((1-p_\epsilon)\D(1-p_\epsilon)) = \eta_\tau(\D),
\end{equation*}
where existence of the limit is part of the statement.
\end{lemma}
\begin{proof}
  We have discussed in Section \ref{sec:unst-rho-invar} that each of
  the integrals defining the eta-invariants exists in the above
  statement. Moreover, for $0<\epsilon<1$,
  \begin{equation*}
    \eta_\tau((1-p_\epsilon)\D(1-p_\epsilon)) =
    \eta_\tau((1-p_1)\D(1-p_1)) +\eta_\tau((p_1-p_\epsilon)\D(p_1-p_\epsilon)).
  \end{equation*}
  Finally observe that $\tau(p_1)\le \tau(5\exp(-\D)^2)<\infty$ and
  that
  $p_1(p_1-p_\epsilon)\D(p_1-p_\epsilon)p_1=
  (p_1-p_\epsilon)\D(p_1-p_\epsilon)$, so that Lemma
  \ref{lem:translate_eta} implies that
  $\eta_\tau((p_1-p_\epsilon)\D(p_1-p_\epsilon))=
  \sgn_\tau((p_1-p_\epsilon)\D(p_1-p_\epsilon))$. Finally,
  $p_\epsilon$ strongly converges to $p_0$, so that
  $q((p_1-p_\epsilon)\D(p_1-p_\epsilon))$ strongly converges to
  $q((p_1-p_0)\D(p_1-p_0))$. By normality of $\tau$, therefore
  \begin{equation*}
    \sgn_\tau((p_1-p_\epsilon)\D(p_1-p_\epsilon))
    =\tau(q((p_1-p_\epsilon)\D(p_1-p_\epsilon)))
    \xrightarrow{\epsilon\to 0} \sgn_\tau((p_1-p_0)\D(p_1-p_0))
  \end{equation*}
  Finally, by definition and Lemma \ref{lem:translate_eta}
  \begin{equation*}
    \eta_\tau(\D)= \eta_\tau((1-p_1)\D(1-p_1)) + \sgn_\tau((p_1-p_0)\D(p_1-p_0)).
  \end{equation*}
\end{proof}

\begin{lemma}\label{lem:eta_conv}
  There are sequences $\epsilon_k>0$ and $\alpha_k>0$ such that
  $\alpha_k$ is small enough for $\D_{\epsilon_k,\alpha_k}$ to be
  smoothing invertible perturbations of $\D$ and such that
  \begin{equation*}
    \begin{split}
      \sgn_\tau(p_{\epsilon_k} \D_{\epsilon_k,\pm\alpha_k}p_{\epsilon_k})
   &   \xrightarrow{k\to\infty} \sgn_\tau(\D_{0,\pm 1})\\
  \end{split}
\end{equation*}
  Here 
  \begin{equation*}
    \D_{0,\alpha} := -i p_0\left( \delta_{0,\alpha}\tau
      L_{0} + \frac{\tau
        L_{0}}{\abs{\tau L_{0}}}
      \delta_{0,\alpha} \frac{\tau
        L_{0}}{\abs{\tau L_{0}}} \tau L_{0}\right) p_0
  \end{equation*}
  with 
  \begin{equation*}
    \delta_{0,\alpha} = p_0
    \begin{pmatrix}
       0 & \alpha T^\dagger_0\\ 0 & 0
    \end{pmatrix}p_0,
 \quad L_{0} =
    \begin{pmatrix}
      1- T_0^\dagger T_0 & \gamma
      T^\dagger_0\\
      -T_0 \gamma  & 1 
    \end{pmatrix},
\qquad T_0=p_0Tp_0.
  \end{equation*}
  In particular, $\D_{0,\alpha}=\alpha \D_{0,1}$ and therefore
  \begin{equation}\label{eq:sgn_change_of_oimit_eta}
    \sgn_\tau(\D_{0,1})= - \sgn_{\tau}(\D_{0,-1}).
  \end{equation}
\end{lemma}

Assuming Lemma \ref{lem:eta_conv}, we can now finish the proof of Theorem
\ref{theo:perturbed_to_unperturbed_eta}. Observe that 
\begin{equation*}
  \eta_\tau(\D_{\epsilon_k,\pm\alpha_k}) =
 \underbrace{
   \eta_\tau((1-p_{\epsilon_k})\D_{\epsilon_k,\pm\alpha_k}
   (1-p_{\epsilon_k}))}_{=\eta_\tau(1-p_{\epsilon_k})\D(1-p_{\epsilon_k})\xrightarrow{k\to\infty} \eta_\tau(\D)} +
 \underbrace{\eta_\tau(p_{\epsilon_k}\D_{\epsilon_k,\pm\alpha_k}
   p_{\epsilon_k})}_{\xrightarrow{k\to\infty} \eta_\tau(\D_{0,
     \pm1})
}
\end{equation*}
where we use Lemma \ref{lem:large_t_eta_contrib} for the first
convergence statement, and Lemma \ref{lem:eta_conv} for the
second. Averaging the two resulting equations for $\alpha_k$ and
$-\alpha_k$ and using Equation \eqref{eq:sgn_change_of_oimit_eta}
gives the statement of the Theorem.
\end{proof}

\begin{remark}\label{rem:exact_perturbed_eta_calc}
  We leave it to the reader to check directly that
  $\eta_\tau(\D_{0,1})=0$. Using the same proof as above, this gives
  the following stronger version of Theorem
  \ref{theo:perturbed_to_unperturbed_eta}:
  \begin{equation*}
    \lim_{k\to\infty} \eta_\tau(\D_{\epsilon_k,\alpha_k}) = \eta_\tau(\D).
  \end{equation*}
\end{remark}

\subsection{Proof of Lemma \ref{lem:eta_conv}}
\label{sec:proof-lemma-}

  Define for $\alpha\ne 0$ 
  \begin{equation*}
 \C_{\epsilon,\alpha}:= \abs{\alpha}^{-1}
  p_\epsilon\D_{\epsilon,\alpha}p_\epsilon
\end{equation*}
and observe  that $    \sgn_\tau(\C_{\epsilon,\alpha})=\sgn_\tau(p_\epsilon \D_{\epsilon,\alpha}p_\epsilon)$.
  Set for $\alpha\ne 0$
  \begin{equation*}
     \widetilde{ \C}_{\epsilon,\alpha}
 := -i \abs{\alpha}^{-1} p_\epsilon\left( \tilde\delta_{\epsilon,\alpha}\tau
      L_{\epsilon,\alpha} + \frac{\tau
        L_{\epsilon,\alpha}}{\abs{\tau L_{\epsilon,\alpha}}}
      \tilde\delta_{\epsilon,\alpha} \frac{\tau
        L_{\epsilon,\alpha}}{\abs{\tau L_{\epsilon,\alpha}}} \tau L_{\epsilon,\alpha}\right)p_\epsilon
  \end{equation*}
  with 
  \begin{equation*}
    \tilde\delta_{\epsilon,\alpha} = 
    \begin{pmatrix}
       0 & \alpha T^\dagger_\epsilon\\ 0 & 0
    \end{pmatrix}, \quad L_{\epsilon,\alpha} = 
    \begin{pmatrix}
      1- T_\epsilon^\dagger T_\epsilon & (\gamma+\alpha y_\epsilon)
      T^\dagger_\epsilon\\
      T_\epsilon (-\gamma -\alpha y_\epsilon) & 1 
    \end{pmatrix}.
  \end{equation*}

 Note that $\widetilde{\C}_{\epsilon,\alpha}$ is in
  general neither self adjoint nor invertible.

  \begin{lemma}\label{lem:good_alphas}
 There is a constant $K>0$ such that the following holds.
    For each $\epsilon>0$ there exists $\alpha_\epsilon>0$ such that
    \begin{equation*}
       \norm{\alpha_\epsilon
         y_\epsilon}<\epsilon,\qquad\norm{L_{\epsilon,\pm\alpha_\epsilon}^{-1}}\le K.
    \end{equation*}
  \end{lemma}
  \begin{proof}
    We first observe that there is a uniform bound for
    $\norm{T_\epsilon}$. This implies a uniform bound for
    $\norm{L_{\epsilon,0}}$ and for $\norm{L_{\epsilon,0}^{-1}}$, using
      $L_{\epsilon,0}= \left(
        \begin{smallmatrix}
          1 & \gamma T_\epsilon^\dagger \\ 0 & 1
        \end{smallmatrix}\right)
\left(
        \begin{smallmatrix}
          1 & 0 \\ -T_\epsilon\gamma & 1
        \end{smallmatrix}\right)$, so that $(L_{\epsilon,0})^{-1}= \left(
        \begin{smallmatrix}
          1 & 0 \\ T_\epsilon\gamma & 1
        \end{smallmatrix}\right)\left(
        \begin{smallmatrix}
          1 & - \gamma T_\epsilon^\dagger \\ 0 & 1
        \end{smallmatrix}\right)$.
      Since each $y_\epsilon$ is bounded and $L_{\epsilon,\alpha}=
      L_{\epsilon,0}+\alpha \left(
        \begin{smallmatrix}
          0 & y_\epsilon T^\dagger_\epsilon\\ -T_\epsilon y_\epsilon
          & 0
        \end{smallmatrix}\right)$
, we can now choose
      $\alpha_\epsilon$ such that all the assertions are fulfilled.
  \end{proof}

  \begin{lemma}\label{lem:tilde_C_convergnece}
    Choose a sequence $\epsilon_k>0$ with
    $\epsilon_k\xrightarrow{k\to\infty} 0$, and choose
    $\alpha_k:=\alpha_{\epsilon_k}>0$ as given by Sublemma
    \ref{lem:good_alphas}. Then $\widetilde{\C}_{\epsilon_k,\pm\alpha_k}$
    strongly converges to $\D_{0,\pm 1}$.
  \end{lemma}
  \begin{proof} For $\epsilon\to 0$, $p_\epsilon$ strongly converges
    to $p_0$,
  $T_\epsilon$ strongly converges to $T_0$, and $T_\epsilon^\dagger$
  strongly converges to $T_0^\dagger$.  $\norm{\alpha_\epsilon
    y_\epsilon}\xrightarrow{\epsilon\to 0} 0$. Since products of
  strongly convergent {sequences} strongly converge to the
  product of the limits, by its construction $L_{\epsilon_k,\alpha_k}$
  strongly converges to $ L_{0} $. Moreover, since
  $\norm{(\tau L_{\epsilon_k,\alpha_k})^{-1}}\le K$ independent of $k\in
  \naturals$, $\frac{\tau L_{\epsilon_k,\alpha_k}}{\abs{\tau L_{\epsilon_k,\alpha_k}}}=f(\tau L_{\epsilon_k,\alpha_k})$
  (and $\frac{\tau L_{0,0}}{\abs{\tau L_{0,0}}}=f(\tau L_{0,0})$) for
  any continuous function $f$ with $f(x)=-1$ for $x\le -1/K$, $f(x)=1$
  for $x\ge 1/K$. Since bounded continuous functions of strongly
  convergent sequences strongly converge by \cite[Theorem
  VIII.20]{MR751959}, $\frac{\tau
    L_{\epsilon_k,\alpha_k}}{\abs{L_{\epsilon_k,\alpha_k}}}$ strongly
  converges to $\frac{\tau L_{0,0}}{\abs{\tau L_{0,0}}}$. Putting
  everything together, and using in particular that
  $\abs{\alpha}^{-1}\tilde \delta_{\epsilon,\alpha}=\abs{\alpha}^{-1}
  \left(
    \begin{smallmatrix}
      0 & \alpha T_\epsilon^\dagger \\ 0 & 0 
    \end{smallmatrix}\right) =  \left(
    \begin{smallmatrix}
      0 & \pm T_\epsilon^\dagger \\ 0  & 0
    \end{smallmatrix}\right) $ ,
the result follows.
  \end{proof}

  \begin{lemma}\label{lem:non_tilde_C_convergence}
    There is a sequence $\epsilon_k>0$, with corresponding $\alpha_k$
    as in Sublemma \ref{lem:tilde_C_convergnece}, such that
    $\C_{\epsilon_k,\pm\alpha_k}$ converges in the strong resolvent
    sense to $\D_{0,\pm 1}$. Here, we will use the positive finite 
faithful 
    normal trace $\tau_A\colon A\to \complexs$.
  \end{lemma}
  \begin{proof}
    We will construct a sequence of monotonously increasing
    projections $Q_1\le Q_2\le \cdots$ with $\sup_k Q_k=p_1$. The
    latter property will be ensured showing  that for all $k$ $Q_k\le p_1$,
    and that $\tau_A(Q_k)\xrightarrow{k\to\infty}\tau_A(p_1)<\infty$,
    therefore $\tau_A(p_1-(\sup_{k}Q_k))=0$, so that $p_1=\sup_k Q_k$, using
      faithfulness and normality of $\tau_A$.

   The main property of the $Q_k$ will be that 
   \begin{equation}
\tilde   
\C_{\epsilon_n,\pm\alpha_n} Q_k = \C_{\epsilon_n,\pm\alpha_n}Q_k\qquad\text{for
   all }n\ge k.\label{eq:very_good_re}
 \end{equation}
   The assertion of the Sublemma follows then from Sublemma
   \ref{lem:tilde_C_convergnece}. To see this, we can represent $A$ as
   bounded operators on a Hilbert space $H$. Since $\tilde
   \C_{\epsilon_k,\pm\alpha_k}(1-p_1)=
   \C_{\epsilon_k,\pm\alpha_k}(1-p_1)=  \D_{0,\pm1}(1-p_1)=0$, it
   suffices to study the restriction to $\im(p_1)$. Now all operators
   in question are bounded, therefore $\bigcup_{k\in\naturals}
   (\im(Q_k))$ will be a common core for all of them (since
   $\sup_{k}Q_k=p_1$, $\bigcup_{k\in \naturals} \im(Q_k)$ is dense in
   $\im(p_1)$). A reference for the lattice of projections and
   properties of it we are using here is
   \cite{MR1468229,MR98f:46001b}, in particular Section 2.5.
Finally, for each $v\in \bigcup_{k}\im(Q_k)$ there is
   a $k\in\naturals$ such that $v\in \im(Q_k)$. Then, for $n\ge k$,
   $\C_{\epsilon_n,\pm\alpha_n}v = \tilde
   \C_{\epsilon_n,\pm\alpha_n}v\xrightarrow{n\to\infty} \D_{0,\pm
     1}v$. By \cite[Theorem VIII.25]{MR751959},
   $\C_{\epsilon_n,\pm\alpha_n}$ converges in the strong resolvent
   sense to $\D_{0,\pm 1}$..

   We now tackle the construction of the projections $Q_k$. 
   Choose $\epsilon_k$ such that $\tau_A(p_{\epsilon_k}-p_0)<10^{-k}$.

   Restrict now all operators in question to
   $\im(p_{\epsilon_k})$. Note that $A_k:=\tau
   L_{\epsilon_k,\alpha_{\epsilon_k}}$ now is an invertible operator,
   mapping $\im(p_{epsilon_k})$ to itself, and the same is true for
   $B_k:=\frac{L_{\epsilon_k,\alpha_{\epsilon_k}}}
   {\abs{L_{\epsilon_k,\alpha_{\epsilon_k}}}}
     L_{\epsilon_k,\alpha_{\epsilon_k}}$. Choose maximal projections
     $Q_{A_k}$ and $Q_{B_k}$ such that 
     \begin{equation}
p_0A_kQ_{A_k}=A_kQ_{A_k}\text{ and }
    p_0B_kQ_{B_k}=B_kQ_{B_k}.\label{eq:proj_rel}
   \end{equation}
In other words, $Q_{A_k}$ is the
     projection onto the inverse image under $A_k$ of $\im(p_0)$. Set
     $Q_{k}^0:=\inf\{p_0,Q_{A_k},Q_{B_k}\}$.

     Observe that, by construction, $\tilde
     \delta_{\epsilon,\alpha}p_0 = \delta_{\epsilon,\alpha}
     p_0$. Consequently, 
     \begin{equation}\label{eq:good_rel}
\tilde \C_{\epsilon_k,\alpha_k}Q_k^{0} =
     \C_{\epsilon_k,\alpha_{k}} Q_k^{0}\qquad\forall k\in\naturals.
   \end{equation}

     Finally, we set $Q_k:= (\inf_{l\ge k} Q^0_l) +
     p_1-p_{\epsilon_k}$. Observe that \eqref{eq:very_good_re}
     immediately follows from \eqref{eq:good_rel} and the fact that $\tilde \C_{\epsilon_k,\alpha_k}(p_1-p_{\epsilon_k})
    = \C_{\epsilon_k,\alpha_k}(p_1-p_{\epsilon_k})=0$. Note that
     $\tau_A(Q_{A_k})=\tau_A(p_0)=\tau_A(Q_{B_k})$. This follows since
     instead of $Q_{A_k}$ we can first construct the idempotent
     $e_k:=A_k^{-1}p_0A_k$; it satisfies the relation \eqref{eq:proj_rel}
     and $\tau_A(e_k)=\tau_A(p_0)$. Use then the standard
     construction of a selfadjoint projection $Q_{A_k}$ with
     $e_kQ_{A_k}=Q_{A_k}$ and $Q_{A_k}e_k=e_k$ (compare
     e.g.~\cite[Theorem 2.1]{MR1650345}), then
     $\tau_A(Q_{A_k})=\tau_A(e_kQ_{A_k})=\tau_A(Q_{A_k}e_k)=
     \tau_A(e_k)$, and $p_0A_kQ_{A_k}=p_0A_k e_kQ_{A_k}=A_ke_kQ_{A_k}=A_kQ_{A_k}$.

     Note that each of the projections $Q_{A_k},Q_{B_k},p_0$ are
     bounded above by $p_{\epsilon_k}$, with
     $\tau_A(Q_{A_k})=\tau_A(Q_{B_k})=\tau_A(p_0)$,
     $\tau_{A}(p_{\epsilon_k})<\tau(p_0)+10^{-k}$. The Kaplansky
     formula \cite[Theorem 6.1.7]{MR98f:46001b}
     $$p_0-\inf\{Q_{A_k},p_0\}=\sup\{p_0,Q_{A_k}\}-Q_{A_k}\le
     p_{\epsilon_k}-Q_{A_k},$$ 
linearity and positivity of $\tau_A$ imply $\tau_A(Q^0_k)\ge
     \tau_A(p_0)-2\cdot 10^{-k}$, and similarly, 
     \begin{equation*}
\tau_A(\inf_{l\ge k}
     Q^0_l) \ge \tau_A(p_0)- 2\cdot 10^{-k} \sum_{l=0}^\infty 10^{-l}
     \ge \tau_A(p_0) - 4 \cdot 10^{-k}.
   \end{equation*}
   It follows that
   \begin{equation*}
     \tau_A(Q_k) \ge \tau_A(p_0)- 4\cdot 10^{-k} +
     \tau_A(p_1-p_{\epsilon_k})\xrightarrow{k\to\infty}
     \tau_A(p_0)+\tau_A(p_1-p_0)= \tau_A(p_1).
   \end{equation*}
   Since $Q_k\le p_1$, the properties of $Q_k$ and therefore Sublemma
   \ref{lem:non_tilde_C_convergence} follow.
  \end{proof}

  \begin{sublemma}\label{lem:FA_str_res_conv}
    Assume that self adjoint operators $A_n$ converge in the strong
    resolvent sense to the self adjoint operator $A$. Assume that
    $\chi_{\{0\}}(A_n)$ strongly converges to $\chi_{\{0\}}(A)$. Then
      $q(A_n)$ strongly converges to $q(A)$, where $q=\chi_{(0,\infty)}-\chi_{(-\infty,0)}$.
  \end{sublemma}
  \begin{proof}
    By assumption, $\chi_{\{0\}}(A_n)+A_n $ converges in the strong
    resolvent
sense to  $\chi_{\{0\}}(A) +A$ (this is best seen using the Trotter
criterion \cite[Theorem VIII.21]{MR751959}). Thus by \cite[Theorem
  VIII.24]{MR751959} we know that $
\chi_{[0,\infty)}(\chi_{\{0\}}(A_n)+A_n)$ strongly converges to  $
\chi_{[0,\infty)}(\chi_{\{0\}}(A)+A)$ and similarly for
$\chi_{(-\infty,0]}$.
Hence $
\chi_{(0,\infty)}(A_n) + \chi_{\{0\}}(A_n)$ converges strongly to 
  $
\chi_{(0,\infty)}(A) + \chi_{\{0\}}(A)$ which implies that 
 $
\chi_{(0,\infty)}(A_n) $ converges strongly to  $
\chi_{(0,\infty)}(A).$ Similarly, $\chi_{(-\infty,0)} (A_n) $ converges strongly to  $
\chi_{(-\infty,0)}(A)$. The proof of the sublemma is complete.
\end{proof}

With all the Sublemmas in place, it is now easy to finish the proof of
Lemma \ref{lem:eta_conv}. We can restrict to $\im(p_1)$ (by
multiplication with $p_1$ from the left), since all operators in
question commute with $p_1$ and vanish on the complement of $\im(p_1)$
(i.e.~when multiplied with $(1-p_1)$). Observe that $\tau$, restricted to $\im(p_1)$, is
strongly continuous
(because $\tau(p_1)<\infty$).

Since $\D_{\epsilon,\alpha}$ is invertible and $\epsilon$-spectrally
concentrated near zero, $\chi_{\{0\}}(\C_{\epsilon,\pm \alpha}p_1) =
p_1-p_\epsilon$. In the same way, when restricted to $\im(p_0)$,
$\D_{0,\alpha}$ is invertible and therefore
$\chi_{\{0\}}(\D_{0,\alpha} p_1) = p_1-p_0$. Consequently,
$\chi_{\{0\}}(\C_{\epsilon_k,\pm\alpha_k})$ converges strongly to
  $\chi_{\{0\}}(\D_{0,\pm 1})$. Sublemmas
  \ref{lem:non_tilde_C_convergence} and \ref{lem:FA_str_res_conv}
  imply that $q(\C_{\epsilon_k,\pm\alpha_k})$ converge strongly to
  $q(\D_{0,1})$, so that finally, using Lemma \ref{lem:translate_eta},
  \begin{equation*}
    \eta_\tau(p_{\epsilon_k}\D_{\epsilon_k,\pm\alpha_k}p_{\epsilon_k})
    = \tau(q(\C_{\epsilon_k,\pm\alpha_k}))\xrightarrow{k\to\infty}
    \tau(q(\D_{0,\pm 1}))= 
    \eta_\tau(\D_{0,\pm 1}).
  \end{equation*}

\begin{remark}\label{rem:discrete_case_much easier}
  We want to observe that, if $0$ is isolated in the spectrum of $\D$,
  the proof of Lemma \ref{lem:eta_conv} is almost trivial. By
  construction, for sufficiently small $\epsilon$ in this case
  $\D_{\epsilon,\alpha}$ is equal to $\D_{0,\alpha}$.
\end{remark}

\section{Homotopy invariance of unstable rho-invariants}
\label{sec:homot-invar-unst}

We are now in the position to prove the statements concerning homotopy
invariance of the introduction, along with some more general results.

\begin{theorem}\label{theo:stable_rho_equal_unstalbe_rho_general_A}
  Let $M$ be a closed oriented Riemannian manifold of odd
  dimension. Assume that $\Gamma$ is a torsion free
  discrete group such that the maximal Baum-Connes map
\begin{equation*}
   \mu_{max}\colon K_*(B\Gamma)
\to K^*(C^*\Gamma)
  \end{equation*}
  is an isomorphism. Let $u\colon M\to B\Gamma$ be a map classifying a
  covering $\tilde M= u^* E\Gamma$. Let $A_j$ be unital von Neumann
  algebras admitting positive finite faithful normal traces
  $\tau_j\colon A_j\to\complexs$. Let
  $\alpha_j\colon C^*\Gamma\to A_j$ be homomorphisms (with
  $j=1,\dots,r$), and $\tau_1,\dots,\tau_r\colon
  A_j\to Z_j$ positive normal traces with
  values in abelian von Neumann algebras $Z_j$. Assume that
  $\beta_1,\dots,\beta_r\colon  Z_j\to V$ are continuous homomorphisms
  to a fixed topological
  vector space $V$ such that
  \begin{equation}\label{condition}
\sum_{j=1}^r \beta_j \tau_j(\alpha_j(1)) = 0 \in V.
\end{equation}

  Let $\L_j:= \tilde M\times_\Gamma A_j$ be the associated Hilbert
  $A$-module bundle, where $\Gamma$ acts on $A_j$ via $\Gamma\to
  C^*\Gamma\xrightarrow{\alpha_j} A_j$.

  Then $\rho_{(\tau_j,\beta_j,L_j)}(M) \in V$ is a homotopy invariant.
\end{theorem}

\begin{remark}
  We leave it to the reader to remove the hypothesis that the $A_j$ admit
  a finite faithful trace in Theorem \ref{theo:stable_rho_equal_unstalbe_rho_general_A}.
\end{remark}

\begin{proof}
  Let $f\colon M'\to M$ be a homotopy equivalence, $\L_j':= f^*\L_j$.
  Let $\mathcal{V}:=\tilde M\times_\Gamma C^*\Gamma$ be the
  Mishchenko-Fomenko line bundle associated
  to $u$ on $M$. Then $f^*\mathcal{V} =: \mathcal{V}'$ is the
  Mishchenko-Fomenko line bundle associated to $f\circ u$ on $M'$.

  Let $\D_{\epsilon_k,\pm\alpha_k}$ be the smoothing perturbation of
  $D_{\mathcal{V}}\disjointunion -\D'_{\mathcal{V}'}$ on
  $M\disjointunion (-M')$ as in Theorem \ref{theo:perturbed_to_unperturbed_eta}. Associated
  to these perturbation we define the stable rho-invariant
  \begin{equation*}
    \rho^s_{[0]}(M\disjointunion -M') = [\eta_{[0]}(\D_{\epsilon_k,\pm\alpha_k})] \in
    (C^*\Gamma)_{{\rm ab}}/<[1]> \,,
  \end{equation*}
where we recall that  $(C^*\Gamma)_{{\rm ab}}:=
(C^*\Gamma)/\overline{[C^*\Gamma,C^*\Gamma]}$.
  By Theorem \ref{theo:vanish-stable}, $\rho^s$ does not depend on the
  particular perturbation and vanishes identically.

  From Lemma \ref{lem:stable_eta_preserved}
  we conclude that
  \begin{equation}\label{eq:naturality_of_eta}
   \alpha_j(\eta_{[0]}(\D_{\epsilon_k,\pm\alpha_k})) =
    \eta_{[0]}(\D^{A_j}_{\epsilon_k,\pm\alpha_k}) \;\in A_j/\overline{\Commutator{A_j}{A_j}},
  \end{equation}
  where $\D^{A_j}_{\epsilon_k,\pm\alpha_k}$ is the  perturbation of the $A_j$-twisted
  signature operator $D_{\L_j}\disjointunion
  D'_{\L_j'}$ on $M\disjointunion (-M')$ given by the corresponding recipe. On the other hand, by Theorem
 \ref{theo:perturbed_to_unperturbed_eta}
  \begin{equation*}
    \frac{\eta_{\tau_j}(\D^{A_j}_{\epsilon_k,\alpha_k})
      +\eta_{\tau_j}(\D^{A_j}_{\epsilon_k,-\alpha_k})}{2}\xrightarrow{k\to \infty} \eta_{\tau_j}(\D^{A_j}) =  
    \eta_{\tau_j}(D_{\L_j})-\eta_{\tau_j}(D'_{\L_j'}).
  \end{equation*}

  Therefore
  \begin{equation}\label{eq:limit}
    \sum_{j=1}^r \beta_j \frac{\eta_{\tau_j}(\D^{A_j}_{\epsilon_k,\alpha_k})+ \eta_{\tau_j}(\D^{A_j}_{\epsilon_k,-\alpha_k})}{2}
    \xrightarrow{k\to\infty}    \sum_{j=1}^r \beta_j
    \eta_\tau(D_{\L_j})-\sum_{j=1}^r \beta_j \eta_\tau(D'_{\L_j'}).
  \end{equation}
  By Equation \eqref{eq:naturality_of_eta}, for the left hand side we
  can write
  \begin{equation*}
    \sum_{j=1}^r \beta_j \eta_{\tau_j}(\D^{A_j}_{\epsilon_k,\pm\alpha_k}) =
    \left(\sum_{j=1}^r \beta_j \tau_j\alpha_j\right)
    \left(\eta_{[0]}(\D_{\epsilon_k,\pm\alpha_k})\right)
  \end{equation*}
  The maps $ \tau_j\alpha_j\colon
  C^*\Gamma\to Z_j$ are traces. Therefore, they factorize through
$(C^*\Gamma)_{{\rm ab}}$. After composing with $\beta_j$ and summing,
we get that $\sum_{j=1}^r \beta_j\tau_j\alpha_j\colon
C^*\Gamma\to V$ maps $1$ to $0$, therefore this map
factorizes through $C^*_{\rm{ab}}/<[1]>$.

But the
  projection of $\eta_{[0]}(\D_{\epsilon_k,\pm\alpha_k})$ to
  $ (C^*\Gamma)_{{\rm ab}}/<[1]>$ is the stable
  rho-invariant of $M\disjointunion(-M')$ which  is identically
  zero.

  The assertion of the theorem follows now from \eqref{eq:limit}.
\end{proof}

\begin{corollary}\label{corol:proof_of_aps_rho_sign}
The statements about the signature operator in
Theorem
  \ref{theo:APS_rho}, Theorem
  \ref{theo:L2_rho} and Theorem \ref{theo:finite_conjugacy} are true.
\end{corollary}

\begin{proof}
It suffices to specialize the result of our main theorem to the
Examples \ref{example:finite-repres} and \ref{example:delocalized}.
Notice that in this way we have established the result about the delocalized eta invariant only
under the Baum-Connes assumption for the maximal $C^*$-algebra.
In order to sharpen this result to the reduced $C^*$-algebra
we simply observe that if in the statement of Theorem \ref{theo:stable_rho_equal_unstalbe_rho_general_A}
we take $j=1$, $A=\NeumannN\Gamma$, $\alpha:C^*_r \Gamma\to
\NeumannN\Gamma$ the natural map, $Z=\CC$, $\tau=\tau_{<g>}$
and $\beta_1=1$, then the basic condition
\eqref{condition} is satisfied, $\rho_{(\tau,\beta,\L)} (M)=\eta_{<g>}
(\tilde{D}^{{\rm sign}})$ and the proof carries over.
\end{proof}

\begin{remark}
  Let $M$ be a closed odd-dimensional manifold. We want to stress here
  that the mere vanishing of the signature index class in $K_1 (C^*
  \Gamma)$ does not imply the vanishing of the corresponding $L^2$-rho
  invariant: we wish to clarify this point.
 There certainly exist allowable perturbations to define the
  corresponding stable rho-invariant (by \cite{LP03})); under our
  standard assumptions on the fundamental group and the Baum-Connes
  map this stable rho invariant will be equal to zero.
However, in this generality it cannot be guaranteed that {\it
  stable=unstable}.
Indeed, there
  are examples where the index class is trivial but where
the  rho invariants are
  non-trivial (even for a fundamental groups as simple as $\integers$);
  we shall construct them
  in Section \ref{sec:examples-non-trivial}.
From this point
  of view, the case
  $M=X\sqcup (-X^\prime)$, with $X$ and $X^\prime$ homotopy equivalent,
  is indeed very  special.
\end{remark}

\section{Vanishing results
on spin manifolds with positive scalar curvature}\label{sec:vanishing-spin}

Let $(M,g)$ be a Riemannian manifold. We assume that $M$ is spin
and  that $\dim M=m=2k -1, k\geq 1$. We fix a spin structure and
we let $\Di$ be the associated Dirac operator.


\begin{theorem}\label{theo:unstable_rho_vanihs_general_A}
  Assume that $M$ has positive scalar curvature.
Under the same assumptions on $\Gamma$, $A_j$, $\alpha_j$, $\tau_j$
and $\beta_j$
as in Theorem \ref{theo:stable_rho_equal_unstalbe_rho_general_A} we have
 \begin{equation}\label{vanishing-spin}
\rho_{(\tau_j,\beta_j,L_j)}(M) =0 \in Z\,.
\end{equation}
\end{theorem}
\begin{proof}
  The proof is parallel to the proof of Theorem
  \ref{theo:stable_rho_equal_unstalbe_rho_general_A}, but much easier
  since under the assumption of positive scalar curvature the unstable
  rho-invariant coincides by definition with the stable
  one.
\end{proof}

\begin{corollary}\label{corol:proof_of_aps_rho_Dirac}
The statements about the spin Dirac operator in
Theorem
  \ref{theo:APS_rho}, Theorem
  \ref{theo:L2_rho} and Theorem \ref{theo:finite_conjugacy} are true.
\end{corollary}

\begin{proof}
  The proof is word by  word parallel to the proof of Corollary
  \ref{corol:proof_of_aps_rho_sign}, with the same remark applying in
  order
to pass from $C^* \Gamma$ to $C^*_r \Gamma$ for the delocalized eta invariant.
\end{proof}

\begin{remark}\label{remark:invisible}
Let $M$ be a closed oriented odd dimensional manifold and let
$\Gamma\to \tilde{M}\to M$ be a Galois covering, with $\Gamma$
torsion-free. Let us assume that the signature operator on
$\tilde{M}$ is $L^2$-invertible; thanks to the work of
Farber-Weinberger \cite{MR2003b:58044} and Higson-Roe-Schick
\cite{MR2002m:57031} we know that there are plenty of examples
of such coverings. Then our arguments
imply that $$\rho_{\lambda_1 - \lambda_2} (D^{{\rm
    sign}})=0\quad\text{and}\quad
\rho_{(2)} (\tilde{D}^{{\rm sign}})=0$$ provided the
max-Baum-Connes map is bijective and that $$\eta_{<g>}
(\tilde{D}^{{\rm sign}})=0$$ provided $<g>$ is of polynomial
growth and the reduced-Baum-Connes map is bijective. {\it
Conclusion}: for such $L^2$-invisible coverings our
signature-rho-invariants are zero, precisely as the spin
rho-invariants in the presence of positive scalar curvature.
\end{remark}

\section{Further remarks on delocalized invariants}\label{sec:remarks}


\subsection{Infinite conjugacy classes}
\label{sec:infin-conj-class}

Let $<g>$ be a non-trivial conjugacy class in $\Gamma$, not necessarily
of finite cardinality.

\begin{definition}\label{def:delocalized_eta_general}
If the limit
\begin{equation}\label{eq:delocalized_eta-truncated}
 \lim_{T\to +\infty} \frac{1}{\sqrt{\pi}} \int_0^T t^{-1/2}
  \Tr_{<g>}(\tilde D\exp(-t\tilde D^2) dt.
\end{equation}
exists, then its value is, by definition, the delocalized
eta invariant $\eta_{<g>}(\tilde D)$ associated to $\tilde D$ and $<g>$.
\end{definition}

Lott \cite{MR2000k:58039} \cite{lott-erratum} establishes the
convergence  if the operator $\tilde D$ on $\tilde M$ has a gap near zero in its spectrum
(i.e.~zero is isolated in the spectrum) and if
the trace $\tau_{<g>}:\CC\Gamma\to\CC$ associated to $<g>$ has a certain
extension-property.
Let us recall his results. Let $\B^\infty_\Gamma$, $\CC\Gamma\subset \B^\infty_\Gamma \subset C^*_r \Gamma$,
be the Connes-Moscovici algebra \cite{CM} and let $\D^\infty$ denote the Dirac operator
on $M$ twisted by the flat bundle $\tilde M\times_\Gamma \B^\infty$.
First of all, if $\tilde D$ has a gap in its spectrum,
then the integral defining the $0$-degree eta-invariant of $\D^\infty$
\begin{equation}\label{eta-cm}
\eta_{[0]} (\D^\infty):=\frac{1}{\sqrt{\pi}} \int_0^\infty t^{-1/2}
\TR (\D^\infty \exp(-t (\D^\infty)^2) dt \;\;\in\;\; \B^\infty_\Gamma/
\overline{\Commutator{\B^\infty_\Gamma}{\B^\infty_\Gamma}}
\end{equation}
converges in the natural Fr{\'e}chet topology induced by $\B^\infty_\Gamma$.
We set as usual $(\B^\infty_\Gamma)_{{\rm ab}}:=\B^\infty_\Gamma/
\overline{\Commutator{\B^\infty_\Gamma}{\B^\infty_\Gamma}}$.
Next, if
the trace $\tau_{<g>}:\CC\Gamma\to\CC$ {\it extends} from $\CC\Gamma$ to
a continuous trace on $\B^\infty_\Gamma$,
then one can prove (see Remark \ref{rem:a-la-Lott} in Appendix \ref{sec:twist-with-l2gamma}) that
$$\tau_{<g>}(\TR (\D^\infty \exp(-t (\D^\infty)^2)))=\Tr_{<g>}(\tilde D\exp(-t\tilde D^2))\,.$$
It is now clear that if $\tilde D$ has a gap in its spectrum
and if $\tau_{<g>}$ extends, then the integral in
(\ref{eq:delocalized_eta-truncated}) converges and the delocalized
eta invariant  $\eta_{<g>}(\tilde D)$ is well defined.
Notice that under our two assumptions
\begin{equation}\label{deloc=tau-zero}
\eta_{<g>}(\tilde D)= \tau_{<g>}( \eta_{[0]} (\D^\infty))
\end{equation}

For example, groups of polynomial growth have the property
that $\tau_{<g>}$ extends for each conjugacy class $<g>$.
More generally:

\begin{proposition}\label{prop:pol-growth}
If $<g>$ is of polynomial growth with respect to a word metric,
then $\tau_{<g>}$ has the extension property.
\end{proposition}

\begin{proof}
This follows from \cite[Lemma 6.4]{CM} and the H{\"o}lder
inequality. Compare also \cite{Ji}
\end{proof}

\begin{remark}
  It is in general very hard to check when a conjugacy class has
  polynomial growth. In particular, it seems likely that a group which
  has a exponential growth will have a conjugacy class whose growth is
  not polynomial. On the other hand, the conjugacy classes of
  commutators in a metabelian group
are in most cases of polynomial growth
\end{remark}

\begin{theorem}
  Assume that $\Gamma$ is torsion-free, $\mu_{red}\colon K_*(B\Gamma)\to
  K_*(C^*_{red} \Gamma)$ is an isomorphism and that $<g>$ is a
  non-trivial conjugacy class in $\Gamma$ of polynomial growth.
If $D=\Di$ is the spin Dirac operator on a spin manifold with
positive scalar
  curvature, then $\eta_{<g>}(\tilde \Di)$ is defined and
  $\eta_{<g>}(\tilde \Di) =0$.
\end{theorem}

\begin{proof}
We shall freely use a $
(\B^\infty_\Gamma)_{{\rm ab}}$-valued
 APS-index theory. This is the zero-degree part of the higher theory
developed in
\cite{LPGAFA} \footnote{The $\B^\infty_\Gamma$ higher APS-index
theory developed in \cite{LPGAFA}
assumes the group to be of polynomial growth. Arbitrary groups are
treated in \cite{LPBSMF} in the invertible case, based on results of
Lott \cite{Lott-torsion}. In order to extend the general  $\B^\infty_\Gamma$-theory
of \cite{LPGAFA} from groups of polynomial growth to arbitrary groups
one simple needs to prove that if $\Ind(\D^\infty)=0$
in $K_1 (\B^\infty_\Gamma)=K_1 (C^*_r \Gamma)$ then there exist trivializing
perturbations in the $\B^\infty_\Gamma$-Mishchenko-Fomenko calculus.
However, this is the consequence of a simple density argument. We
thank
Victor Nistor for pointing this out.}.
Since we are in the presence of positive scalar curvature, the fact that
$\eta_{<g>}(\tilde \Di)$ is well defined follows from Lott's results
and Proposition \ref{prop:pol-growth}.
The vanishing of  $\eta_{<g>}(\tilde \Di)$  follows
 from (\ref{deloc=tau-zero})  and from the proof of Theorem
 \ref{theo:vanish_Dirac_stable}
and Theorem \ref{theo:vanishing-higher} below,
which show that $[\eta_{[0]}(\cDi^\infty)]\in (\B^\infty_\Gamma)_{{\rm ab}}/<[1]>$
vanishes. Notice that $\tau_{<g>}:(\B^\infty_\Gamma)_{{\rm ab}}\to \CC$
factors through $(\B^\infty_\Gamma)_{{\rm ab}}/<[1]>$.
\end{proof}

\begin{remark}\label{rem:infinite-signature}
For the signature operator one can prove the following
statement:\\   {\it If $D=D^{{\rm sign}}$ is the signature operator of an oriented Riemannian
manifold, and if 0 is isolated in the $L^2$-spectrum
of $\tilde{D}^{{\rm sign}}$,
then $\eta_{<g>}(\tilde{D}^{{\rm sign}})$ is defined for each
  pair $(M^\prime,u^\prime:M^\prime \to B\Gamma)$
   $\Gamma$-homotopy equivalent to $M$ and it is an oriented
   $\Gamma$-homotopy invariant.}\\
That the delocalized eta invariants are  well defined follows from
Lott's results, form Proposition \ref{prop:pol-growth} and from the homotopy invariance
of the Novikov-Shubin numbers, see Gromov-Shubin \cite{GShubin}.
The proof of the homotopy invariance given in the previous sections
can be modified so as to cover this case as well. We omit the details.
In fact, under the gap assumption one can argue in the following
alternative and more general way.

Because of the gap assumption the signature index class
associated to $(M,u:M\to B\Gamma)$ is zero in $K_1
(\B^\infty_\Gamma)=K_1 (C^*_r \Gamma)$.
Since the Baum-Connes assembly map is by assumption bijective, we conclude
from Theorem \ref{theo:vanish-stable} that the {\it stable}
rho-invariant $\rho^s_{[0]} (M,u)$ is equal to zero in
$(\B^\infty_\Gamma)_{{\rm ab}}/<[1]>$.
However, in this case there is a very natural perturbation available
for the definition of the stable rho-invariant,
namely the projection $\Pi_{\Ker (\D^\infty)}$ onto the null space of $\D^\infty$, which is
a finitely
generated projective $\B^\infty_\Gamma$-module because of the gap assumption. Thus
$\rho^s_{[0]} (M,u)=[\eta_{[0]} (\D^\infty+\Pi)]\, =0$. On the other
hand, in general, it can be proved \cite{LPGAFA} that
\begin{equation*}
\eta_{[0]} (\D^\infty+\Pi)=\eta_{[0]} (\D^\infty)
+[\Ker (\D^\infty)]_{[0]} \equiv \eta_{[0]} (\D^\infty)
+\TR(\Pi_{\Ker (\D^\infty)}) \;\; \text{in}\;\;
\B^\infty_\Gamma\,.
\end{equation*}
We conclude that  under our three assumptions (gap + Baum-Connes +
polynomial
growth of $<g>$)
the following formula holds:
\begin{equation*}
\eta_{<g>} (\tilde D)=\tau_{<g>}(\TR (\Pi_{\Ker (\D^\infty)}))
\end{equation*}
On the other hand, because of the gap assumption, one can
establish a Hodge theorem, proving that the whole null space
$\Ker (\D^\infty)$ is in this case a homotopy invariant, being
isomorphic to the cohomology of $M$ with values
in the local system $\tilde M \times_\Gamma \B^\infty_\Gamma$.
Thus, under our assumption, the homotopy invariance of
$\eta_{<g>} (\tilde D^{{\rm sign}})$ is a consequence of
a much more general result.
\end{remark}

\subsection{Higher rho-invariants.}\label{subsect:higher}

Let us consider $(M,u\colon M\to \B\Gamma)$ with $M$ spin and with a
 metric with positive scalar curvature. In this section $M$ is not necessarily of odd
dimension.
Let $\cDi$ be the Dirac operator twisted by the
Mishchenko-Fomenko bundle $u^* E\Gamma\times_\Gamma \B^\infty_\Gamma$.
In this case, Lott's  higher eta invariant
\begin{equation*}
\tilde\eta (\cDi)\;\in \;(\Omega_* (\B^\infty_\Gamma))_{{\rm
    ab}}:= \Omega_* (\B^\infty_\Gamma)/\overline{[\Omega_*
  (\B^\infty_\Gamma),\Omega_* (\B^\infty_\Gamma)]}
\end{equation*}
is well defined, see \cite{Lott-torsion}, \cite[Appendix]{LPGAFA}.
Higher rho-invariants are obtained by pairing this noncommutative
differential form with suitable closed graded traces on $\Omega_*
(\B^\infty_\Gamma)$.
Let us describe these traces. We start with a closed graded
trace $\Phi$ on $\Omega_* (\CC\Gamma)$; we assume that
$\Phi$ is identically zero on the noncommutative differential forms concentrated
at the identity conjugacy class, i.e. on elements of the form
$$ \sum_{\gamma_0, \dots,\gamma_k; \gamma_0\cdots
  \gamma_k=1}
\omega_{\gamma_0,\gamma_1,\dots,\gamma_k}\gamma_0 d\gamma_1 \cdots
d\gamma_k\,.$$
There are examples of such traces, see \cite[p.~209]{LottI}.
 We
briefly refer to $\Phi$ as a {\it delocalized closed graded trace}.
We {\it assume} that $\Phi$ {\it extends} to a closed graded trace
$\Phi_\infty$
on  $\Omega_* (\B^\infty_\Gamma)$. The {\it higher rho-invariant}
associated to $\cDi$ and $\Phi$ is by definition the complex
number
$$\tilde\rho _{\Phi} (\cDi):=<\tilde\eta (\cDi) , \Phi_\infty>\,.$$
It is clear that the delocalized eta invariant of Lott is a special case of this
construction. We shall also use the notation  $\tilde\rho_{\Phi} (M,u):=
\tilde\rho_{\Phi} (\cDi)$.

It is proved in \cite[Proposition 4.2]{LPPSC} that the higher rho invariants
defines maps $\tilde \rho_\Phi\colon  {\rm Pos}^{{\rm spin}}_n (B\Gamma)\to
\CC$. If $\Gamma$ has torsion, the latter
 information is used in \cite{LPPSC} in order to distinguish
metrics of positive scalar curvature (under suitable assumptions on
$\Gamma$). In the torsion-free case, on the other hand, we can prove
the following

\begin{theorem}\label{theo:vanishing-higher}
Assume that $\Gamma$ is {\it torsion-free} and that the assembly map
$\mu_{red}\colon K_*(B\Gamma)\to
  K_*(C^*_{red} \Gamma)$ is an isomorphism. Let $\Phi\colon \Omega_*
  (\CC\Gamma)\to\CC$ be
a delocalized closed graded trace. Assume that $\Phi$ {\it extends}
to a closed graded trace on $\Omega_* (\B^\infty_\Gamma)$.
If  $M$ is  a spin manifold with
positive scalar
  curvature and  $u\colon M\to B\Gamma$ is a classifying map,
then for the associated Dirac operator $\cDi$,
the higher rho invariant vanish:
\begin{equation*}
\tilde\rho_{\Phi} (M,u):=
\tilde\rho_{\Phi} (\cDi)=0.
\end{equation*}
\end{theorem}

\begin{proof}
Since $M$ has positive scalar curvature the index class of $\cDi$
is equal to zero. Thus, using Proposition \ref{prop:spin_bordism}
and the injectivity of $\mu_{red}$, we conclude that for some
$d\in \NN\setminus \{0\}$, $d(M,u\colon M\to B\Gamma)$ is bordant
to $$\bigcup_{j=1}^k (A_j\times B_j,r_j\times 1:A_j\times B_j\to
B\Gamma)$$ with $\dim B_j=4b_j$, $\pi_1 (B_j)=1$ and $<\hat{A}
(B_j),[B_j]>=0$; we denote the latter manifold with classifying
map by $(N,v\colon N\to B\Gamma)$. We can and we shall endow $N$
with a metric of positive scalar curvature. Let $(W,F:W\to
B\Gamma)$ be the manifold with boundary realizing the bordism.
Since there is a metric of positive scalar curvature on $\partial
W$ the $b$-index class $\Ind_b (\cDi_W)$ is well defined in $K_*
(C^*_r \Gamma)$. By the surjectivity of $\mu_{red}$ we know that
$\Ind_b (\cDi_W)=\Ind (\cDi_X)$ with $X$ a {\it closed} spin
manifold with classifying map $r\colon X\to B\Gamma$. In
particular, using Lott's treatment of the Connes-Moscovici higher
index theorem \cite{CM}, \cite{LottI}, we see that the
Karoubi-Chern character of $\Ind_b (\cDi_W)$ is concentrated in
the trivial conjugacy class. This means that $<\Ch (\Ind_b
(\cDi_W)),\Phi_\infty>=0$. On the other hand, by the higher
APS-index theorem in \cite{LPMEMOIRS} we know that $$<\Ch(\Ind_b
(\cDi_W)),\Phi_\infty>=d\tilde\rho_\Phi (M,u)-\sum_{j=1}^k
\tilde\rho_\Phi (A_j\times B_j, r_j\times 1),$$ since the local
part in the index formula is concentrated in the trivial conjugacy
class and it is thus sent to zero by $\Phi_\infty$. Hence using
the bijectivity of $\mu_{red}$ we have proved that
$$\tilde\rho_\Phi (M,u)=\frac{1}{d}\sum_{j=1}^k \tilde\rho_\Phi
(A_j\times B_j, r_j\times 1).$$ Using now the product formula for
{\it higher} eta-invariants proved in \cite[Proposition
2.1]{LPPSC}  we see that $\tilde\rho_\Phi (A_j\times B_j,
r_j\times 1)=0$ $\forall j$. Thus $\tilde\rho_\Phi (M,u)=0$ and
the Theorem is proved.
\end{proof}

\section{The center valued $L^2$-signature formula for manifolds with
  boundary}
\label{sec:center-valued-l2}

In this section, we relate the eta- and rho invariants of the
signature operator, which show up in the APS-index theorem for the
signature operator, to the signature of the manifold with boundary
(the signature of some possibly degenerate intersection form).

Note that, even in the compact case, the ordinary signature formula for manifolds with
boundary does not immediately follow from the Atiyah-Patodi-Singer
index theorem for the signature operator. It is a non-trivial result
of \cite{APS1}
to connect the APS-index of the signature operator to the signature of
the manifold with boundary. This is much more complicated in the
$L^2$-case, since in \cite{APS1} eigenvalue decompositions of the
space of $L^2$-sections on the boundary are used, which are not
available in our setting. This is overcome in \cite{LS03} together
with \cite{Vaillant} for the numerical $L^2$-signature. We now explain
how this is done  and  how it generalizes to the situation we
are considering in Section \ref{sec:examples-non-trivial}.

Therefore, let $W$ be a compact oriented Riemannian manifold of
dimension $4k$ with
boundary $\boundary W=M$. Let $A$ be a von Neumann algebra and
$\tau\colon A\to Z$ a finite positive normal trace with values in an abelian
von Neumann algebra $Z$. We are thinking here in particular of the von Neumann
algebra of a discrete group $\Gamma$ with its canonical trace or with
its center valued trace.

Let $\mathcal{L}$ be a \emph{flat} bundle of finitely generated projective
Hilbert $A$-modules on $W$ (giving rise to a local coefficient system
of such modules). Recall that $\mathcal{L}$ is given by a
representation of $\Gamma$ in a finitely generated projective Hilbert
$A$ module (which we call $\mathcal{L}_x$ here). We assume that everything involved
(Riemannian metric, bundle, connection) is of product type near the
boundary.

We can now define three kinds of intersection forms on $W$,  using the
twisting bundle $\mathcal{L}$, and with a signature in $K_0(A)$.

The most computable one is obtained combinatorially: we consider a
triangulation (or more generally a CW-decomposition) of $W$. This
defines a cellular chain complex $C_*(W;\mathcal{L})$ of finitely
generated free Hilbert $A$-modules, with coefficients in the local
coefficient system $\mathcal{L}$. There is a
Poincar{\'e} duality chain homotopy equivalence to the relative cochain
complex $C^{4k-*}(W,\boundary W;\mathcal{L})$. From there we can restrict
to $W$ to get a map to $C^{4k-*}(W;\mathcal{L})$. Since the second
map is not a chain homotopy equivalence in general, neither is the
composition. But it is self-dual (note that the cochain complex is
dual to the chain complex).

Now one can pass to the projective part of the homology and cohomology
of these Hilbert $A$-module (co)chain complexes. This passage to the
projective part involves some additional consideration in homological
algebra  ---special to finite von Neumann algebras---, developed
in different languages and independently by Farber
\cite{MR2000b:58041} and L{\"u}ck \cite{MR99k:58176}. We will
later only look at the special example where the homology and
cohomology in middle degree is itself a finitely generated free
Hilbert $A$-module (in fact equal to the chain- and cochain
module). Then the projective part in middle degree is equal to the
whole (co)homology.
The Poincar{\'e} duality chain homotopy equivalence  composed with
restriction to the boundary will then define a self dual map
\begin{equation*}
  H_{2k}(W;\mathcal{L}) \to H^{2k}(W;\mathcal{L})\to
  \Hom_{\mathcal{A}}(H_{2k}(W;\mathcal{L}),\mathcal{A}).
\end{equation*}

In the special situation we are going to consider this will be given
as follows:
There is a free finitely generated $\integers\Gamma$-module
  $V\iso (\integers\Gamma)^l$
 with (possibly singular) self dual map $\sigma\colon V\to
  V^*:=\Hom_{\integers\Gamma}(V,\integers\Gamma)$ of the form
  $\sigma=\psi+(-1)^k \psi^*$ (i.e.~$\sigma$ has a quadratic
  refinement).
  Identifying $V^*$ with $V$ using the given basis, $\sigma$ is
  represented by a matrix $B=A+A^*$ with $A\in M_n(\integers\Gamma)$.
The self dual map $H_{2k}(W;\mathcal{L})\to H_{2k}(W;\mathcal{L})^*$
is then obtained by tensoring $\sigma\colon V\to V^*$ with the
$\Gamma$-representation $\mathcal{L}_x$, i.e.~$H_{2k}(W;
\mathcal{L})\iso V\tensor_{\integers\Gamma}\mathcal{L}_x\iso
\mathcal{L}_x^l$, and the map is given as $B':=B\tensor \id_{\mathcal{L}_x}$.

The combinatorial signature $\sgn(W,\mathcal{L})\in
K_0(\mathcal{A})$ is then defined by
\begin{equation*}
\sgn(W,\mathcal{L}):=
\sgn_{\mathcal{A}}(B'):=\chi_{(0,\infty)}(B')-\chi_{(-\infty,0)}(B') \in K_0(\mathcal{A}).
\end{equation*}
Note
that, for manifolds with boundary, this can only be defined for a
von Neumann algebra, where measurable functional calculus is available,
since $0$ might well be contained in the spectrum of $B'$ (in contrast
to the closed case, where the topology implies that an intersection
form is necessarily invertible).

Using the trace $\tau\colon \mathcal{A}\to Z$, we can then define the $Z$-valued
combinatorial signature
\begin{equation*}
  \sgn_\tau(W,\mathcal{L}):= \sgn_\tau(B'):=\tau(\sgn_A(B')).
\end{equation*}

The second version of $A$-signature is obtained by replacing the
combinatorial chain complex by the chain complex of differential forms
with values in the flat bundle $\mathcal{L}$. Here, the cup product of
forms together with the $A$-valued inner product in the fibers of
$\mathcal{L}$ induces an $A$-valued intersection pairing on the
de Rham cohomology in middle degree with coefficients in
$\mathcal{L}$. This intersection pairing is again defined by a self adjoint operator
$B_{dR}$, and then $\sgn_{\mathcal{A}}(B_{dR}):=
\chi_{(0,\infty)}(B_{dR})-\chi_{(-\infty,0)}(B_{dR}) \in K_0(\mathcal{A})$ defines the de Rham
signature of $W$ with coefficients in $\mathcal{L}$.

Finally, we can attach infinite cylinders to the boundary of $W$, and
then carry out the construction as above, but now with square
integrable differential
forms (with values in $\mathcal{L}$) on the enlarged non-compact
manifold. We get a third signature invariant in $K_0(A)$.

Application of $\tau$ gives three signatures in $Z$.

One can now use the arguments of \cite{Lueck-Schick(2001b)} in order
to show that
these three invariants in fact always coincide. In
\cite{Lueck-Schick(2001b)}, the corresponding result is proved for the
ordinary $L^2$-signature, i.e.~where $\mathcal{L}=\NeumannN$ is the
Mishchenko-Fomenko line bundle, and
$\tau\colon\NeumannN\Gamma\to\complexs$ is the canonical trace. The techniques carry
over. The relevant properties are that  $\tau$ is a positive and normal
trace. Note that in \cite{Lueck-Schick(2001b)}, one works with $\mathcal{A}$-Hilbert
spaces instead of the Hilbert
$\mathcal{A}$-modules used here. However, one can translate between
these two settings as
explained in \cite[Section 8.6]{Schick03}.

Let $D_\mathcal{L}$ be the signature operator on $W$ twisted by $\mathcal{L}$. It has
boundary operator $D_{M,\mathcal{L}}$. Then we can express the three
equal higher
signature of the bordism $W$ which were defined above  using the signature operator.
\begin{theorem}\label{theo:A-valued_boundary_signature_theorem}
    \begin{equation*}
  \sgn_\tau(W;\mathcal{L}) = \left(\int_W AS(D)\cdot \tau(\mathcal{L}_x)\right) -
  \frac{\eta_{\tau}(D_{M,\mathcal{L}})}{2}\;\in Z.
\end{equation*}
The Atiyah-Singer integrand $AS(D)$ for the signature operator is of
course given by the Hirzebruch $L$-form of the Riemannian manifold
$W$. $\tau(\mathcal{L}_x)$ is the trace of the projection onto the finitely
generated projective fiber $\mathcal{L}_x$, this is a locally constant
$A$-valued function on $W$.
\end{theorem}
\begin{proof}
  This is proved for the ``ordinary'' cylindrical end $L^2$-signature
  by Vaillant \cite{Vaillant}, using the $L^2$-Atiyah-Patodi-Singer
  index formula. His proof only uses the formal
  properties
  of the canonical $L^2$-trace of being positive and normal and
  therefore carries over to prove the asserted equality. As explained
  above, from the work in
  \cite{Lueck-Schick(2001b)} it follows that the formula also holds for the
  combinatorially defined signature.
\end{proof}

\section{Examples of non-trivial rho-invariants}
\label{sec:examples-non-trivial}

In this section, we show that there are many examples where the
rho-invariants considered in this paper (know to be homotopy invariants)  show
that certain manifolds are not homotopy equivalent.
The fundamental group can be as simple as $\integers$. Similar explicit
calculations (without using the notion of $\rho$-invariant) have been
carried out in \cite[Section 5]{Cochran-Orr-Teichner}. In the latter paper,
these invariants are used to detect certain knots which do not have
the slice property
(which implies that certain types of bordisms can not exist).

We use the surgery construction of \cite{LS03}
employed there to construct manifolds with boundary with very general
intersection form.

Recall from Section \ref{sec:center-valued-l2} that the
combinatorially defined $L^2$-signatures of a triangulated
manifold $M$ (with or without boundary) is obtained as follows. Assume
that $\Gamma=\pi_1(M)$. We
have an ``intersection form'' on the combinatorial chain complex. This
can be understood as a matrix $B$ with entries in
$\integers\Gamma$. Actually, because of the symmetry properties of the
intersection form, $B=A+A^*$ for $A\in M_n(\integers\Gamma)$.
We now use $\mathcal{A}=\NeumannN\Gamma$ and
$\mathcal{L}=\NeumannN=\tilde M\times_\Gamma\NeumannN\Gamma$. In this
situation, since $\integers\Gamma$ is a subset of $\NeumannN\Gamma$,
we understand the matrices $A$ and $B$ to be matrices also over
$\NeumannN\Gamma$ (i.e.~we write $B$ instead of $B'$ in the notation
of Section \ref{sec:center-valued-l2}).

Then $\sgn(M,\NeumannN)=
\sgn_{\NeumannN\Gamma}(B)=\chi_{(0,\infty)}(B)-\chi_{(-\infty,0)}(B)$.

Applying the center-valued trace $\tau\colon \NeumannN\Gamma\to Z$, we
get then $\sgn(M,\NeumannN,\tau)= \tau(\sgn(M,\NeumannN))$. Set
$\sgn_{\tau}(B):=\tau(\sgn_{\NeumannN\Gamma}(B))$. Applying the
canonical trace, from this we get $\sgn_{(2)}(B)\in \complexs$, and
applying the delocalized traces corresponding to a finite conjugacy
class $<g>$, we get $\sgn_{<g>}(B)\in\complexs$.

Moreover, if $\lambda\colon\Gamma\to U(d)$ is a finite dimensional
representation of $\Gamma$, set $\lambda(B)\in
M_n(M_d(\complexs))=M_{nd}(\complexs)$, the matrix obtained by
applying $\lambda$ entry-wise. Then $\sgn(\lambda(B))$ is the $\lambda$-twisted
signature of $M$. Recall that $\lambda(B)$ is a possibly indefinite
self-adjoint matrix. $\sgn(\lambda(B))$ is the difference of the
number of positive and negative eigenvalues, the eigenvalue $0$ is ignored.

We will use the following Proposition 1.1 of
\cite{Lueck-Schick(2001b)}.
\begin{proposition}\label{prop:bordismsP_with_prescribed_signature}
    Fix a dimension $4k\ge 6$ and a finitely presented group $\Gamma$. Let $X$ be a closed
  $(4k-1)$-dimensional manifold with fundamental group $\Gamma$ and with
  Morse decomposition without a $k$-handle.  Choose $A\in
  M_n(\integers\Gamma)$. Let $B=A+A^*$.

  Then there is a compact manifold with boundary $(W; X, Y)$ of
  dimension $4k$ with boundary
  $\boundary W=X\disjointunion Y$ and with fundamental group $\Gamma$,
  such that the Morse chain complex
  $C_*(\tilde W)$ of the universal covering $\tilde W$ is isomorphic
  to $C_*(\tilde X)\oplus V$, where $V$ is considered as trivial chain
  complex concentrated in the middle dimension $k$, and with inverse
  Poincar{\'e} duality homomorphism
  \begin{equation*}
    C_{4k-*}(\tilde W) \to C_{4k-*}(\tilde W,\boundary\tilde
  W) \xrightarrow{PD^{-1}}
    C^*(\tilde W)
  \end{equation*}
  which in the middle dimension is  given by $B:=A+A^*$, as explained
  in Section \ref{sec:center-valued-l2}. Here $PD^{-1}$ is
  a chain homotopy inverse to the cup product with $[W,\boundary W]$.

  In particular,
  $\sgn_\tau(W,\NeumannN)=
  \tau(\chi_{(0,\infty)}(B))-\tau({(-\infty,0)}(B))$, where
  $\tau\colon\NeumannN\Gamma\to Z$ is the center valued trace.
\end{proposition}

Recall that, on the other hand, for the constructed manifold $W$ the
ordinary signature is the signature of $\rho(B)$, where $\rho\colon
\Gamma\to \{1\}$ is the trivial representation.

Using the $L^2$-signature theorem
\ref{theo:A-valued_boundary_signature_theorem}, we therefore get in
principle a large number of examples of manifolds with interesting
difference of $\rho$-invariants.

\begin{corollary}
  Given any $A\in M_n(\integers\Gamma)$, set $B:=A+A^*$. Then there
  exist manifolds $X,Y$ with $\pi_1(X)=\pi_1(Y)=\Gamma$ such that
  \begin{equation*}
    \rho_{(2)}(X)-\rho_{(2)}(Y) = \sgn_{(2)}(B)-\sgn(\rho(B)).
  \end{equation*}
  If $\lambda_1,\lambda_2\colon\Gamma\to U(d)$ are two
  representations, then
  \begin{equation*}
    \rho_{\lambda_1-\lambda_2}(X) -\rho_{\lambda_1-\lambda_2}(Y) = \sgn(\lambda_1(B))-\sgn(\lambda_2(B)).
  \end{equation*}
  If $g\in\Gamma$ has finite conjugacy class $<g>$, then
  \begin{equation*}
    \rho_{<g>}(X)-\rho_{<g>}(Y) = \sgn_{<g>}(B).
  \end{equation*}
\end{corollary}

\begin{remark}
    Note that $X\disjointunion -Y$ with the reference map to $B\integers$
  is a boundary and therefore has
  signature class $0\in K_0(C^*\integers)$.
  Nevertheless, our  construction shows that the
  rho-invariants of $X\disjointunion -Y$ are non-zero. In particular,
  for such a manifold the stable and the unstable rho-invariant can
  differ (we know that the stable one vanishes in such an example).
\end{remark}

This theorem implies that we have a great freedom at constructing
manifolds $X$ and $Y$ such that the various $\rho$-invariants
differ. There is of course a problem in explicitly calculating these
invariants for a given matrix $B$. Let us therefore recall how this
can be done in the easiest case, i.e.~if $\Gamma=\integers$.

In this case $\NeumannN\integers\iso L^\infty(S^1)$ using Fourier
transform. The subring $\integers[\integers]\subset
\NeumannN\integers$ corresponds to Laurant polynomials on $S^1$.

Therefore, $B=A+A^*\in M_n(\integers\Gamma)$ can be understood as an
$n\times n$-matrix with entries in Laurant polynomials on $S^1$, or
alternatively as a function (a Laurant polynomial) $B(z)$ on $S^1$ with
values in $M_n(\complexs)$. To compute $\sgn_{\NeumannN\integers}(B)$,
one then computes pointwise $\sgn(B(z))$, getting an integer valued
function on $S^1$. Considered as an element of
$L^\infty(S^1)=\NeumannN\integers=Z$, this is exactly $\sgn_\tau(B)$
(since $\NeumannN\integers$ is abelian, the center valued trace is the
identity).

If $\lambda\colon \integers\to U(1)$ is a $1$-dimensional
representation, sending the generator $z$ to $\lambda\in U(1)=S^1$,
then $\sgn_\lambda(B)=\sgn(B(\lambda))$, i.e.~we have to evaluate the
function $\sgn(B(z))$ at the point $\lambda$. In particular, for the
trivial representation, $\sgn_1(B)=\sgn(B(1))$.

If, on the other hand, $g=z^n\in\integers$, then $\sgn_{<g>}(B)$ is the
Fourier coefficient of $g$ for the function $\sgn(B(z))$, i.e.
\begin{equation*}
  \sgn_{<g>}(B) = \int_{S^1} \sgn(B(z))z^{-n}\,dz,\text{ in particular
    }\sgn_{(2)}(B)=\int_{S^1}\sgn(B(z))\,dz.
\end{equation*}

Lets look at a particular example: if $A=z+z^{-1}+1$ (a $1\times 1$-matrix), then
$B(z)=2(z+z^{-1}+1)$, and $\sgn(B(z))=0$ for $z=z_{1,2}=\exp(\pm 2\pi i/3)$,
$\sgn(B(z))=1$ if $z$ is contained in the connected component of $1$
of $S^1\setminus\{z_1,z_2\}$, and $\sgn (B(z))=-1$ for the remaining
points of $z$.

Observe that this is the general pattern: $\sgn(B(z))$ jumps (at most)
at those points on $S^1$ where an eigenvalues of $B(z)$ crosses $0$,
i.e.~where the rank of $B(z)$ is lower than the maximal rank.

It follows e.g.~that $\sgn_{(2)}(B)= 1/3$, but $\sgn_{1}(B)=1$, and
$\sgn_{<z>}(B) = 2(\exp(2\pi i/3)-\exp(-2\pi i/3))=2i\sqrt{3}$.

In particular, since $\sgn_{(2)}(B)-\sgn_{1}(B)\ne 0$, (or since
$\sgn_{<z>}(B)\ne 0$), if we use the matrix $A$ in the construction of
Proposition \ref{prop:bordismsP_with_prescribed_signature}, the
resulting closed manifolds $X$ and $Y$ satisfy
\begin{equation*}
  \rho_{(2)}(Y)-\rho_{(2)}(X) = \sgn_{(2)}(B)-\sgn_{1}(B) = -2/3,
\end{equation*}
and are therefore not homotopy equivalent.

Note that we can apply Theorems \ref{theo:APS_rho}, \ref{theo:L2_rho}
or \ref{theo:finite_conjugacy} since $\integers$ is torsion-free and
satisfies the Baum-Connes conjecture for the maximal $C^*$-algebra.

Unfortunately, this conclusion has one flaw: by \cite[Lemma
2.2]{MR2000a:57089} and the proof of \cite[Theorem
5.8]{MR2000a:57089},
\begin{equation}\label{eq:surgery_homology}
H_{k}(Y;\integers\Gamma)\iso \ker(B)\oplus
H_k(X;\integers\Gamma)\text{ and }
H_{k-1}(Y;\integers\Gamma)\iso\coker(B)\oplus
H_{k-1}(X;\integers\Gamma),
\end{equation}
whereas the homology of $X$ and $Y$ in all other degrees
coincides. This implies that we might as well distinguish $X$ and $Y$
using their homology, which is of course much easier to compute that
the $\rho$-invariants.

One can easily construct more interesting examples as
follows: consider a diagonal matrix $A$ with entries $A_1(z),\dots,
A_m(z)$, and a second diagonal matrix $A'$ with entries $\epsilon_1
A_1(a),\dots,\epsilon_n A_n(z)$, with $\epsilon_i\in \{-1,1\}$. Starting with a manifold $X$ as in
Proposition \ref{prop:bordismsP_with_prescribed_signature}, we then
get two manifolds $Y$ and $Y'$. By \cite[Lemma 2.2]{MR2000a:57089} and
the proof of \cite[Theorem 5.8]{MR2000a:57089}, the homology of $Y$
and $Y'$ is isomorphic in all degree and with arbitrary coefficients.

However, the signatures and therefore the difference of
$\rho$-invariants changes sign if the matrix $A$ changes sign. This
means that we can easily arrange that certain rho-invariants of $Y$
and $Y'$ do not coincide, we have even enough freedom to make sure
that there are examples where they neither coincide nor are negative
of each other.

The conclusion is that
although the homology of $Y$
and $Y'$ is isomorphic in all degree and with arbitrary coefficients,
there is no homotopy equivalence between $Y$ and $Y'$ (not even one which reverses the
orientations). This result has been obtained using rho-invariants,
it can be recast in terms of Blanchfield forms (also called linking
forms), i.e.~in terms of classical methods of advanced algebraic topology.

It is evident that with other matrices $A(z)$, we can get all kinds of further
interesting examples.

This was all done for a group as simple as $\integers$. One is
often interested to get examples for more complicated groups,
compare also the problem of constructing knots of
e.g.~\cite{Cochran-Orr-Teichner}. We can ``induce up'' these
examples by embedding $\integers$ in an arbitrary torsion-free
group $\Gamma$ and then use the fact that the signature
calculations happen entirely in the subgroup $\integers$ and
therefore are unchanged, as long as the $L^2$-signature and the
ordinary signature are involved. The relevant result is  stated
and proved e.g.~in \cite[Proposition 5.13]{Cochran-Orr-Teichner}.
The following is happening: when we have an inclusion
$i\colon\integers\into\Gamma$, then this induces inclusions
$l^2(\integers)\into l^2(\Gamma)$ and $i\colon
\NeumannN\integers\into\NeumannN\Gamma$. The latter one being an
inclusion of von Neumann algebras, it is compatible with
functional calculus. If we therefore start with $B=A^*+A\in
M_n(\integers\Gamma)$, then
$\chi_{(0,\infty)}(i(B))=i(\chi_{(0,\infty)}(B))$. It is
not to us clear how the classical methods of algebraic
topology mentioned above could to be used in this general setting.

Finally, if we want to compute a delocalized trace of an element
$i(x)$, we observe that by definition for $g\in\Gamma$ with
$\abs{<g>}<\infty$.
\begin{equation*}
  \tau_{<g>}(i(x))=\sum_{h\in <g>} \innerprod{e\cdot i(x),h}_{l^2(\Gamma)}.
\end{equation*}
But since $e\in l^2(\integers)$, also $e\cdot i(x)\in
i(l^2(\integers))$, and therefore the inner product is zero if
$h\notin i(l^2(\integers))$, and
\begin{equation*}
  \tau_{<g>}(i(x)) = \sum_{h\in<g>\cap i(\integers)} \innerprod{e\cdot
    i(x),h}_{l^2(\Gamma)} = \sum_{h\in <g>\cap i(\integers)}
  \tau_{i^{-1}(h)} (x).
\end{equation*}
In particular, if $g=e$ then $\tr_{(2)}(i(x))=\tr_{(2)}(x)$, and if
$g\ne e$, one has to analyze which conjugates of $g$ lie in
$i(\integers)$. If $<g>$ is finite, then $g^{-1}$ can be the only
power of $g$ contained in $<g>$. It follows that if $g$ generates
$i(\integers)$, we understand $<g>\cap i(\integers)$ completely.

For example, if one embeds $\integers$ into the center of $\Gamma$,
all delocalized invariants with respect to elements of this embedded
$\integers$ will be the same, independent of the question whether the
corresponding matrix is considered as a matrix over $\integers$ or
over $\Gamma$. With a little extra care one can use the same
``induction process'' to obtain examples of non-trivial delocalized
rho-invariants for finite conjugacy classes with more than one element.

\begin{appendix}

\section{Index theory: proofs.}\label{sec:proof}

This first appendix is devoted to the proof of the two index theorems
stated in Section \ref{sec:degree-zero-as}.

\subsection{Proof of the $A/\overline{[A,A]}$ valued Atiyah-Singer index formula}
\label{subsec:proof-closed}

First of all we need to deal with the existence of the heat kernel
for the Dirac Laplacian $D_{\L}^2$ and its perturbation
$D_{\L,\C}^2\equiv (D_{\L}+\C)^2$.

\begin{lemma}\label{lem:heat semigroup closed}
 For each $t>0$ there exists a well-defined  operator
 $e^{-tD_{\L,\C}^2}\;\in\; \Psi^{-\infty}_A\,.$
 If $v\in C^\infty (M,E\otimes \L)$, then
 $\,u(t,\cdot):=e^{-tD_{\L,\C}^2}v\,$ is a solution of the heat equation
 $(\partial_t + D_{\L,\C}^2)u=0$
 with initial condition $v$ at $t=0$.
 The operators $e^{-tD_{\L,\C}^2}$ form a semigroup.
\end{lemma}

\begin{proof}
We begin by the Dirac Laplacian $D_{\L}^2$. The best way to prove
the above lemma is by use of the {\it heat space} of Melrose, see
\cite{Melrose} Chapter VII. First we make a general comment.
Microlocal analysis can be viewed as a geometrization of operator
theory, the operators of interest being certain distributions on
$M\times M$ with precise singularities on the diagonal (conormal
distributions). Smoothing operators, for example, are  given as
smooth functions on $M\times M$, or, more generally, as smooth
sections of  suitable homomorphism bundles on $M\times M$. It is
not difficult to understand that the Mishchenko-Fomenko
pseudodifferential calculus can be easily developed by simply
considering $A$-valued conormal distributions. For example,
smoothing operators in the Mishchenko-Fomenko calculus, acting for
simplicity on the trivial bundle $\CC^\ell\otimes A$ are nothing
but $C^\infty$ functions on $M\times M$ with values in
$M_{\ell\times \ell} (A)$. From the microlocal point of view the
presence of the $C^*$-algebra $A$ does not affect in a significant
way all the usual arguments culminating in  the existence of the
pseudodifferential calculus. (Needless to say, operators in the
Mishchenko-Fomenko calculus act on $A$-Hilbert modules and not on
Hilbert spaces, it is at this point that much more care is
needed.) This general philosophy will be applied here to the
construction of the heat semigroup $e^{-tD_{\L}^2}$. The advantage
of the treatment given by Melrose through the heat space (a
certain blow-up of $M\times M\times [0,\infty)_t$) is that it is
as geometric as it can be, thus generalizing without any
difficulty to operators acting on sections of bundles of
$A$-modules (such as $E\otimes \L$). This general principle has
been already applied in \cite{LPMEMOIRS} in the case $A=C^*_r
\Gamma$ but it is clear that it extends readily  to an arbitrary
unital $C^*$-algebra $A$.\\ {\it Summarizing}: by following
closely the treatment given by Melrose for ordinary Dirac operator
we can prove the lemma for the heat kernel associated to $D_{\L}
^2$. A standard argument involving a Volterra series can be
applied in order to obtain the lemma for the perturbed operator
$D_{\L,\C}^2$. See \cite{BGV}.
\end{proof}

The problem we encounter with the heat kernel in the
$C^*$-algebraic context, is that it doesn't behave well for
$t\to\infty$, since our operators are in general not invertible in
the Mishchenko-Fomenko calculus. For this reason, we introduce
suitable perturbation which make them invertible.

\begin{definition}
  Instead of $M$ consider $M_+$, the disjoint union
  of the manifold with an additional point $*$. The Dirac type operator on
  the additional point is by definition the $0$ operator on
  $\complexs$. Recall the finitely generated projective modules
$\mathcal{I}_+$, $\mathcal{I}_-$ appearing in the
Mishchenko-Fomenko decomposition theorem. We define $\L_+$, the
twisting bundle on $M_+$ to be
  $\L\disjointunion (\mathcal{I}_+\oplus \mathcal{I}_-)$, where we view
  $\mathcal{I}_{\pm}$ as an abstract finitely generated projective Hilbert
  $A$-module (therefore a bundle over the point).

  Following \cite{LottI}, Section VI,
  we now define for $\alpha\in\reals$ the perturbed operators (of
  $D_{\L}$, but we suppress the $\L$ in the notation)
  \begin{equation*}
    D^+_\alpha :=
    \begin{pmatrix}
      D^+_{\L}|_{\mathcal{I}_+ ^\perp} & 0 & 0\\
      0 & D^+_{\L}|_{\mathcal{I}_+} &  \alpha\\
      0 & \alpha & 0
    \end{pmatrix},
  \end{equation*}
  and  $D^-_\alpha$ as the (formal) adjoint of
  $D^+_{\alpha}$.
  The description is with respect to the Mishchenko-Fomenko
  decomposition of Theorem \ref{lem:Mishchenko_decomposition_closed}
  \begin{equation*}
C^\infty(M_+,(E\tensor \L)_+^+)= \mathcal{I}_+^{\perp}\oplus
\mathcal{I}_+
  \oplus \mathcal{I}_-\quad  C^\infty(M_+,(E\tensor \L)_+^-)=
  D_{\L}(\mathcal{I}_+^{\perp})\oplus I_-
  \oplus \mathcal{I}_+
\end{equation*}
Note that $\mathcal{I}_{\pm}$ has two roles here: first as subset
of the section of $E\tensor \L$ on $M$, secondly as possible
values of the sections at the additional point $*$, where the
fiber is $\mathcal{I}_-\oplus\mathcal{I}_+$.
\end{definition}

\begin{lemma}
  The operator $e^{-tD_\alpha^2}$ is defined for each $t>0$ and each
  $\alpha\in\reals$, and is a smoothing operator in the
  Mishchenko-Fomenko calculus on $M_+$. In particular
  $\STR(e^{-tD_\alpha^2})\in A_{{\rm ab}}$ is defined.
\end{lemma}

\begin{proof}
The heat operator  $e^{-tD_\alpha^2}$ is defined by Duhamel's
expansion, using the fact that $D_\alpha - D_0$ is finite rank as
$A$-linear operator, see \cite[Section VI]{LottI}. Since the
heat-operator of $D_0^2$ is smoothing on $M_+$, the lemma follows.
In fact, by using the information that the orthogonal projection
onto $\mathcal{I}_+$ and the projection onto $\mathcal{I}_-$ along
$D_{\L} (\mathcal{I}_+^\perp)$ are smoothing operators,  it is
possible to check that $D_\alpha-D_0$ is smoothing on $M_+$, i.e.
an element in  $\Psi_A^{-\infty}(M_+, (E\otimes\L)_+,
(E\otimes\L)_+)$.

\end{proof}

The following lemma is clear:

\begin{lemma}\label{lem:alpha=0}
  Let $\alpha=0$, 
  then $$\STR(e^{-tD_0^2}) = \STR_M(e^{-tD_{\L}^2}) +
  [I_-]_{[0]}-[I_+]_{[0]} \in
  A_{{\rm ab}}\,.$$
\end{lemma}


\begin{lemma}\label{lem:large time alpha}
  For $\alpha$ sufficiently large, $D_\alpha$ is invertible in the MF-calculus, and
  therefore
  \begin{equation}\label{large time alpha}
  \STR(e^{-tD_\alpha^2})\xrightarrow{t\to\infty} 0 \in A_{{\rm ab}}.
\end{equation}
\end{lemma}

\begin{proof}
The invertibility is explained in \cite{LottI} Section VI, after
formula (107). One can proceed as in \cite{Melrose} and show that
the heat kernel defined through the heat space is indeed
expressible in terms of the usual integral involving the
resolvent. Using the invertibility of $D^2_\alpha$ one then gets
(\ref{large time alpha}).
\end{proof}

\begin{lemma}\label{lem:independence_of_alpha_no_w}
 For each $\alpha\in\reals$ and $t>0$
\begin{equation*}
  \STR(e^{-tD_\alpha^2}) -\STR(e^{-tD_0^2}) = 0 \in A_{{\rm ab}}.
\end{equation*}
\end{lemma}
\begin{proof}
  Use Duhamel's formula to compute the derivative with respect to
  $\alpha$ of $\STR(e^{-tD_\alpha^2})$. The usual calculations show
  that this is the supertrace of a supercommutator and therefore
  vanishes in $ A_{{\rm ab}}$. Details are as in \cite{BGV}, Chapter 3.
\end{proof}

As a Corollary we obtain Proposition \ref{prop:large time closed no alpha}:

\begin{corollary}\label{cor:large time closed no alpha}
  \begin{equation*}
  \lim_{t\to\infty}\STR_M (e^{-tD^2_{\L}}) = [I_+]_{[0]} - [I_-]_{[0]}\equiv
  \Ind_{[0]} (D_L) \in A_{{\rm ab}},
\end{equation*}
where part of the assertion is that the limit exists.
\end{corollary}

\begin{proof}
By lemmas \ref{lem:alpha=0} and
\ref{lem:independence_of_alpha_no_w} we have in
$A_{{\rm ab}}$, for each $\alpha\in\reals$ and $t>0$:
$$\STR_M (e^{-tD^2_L})=
 \STR (e^{-tD_0^2}) -
  [I_-]_{[0]}+[I_+]_{[0]}= \STR(e^{-tD_\alpha^2}) -
  [I_-]_{[0]}+[I_+]_{[0]}\,.$$
  Taking $\alpha$ large enough we get the Corollary by applying
  lemma \ref{lem:large time alpha}.
  \end{proof}

\subsubsection{The integral operator index}
\label{sec:integr-oper-index}

Next we tackle the problem of connecting the index $\Ind_{[0]}
(D_{\L})\in  A_{{\rm ab}}$ defined using the index class
and the algebraic trace  $\tr^{{\rm alg}}:K_0 (A) \to
A_{{\rm ab}}$ to the integral-kernel-trace, $\TR$, of the
projection operators given by the Mishchenko-Fomenko
decomposition.

\begin{definition}
  Define a second smoothing perturbation of $D$ by
  \begin{equation*}
B^+_{\alpha}:=
  D^+_{\L} -\alpha 
P_- D^+_{\L} P_+\,, 
\end{equation*}
for $\alpha\in\reals$, with
  $B^-_{\alpha}$ the (formal) adjoint of $B^+_{\alpha}$
and with $P_+:= \Pi_{\mathcal{I}_+}$, $P_-:= \Pi_{\mathcal{I}_-}$.
\end{definition}
\begin{remark}
  Observe that $$B_{\alpha}=\begin{pmatrix} 0& B^-_\alpha\\
                    B^+_\alpha & 0 \end{pmatrix}
  $$ is a smoothing perturbation of $D_{\L}$.
Note that $P_+$ is self
  adjoint, but $P_-$ is not necessarily. Thus we will also use its adjoint
  $P_-^*$. Note, finally,  that $\im(1-P_-) = D^+_{\L}(\mathcal{I}_+^\perp)$.
\end{remark}

\begin{lemma}\label{lem:pr_smoothing}
  We have decompositions
  \begin{equation*}
    C^\infty(M,(E\tensor \L)^-)=
       \im(P_-) \oplus \im(1-P_-) = \im(P_-^*)\oplus \im(1-P_-).
    \end{equation*}
    The second decomposition is an orthogonal decomposition, since
    \begin{equation*}
\im(1-P_-)^\perp = \ker(1-P_-^*)=\im(P_-^*).
\end{equation*}
Let $\pr$ be the
    orthogonal projection onto $\im(P_-^*)$.  Then $\pr$ is also a
    smoothing operator.
\end{lemma}
\begin{proof}
  $P_-^*P_- + (1-P_-)=1+(P_-^*P_--P_-)$ is an isomorphism of
  $C^\infty(M,E\tensor
  \L)^-)$ which is diagonal with respect to the two decompositions (it
  is an isomorphism, since $P_-^*$ has kernel $\im(P_-)^\perp$ and is
  surjective, so that $P_-^*\colon \im(P_-)\to \im(P_-^*)$ is an
  isomorphism by the open mapping theorem). Then
  \begin{equation*}
    \pr = (1+(P_-^*P_--P_-))^{-1} P_- (1+(P_-^*P_--P_-))
  \end{equation*}
  is smoothing by Lemma \ref{lem:inverse_of_one_plus_smoothing} below.
\end{proof}

\begin{lemma}\label{lem:inverse_of_one_plus_smoothing}
  If $P$ is smoothing and $1+P$ is invertible in the sense of
  Hilbert-$A$-module morphisms (on the completed space of sections),
  then
  \begin{equation*}
    (1+P)^{-1} = 1+ Q\quad\text{with } Q\text{ smoothing}.
  \end{equation*}
\end{lemma}
\begin{proof}
  Write $(1+P)^{-1}=1+Q$; then $Q$ satisfies $Q=-P-P^2-PQP$ and it
  is therefore smoothing.   Here we are using the fact that smoothing operators
  in $L^2$ are a semi-ideal (a subring $R$ of a ring $R$ is a semi-ideal if $i,j\in I$, $r\in R$
  $\Rightarrow $ $irj\in I$, for all $i,j\in I$, for all $r\in R$).
\end{proof}

\begin{lemma}
  In this lemma we shall suppress the $\L$ subscript
  in the notation of $D_{\L}$.
  For $\alpha=0$, $B_\alpha=D_{\L}$; for $\alpha=1$, we have
  \begin{equation*}
    B^+_{1} =
    \begin{pmatrix}
      D|_{I_+^{\perp}} & 0 \\  0 & 0
    \end{pmatrix}
  \end{equation*}
  with respect to the decomposition
    \begin{equation*}
C^\infty(M,(E\tensor \L)^+)= \mathcal{I}_+^{\perp}\oplus
\mathcal{I}_+
  \qquad  C^\infty(M,(E\tensor \L)^-)= \im(1-P_-)\oplus \im(P_-^*).
\end{equation*}
Observe that this operator decomposes as an invertible operator plus
(direct sum) the zero operator between two finitely generated
projective modules.

Since the decompositions of domain and range are both orthogonal,
the adjoint $B^*_{1}$ decomposes (with respect to the same
decomposition) as
\begin{equation*}
  B^-_{1} =
  \begin{pmatrix}
    D^*|_{I_+^\perp} & 0 \\ 0 & 0
  \end{pmatrix},
\end{equation*}
and again the left upper corner is an isomorphism. From this
\begin{equation*}
  B_{1}^2 =
  \begin{pmatrix}
    D_{I_+^\perp}^*D_{I_+^\perp} & 0 & 0 & 0\\
    0& 0 & 0& 0\\
    0&0 & D_{I_+^\perp}D_{I_+^\perp}^*\\
    0 & 0 & 0 & 0
  \end{pmatrix},
\end{equation*}
where we use the same decomposition as before.

It follows that
\begin{equation*}
  e^{-t B^2_{1}} \xrightarrow{t\to \infty} P_+ \oplus \pr,
\end{equation*}
where $P_+$ is the even part of the operator, and $\pr$ the odd
part.
 For the
supertraces, this implies
\begin{equation*}
  \STR (e^{-t B^2_{1}}) \xrightarrow{t\to\infty} \TR (P_+) - \TR(\pr) \in A_{{\rm ab}}.
  \end{equation*}
  Here, the trace is always taken  in the sense of integration over
  the diagonal. Note that $\pr$ is a
  smoothing operator by Lemma \ref{lem:pr_smoothing}.
\end{lemma}

\begin{lemma}
  \begin{equation*}
    \int_M \tr^{\rm alg}_x \pr(x,x) = \int_M \tr^{\rm alg}_x P_-(x,x) \in A_{{\rm ab}}.
  \end{equation*}
\end{lemma}
\begin{proof}
  By the proof of Lemma \ref{lem:pr_smoothing} (using Lemma
  \ref{lem:inverse_of_one_plus_smoothing}),
  \begin{equation*}
    \pr = (1-P)^{-1} P_- (1- P)=(1-Q) P_- (1-P),
  \end{equation*}
  where $P$ and $Q$ are smoothing operators. The assertion follows
  from the trace property for the integral trace for smoothing operators.
\end{proof}

\begin{lemma}\label{lem:independece_of_w_alpha}
$\STR(e^{-tB^2_{\alpha}})$ is independent of $\alpha$.
\end{lemma}
\begin{proof}
  The independence follows, as before, from Duhamel's
  formula.
\end{proof}

We are now in the position of proving Proposition \ref{prop:alg=int},
which is stated once again here for the convenience of the reader:

\begin{corollary}\label{cor:alg=int}
  The algebraic trace $[I_+]_{[0]}-[I_-]_{[0]}$ of $[I_+]- [I_-]\equiv \Ind (D_{\L})$,
  i.e.~the image under the induced map
$  \tr^{{\rm alg}}\colon K_0(A)\to A_{{\rm ab}}$
of the index class $\Ind (D_{\L})$, can be calculated as
  \begin{equation*}
  [I_+]_{[0]}-[I_-]_{[0]} =   \int_M \tr^{\rm alg} P_+(x,x) - \int_M \tr
  P_-(x,x)\equiv \TR P_+ - \TR P_- \in
  A_{{\rm ab}},
  \end{equation*}
  where $P_+$ and $P_-$ are the projections onto $I_+$ and $I_-$ as
  given by the Mishchenko-Fomenko decomposition.
\end{corollary}
\begin{proof}
  Both expressions are limits for $t\to \infty$ of
$    \STR (e^{-tD_{\L}^2})$ by Lemma
\ref{lem:independece_of_w_alpha} and Lemma
\ref{lem:independence_of_alpha_no_w}.
\end{proof}

\begin{remark}
  Note that the proof shows that $[I_+]_{[0]}-[I_-]_{[0]} $ can also
  be expressed with many of the other integral operators we have used
  throughout the proof. In fact it might be useful to observe that
  the following general proposition holds:
  \begin{proposition}\label{index-via-parametrix}
Let $D_{\L}$ as above and let $\Q\in \Psi^{-1}_A (M,E^-\otimes
\L,E^+\otimes \L)$ be a parametrix for $D_{\L}^+$ with remainders
$S_\pm\in \Psi^{-\infty}_A$. Then $$\ind_{[0]}(D_{\L}) = \TR
(S_+)-\TR (S_-)\,.$$
\end{proposition}
The proof of the proposition is standard.
\end{remark}

Having analyzed the large-time behavior of the $A_{{\rm ab}}$-valued
supertrace of the heat-kernel, we now turn our attention to the short-time behavior.

\begin{lemma}\label{lem:short time closed}
  The local supertrace
  \begin{equation*}
\str^{\rm alg}_x (e^{-tD_{\L}^2}(x,x))\vol(x)
\end{equation*}
has a limit for $t\to 0$
  which is precisely the differential form
  \begin{equation*}
(x\mapsto AS(D)(x)\wedge \ch(E)(x)\wedge \ch{\L}(x)_{[\dim
    M]}) \in \Omega^{\dim M} (M,A_{{\rm ab}}),
\end{equation*}
where the Chern forms are
  defined as usual using Chern-Weyl theory and the curvature of the
  connections.
\end{lemma}
\begin{proof}
  We use Getzer's proof, as in the book by Berline, Getzler and Vergne \cite{BGV}.
\end{proof}

\begin{remark}\label{rem:short-time}
In fact, the same statement
  holds for the smoothing perturbation
  $D_{\L,\C}$ since, as already remarked,
  for the heat kernel $$ e^{-t
(D_{\L}+\C)^2}= e^{-t D_{M,\L}^2}+ t
C^{\infty}([0,\infty),\Psi^{-\infty}_A(M,E\otimes \L,E\otimes
\L))\,$$ (a consequence of Duhamel's formula together with the
fact that $\C\in\Psi_A^{-\infty}(M,E\tensor \L,E\tensor \L)$).
\end{remark}

\begin{lemma}
  \begin{equation*}
    \frac{d}{dt} \STR(e^{-tD_{\L}^2}) = 0 \;\;\text{in}\;\;A_{{\rm ab}}.
  \end{equation*}
\end{lemma}
\begin{proof}
This is, once again,  a consequence of Duhamel's formula.
\end{proof}

We can finally give a complete proof of
the $A_{{\rm ab}}$-valued Atiyah-Singer index theorem
\ref{theo:general_AS}:

\begin{proof} Only the index formula itself remains to be established.
 We integrate from $0<\epsilon<1$ to
$1/\epsilon$  the derivative $\frac{d}{dt} \STR(e^{-tD_{\L}^2})$;
we apply the fundamental theorem of calculus, we let
$\epsilon\downarrow 0$ and use the large and short time behavior
of $\STR(e^{-tD_{\L}^2})$ as in \ref{lem:short time closed} and in
\ref{cor:large time closed no alpha}.
\end{proof}

\subsection{Proof of the $A/\overline{[A,A]}$-valued
APS index formula}
\label{sec:aps-index-theory-proof}

In this subsection we shall recall and complement (a special case
of) the Atiyah-Patodi-Singer index theory developed by Leichtnam and
the first author in
\cite{LPGAFA} \cite{LP03}. In order to simplify the exposition we
shall only consider even dimensional manifolds with boundary.
Fundamental in our treatment is the extension of Melrose'
$b$-calculus to the $C^*$-algebraic setting. This is developed in
\cite{LPMEMOIRS}, \cite{LPGAFA} when the $C^*$-algebra $A$ is
equal to the reduced $C^*$-algebra of a discrete group; exactly
the same arguments work when $A$ is an arbitrary (unital)
$C^*$-algebra. We shall not enter into the precise definition of
the Mishchenko-Fomenko $b$-calculus with bounds
$\Psi^{*,\epsilon}_{b,A}$, $\epsilon>0$; we only recall that
operators in $\Psi^{*,\epsilon}_{b,A}$ are characterized by their
behavior, as distributions, on $W\times W$ or, more precisely, by
the behavior of their lifts on the so-called $b$-stretched
product $W\times_b W$ (the manifold with corners obtained by
blowing up $\partial W\times \partial W$ in $W\times W$); the
$b$-calculus with bounds is the sum of three spaces of operators
$$\Psi^{*,\epsilon}_{b,A}=
\Psi^{*}_{b,A}+\Psi^{-\infty,\epsilon}_{b,A}+\Psi^{-\infty,\epsilon}_{A}\,.$$
The first space on the right hand side is the {\it small}
$b$-calculus; it is an algebra and contains as a sub-algebra the
space of $b$-differential operators. Elements in the second space
are called (smoothing) boundary terms whereas the operators in the
third space are called {\it residual} and are directly
characterized on $W\times W$. The residual elements play the role
of the smoothing operators in the closed case. There are natural
composition rules for elements in the $b$-calculus with bounds.
Finally, elements in the $b$-calculus are bounded on suitable
$b$-Sobolev Hilbert modules: $$P\in
\Psi^{m,\epsilon}_{b,A}(W,E\otimes \L,F\otimes \L)\Rightarrow\;\;
P\colon H^\ell_b (W,E\otimes \L)\to H^{\ell-m}_b (W,F\otimes
\L)\;\;\text{is bounded}\,,$$ where the subscript $b$ in the
Sobolev modules indicates that these modules are defined using a
$b$-metric and $b$-differential operators. We set $H^\infty_b:=
\cap_{k\in\NN} H^k_b$. It is also important to consider weighted $b$-Sobolev
modules $x^\delta H^m_b$, with $x$ a boundary defining function;
in fact the inclusion $x^\delta H^{\ell+\epsilon}_b
\hookrightarrow H^\ell_b$ is an $A$-compact operator $\forall
\epsilon>0, \forall \delta>0$ .

\medskip
\noindent {\bf A short guide to the literature:} the basic
reference for the $b$-calculus and its applications to index
theory on manifolds with boundary is the book by Melrose
\cite{Melrose}. Short but rather complete introductions to the
theory are given in the surveys \cite{Mazzeo-Piazza} and
\cite{grieser} and also in the Appendix of \cite{MPI}. The
existence of spectral sections, and thus of trivializing
perturbations, for a $A$-linear Dirac-type operator $D_{M,\L}$ on
a closed manifold $M$ with vanishing index class in $K_{\dim M}
(A)$ (see \ref{eq:non_empty_perturb}) is established in
\cite{WuI} \cite{LP03}, building on work of Melrose-Piazza \cite{MPI}. The
$b$-calculus in the $C^*$-algebraic context (including $b$-Sobolev
modules) is studied in \cite{LPMEMOIRS}. If $M=\partial W$ and
$\C$ is a trivializing perturbation for $D_{\partial W, \L}$, then
the lifted perturbation $\C_W\in \Psi^{-\infty}_{b,A}$ is defined
in \cite{MPI} for families and in \cite{LPGAFA} in the
$C^*$-context. Using the $b$-calculus and the invertibility of
$D_{\partial W,\L}+\C$ one proves that the operator
$D_{W,\L}+\C_W$ is invertible modulo residual operators; since a
residual operator extends to a $A$-compact operator on $b$-Sobolev
modules we see that $D_{W,\L}+\C_W$ has a well defined index class
in $K_0 (A)$ (see below for more on this point).

\medskip
Our main interest is thus in the perturbed  $b$-operator
$D_{\L}+\C_W\in \Psi^1_{b,A}(M,E\otimes \L,E\otimes \L)$. We begin
by the existence of the heat kernel.

\begin{lemma}
  The  operator $H(t):=e^{-t{(D_{\L}+\C_W)^2}}$ is a smoothing operator in
  the (small)
Mishchenko-Fomenko $b$-calculus on $W$; $e^{-t(D_{\L}+\C_W)^2}\in
\Psi^{-\infty}_{b,A}$.
   The heat operator $H(\cdot)$ satisfies the heat equation with
  initial condition  $\lim_{t\downarrow 0} H(t)=\Id$.
\end{lemma}

\begin{proof}
 The result for the (unperturbed) $b$-differential operator $D_{\L}^2$
is obtained by employing the $b$-heat space as in \cite{Melrose}
and applying the same reasoning as in the proof of Lemma
\ref{lem:heat semigroup closed}. For the perturbed operator
$(D_{\L}+\C_W) ^2$ we apply the same arguments used in the proof
of Proposition 8 in \cite{MPI}.
\end{proof}

The following result is fundamental in developing a higher
APS-index theory

\begin{lemma}\label{lem:Mishchenko_decomposition_boundary}
There is a Mishchenko-Fomenko decomposition of the space of
sections
  of $E\tensor \L$ with respect to $D^+_{\L}+\C^+_W$, i.e.
  \begin{equation*}
    H_b^\infty(W,(E\tensor \L)^+) = \mathcal{I}_+ \oplus
    \mathcal{I}_+^{\perp},
 \quad H_b^\infty(W,(E\tensor
    \L)^-)=\mathcal{I}_- \oplus (D^+_{\L}+\C^+_W)(\mathcal{I}_+^{\perp}).
  \end{equation*}
In this decomposition, $\mathcal{I}_{\pm}\subset x^\epsilon
H_b^\infty$, with $\epsilon> 0$.
By completion, this
decomposition gives a decomposition of the Hilbert
  $A$-modules $H^m_b (W,(E\tensor \L)^{\pm})$, $m\in \NN$.

  The second decomposition is not a priori
  orthogonal. However, $D^+_{\L}+\C^+_W$ induces an isomorphism (in the Fr{\'e}chet
  topology) between $\mathcal{I}_+^{\perp}$ and
  $(D^+_{\L}+\C^+_W)(\mathcal{I}_+^{\perp})$. Moreover,
  $\mathcal{I}_+$ and $\mathcal{I}_-$
  are finitely generated projective Hilbert $A$-modules; the projections
  $\Pi_{\mathcal{I}_{+}}$   onto
  $\mathcal{I}_+$ (orthogonal) and  $\Pi_{\mathcal{I}_{-}}$   onto
  $\mathcal{I}_- $ (along
  $(D^+_{\L}+\C^+_W)(\mathcal{I}_+^{\perp})$) are residual,
i.e.~belong to $\Psi^{-\infty,\epsilon}_{A}(W; E\tensor \L,
E\tensor \L)$. $\mathcal{I}_{\pm}$ are already
  complete finitely generated projective Hilbert $A$-modules,
  i.e.~unchanged when passing to any  completion $H^m_b (W,E\tensor \L)$.
 The  index class associated to $D_{\L}^++\C_W^+$, denoted $\Ind_b (D_{\L},\mathcal{C})$,
 is by definition
  \begin{equation}\label{explicit b-index}
\Ind_b (D_{\L},\mathcal{C}):=[\mathcal{I}_{+}] - [\mathcal{I}_{-}]
\,\in\, K_0 (A)
  \end{equation}
\end{lemma}

\begin{proof}
The Lemma is proved in \cite[Appendix B]{LPMEMOIRS}  for an operator
$D_{\L}$ with invertible boundary operator. As explained in
\cite{LPGAFA}, Theorems 6.2 and 6.5,  the same proof applies to
perturbed operators, such as $D_{\L}+\C_W$, with invertible
indicial family.
\end{proof}

\begin{remark}
Suppose now that
 $\mathcal{C}$ is defined from a spectral section
$\mathcal{P}$ associated to $D_{\partial W, \mathcal{L}}$ : thus
$\mathcal{C}\equiv \mathcal{C}_{\mathcal{P}}$ for some
$\mathcal{P}$. Then $$ \Ind_b ({D}_{\L},
\mathcal{C}_{\mathcal{P}})=\Ind_{{\rm
    APS}}({D}_{\L},\mathcal{P})\;\;\text{in}
\;\;K_0 (A)$$
where on the right-hand side a suitable generalization of the
Atiyah-Patodi-Singer boundary condition appears. See \cite{WuI},
\cite{LP03}
\end{remark}

Recall that our goal is to prove an index formula for $$\tr^{{\rm
alg}}(\Ind_b (D_{\L},\C)):= \Ind_{b,[0]}
(D_L,\mathcal{C})\,\equiv\, [\mathcal{I}_{+}]_{[0]} -
[\mathcal{I}_{-}]_{[0]} \in A_{{\rm ab}}\,.$$

\subsubsection{The algebraic index perturbation}
\label{sec:algebr-index-pert}

\begin{definition}
  Instead of $W$ consider $W_+$, the disjoint union
  of the manifold with an additional point. The Dirac type operator on
  the additional point is by definition the $0$ operator on
  $\complexs$. Note that the boundary of $W_+$ is still $M$. We define
  $L_+$, the twisting bundle on $W_+$ to be
  $L\disjointunion (\mathcal{I}_+\oplus \mathcal{I}_-)$, where we view
  $\mathcal{I}_{\pm}$ as an abstract finitely generated projective Hilbert
  $A$-module (therefore a bundle over the point).

  Now we define for $\alpha\in\reals$ the perturbed operators (of
  $D^+_{\L}+\C^+_W$, but we suppress $\L$ and $\C_W$ in the notation
  of the perturbation)
  \begin{equation*}
    D^+_\alpha :=
    \begin{pmatrix}
      (D_{\L}^+ +\C_W^+)|_{\mathcal{I}_+^\perp} & 0 & 0\\
      0 & (D_{\L}^+ +\C_W^+)|_{\mathcal{I}_+} &  \alpha\\
      0 & \alpha & 0
    \end{pmatrix},
  \end{equation*}
  and  $D^-_\alpha$ as the (formal) adjoint of
  $D^+_{\alpha}$.
  The description is with respect to the
  b-Mishchenko-Fomenko decomposition of Lemma \ref{lem:Mishchenko_decomposition_boundary},
  \begin{align*}
H_b^\infty(W_+,(E\tensor L)_+^+)=& \mathcal{I}_+^{\perp}\oplus
\mathcal{I}_+
  \oplus \mathcal{I}_- \\  H_b^\infty(W_+,(E\tensor L)_+^-)=&
  (D^+_{\L}+\C_W^+)(\mathcal{I}_+^{\perp})\oplus \mathcal{I}_-
  \oplus \mathcal{I}_+
\end{align*}
Note that $\mathcal{I}_{\pm}$ has two roles here: first as subset of
the section of $E\tensor L$ on $W$, secondly as possible values of the
sections at the additional point $*$, where the fiber is
$\mathcal{I}_-\oplus\mathcal{I}_+$.
\end{definition}

\begin{lemma}
 For each $\alpha\in\reals$, $D_\alpha$ is a bounded
perturbation
  of $D_0$. More precisely $D_\alpha-D_0$ belongs to the residual
  space
  $\Psi^{-\infty,\epsilon}_{A}(M_+, (E\otimes L)_+, (E\otimes L)_+)$, for some $\epsilon> 0$.
  \end{lemma}

  It follows by Duhamel expansion that $e^{-tD_\alpha^2}$ is defined for each $t>0$ and each
  $\alpha\in\reals$, and is an element in
  $\Psi^{-\infty,\epsilon}_{b,A} (M_+, (E\otimes L)_+, (E\otimes L)_+)$. In particular
  the $b$-supertrace
$\bSTR(e^{-tD_\alpha^2})\in A_{{\rm ab}}$ is defined. Let
us recall the definition of the $b$-supertrace, from
\cite{Melrose}. First we
 recall the definition of the regularized integral on the manifold $W$ endowed
 with a product $b$-metric. Let ${\rm vol}_{b,W}$ be the associated
 volume form.
 Let us fix once and for all
a  trivialization $\nu \in C^\infty(\partial W, N_+\partial W)$ of
the inward pointing normal bundle to the boundary and let $x \in
C^\infty(W)$ a boundary defining function for $\partial W$ such
that $dx\cdot\nu=1$ on $\partial W$. For any function $f \in
C^\infty(W)$ we set
\def  \nuint {\raise10pt\hbox{$\nu$}\kern-6pt\int}
\begin{equation*}
\nuint_{W} f\,:=\,
 \lim_{\epsilon \rightarrow 0^+} [ \int_{x> \epsilon}
f\, {\rm vol}_{b,W} \; +\;\log \epsilon \int_{\partial W }
f|_{\partial W}\, {\rm vol}_{\partial W} ]
\end{equation*}
 Consider now an element  $K\in \Psi^{-\infty,\epsilon}_{b,A}$
and its restriction to the lifted diagonal $\Delta_b\subset
M\times_b M$. Using the identification $\Delta_b\equiv W$  and
taking the $A_{{\rm ab}}$-valued
  supertrace of the endomorphism of $(E\tensor L)_x$, denoted $\str_x$, we define
  \begin{equation*}
\bSTR_\nu ( K) \,:=\, \nuint_{W} \str^{\rm alg}_x K|_{\Delta_b} (x)\,.
\end{equation*}
One can easily compute the effect of changing the trivialization
$\nu$; from now on we shall suppress the subscript $\nu$ in the
notation of the $b$-supertrace. Needless to say, in the ungraded
case we can similarly define the $b$-trace. The $b$-trace does not
vanish on commutators; however, there is an explicit formula due
to Melrose, see \cite{Melrose}, for computing the defect. The
formula involves the indicial families of the two operators:

\begin{proposition}\label{prop:commutator for b trace}
For any $K,K^\prime$ belonging to $ \Psi^{-\infty,
\epsilon}_{b,A}$ one has: \begin{equation}\label{commutator for b
trace} \bSTR [K,K^\prime]\,=\, \frac{\sqrt{-1}}{2\pi} \int_\RR
\STR \bigl(\, \frac{\partial|}{\partial \lambda}I(K,\lambda)\circ I(K^\prime,\lambda)\,\bigr)
d\lambda. \end{equation}
 If we replace $K$ by a differential
operator in ${\rm Diff}^*_{b, A}$ and $K^\prime$ by the
composition of $K^\prime$ with an element of the calculus with
bounds $\Psi^{*,\epsilon}_{b,A}$ then the same commutator formula
is valid.
\end{proposition}

\begin{lemma}
The following formula is clear from the definition:
$$\bSTR(e^{-tD_0^2}) = \bSTR_W(e^{-t(D_{\L}+\C_W)^2}) +
  [I_-]_{[0]}-[I_+]_{[0]} \in
  A_{{\rm ab}}\,.$$
\end{lemma}

\begin{lemma}
  For $\alpha$ sufficiently large, $D_\alpha$ is invertible in the MF $b$-calculus
  with bounds, and
  therefore
  \begin{equation*}
  \bSTR(e^{-tD_\alpha^2})\xrightarrow{t\to\infty} 0.
\end{equation*}
\end{lemma}

\begin{proof}
Directly from the $b$-Mishchenko-Fomenko decomposition
theorem we prove the invertibility of $D_\alpha$ , $\alpha$ large,
exactly as in the closed case, see \cite{LottI}, Section VI.
\end{proof}

\begin{lemma}\label{lem:independence_of_alpha_no_w_with_boundary}
\begin{equation*}
  \bSTR(e^{-tD_\alpha^2}) - \bSTR(e^{-tD_0^2}) = F(\alpha,t) \quad (=0),
\end{equation*}
where for each fixed $\alpha$ and $t$, $F(\alpha,t)= 0$ in
$A_{{\rm ab}}$.
\end{lemma}

\begin{proof}
  Use Duhamel's formula to compute the derivative with respect to
  $\alpha$ of $\bSTR(e^{-tD_\alpha^2})$. The usual calculations show
  that this is the $b$-supertrace of a supercommutator.
  Using \ref{commutator for b
trace} we see that this is an
  explicit term $F(\alpha,t)$, localized on the boundary.
  As explained in \cite{LPMEMOIRS}, formula (14.16),
  one can show that $F(\alpha,t)=0$ in $A_{{\rm ab}}$; this is a consequence of the particular
  structure of the two projections onto $\mathcal{I}_\pm$, namely, that they are residual.
\end{proof}

\begin{corollary}\label{cor:large-b}
  \begin{equation*}
 \lim_{t\to\infty}\bSTR_W (e^{-t(D_{\L}+\C_W)^2})
= [I_+]_{[0]} - [I_-]_{[0]}\equiv \Ind_{b,[0]} (D_{\L},\C)\,,
\end{equation*}
where part of the assertion is that the limit exists.
\end{corollary}

\subsubsection{The integral operator $b$-index}
\label{sec:integral-index-APS}

As in the closed case, we wish to connect the index $\Ind_{b,[0]}
(D_{\L},\C)\in A_{{\rm ab}}$ defined using the index class
and the algebraic trace  $\tr^{{\rm alg}}:K_0 (A) \to
A_{{\rm ab}}$ to the integral-kernel-trace, $\TR$, of the
projection operators given by the $b$-Mishchenko-Fomenko
decomposition.

\begin{definition}
We set $P_+ := \Pi_{\mathcal{I}_+}$, $P_- :=
\Pi_{\mathcal{I}_-}$.
  We define a second smoothing perturbation of $(D_{\L}+\C_W)^+$ by
  \begin{equation*}
B^+_{\alpha}:=
  (D^+_{\L}+\C_W^+) -\alpha 
P_- (D^+_{\L}+\C^+) P_+\,, 
\end{equation*}
for $\alpha\in\reals$; define
  $B^-_{\alpha}$ as the (formal) adjoint of $B^+_{\alpha}$.
\end{definition}

\begin{remark}
  Observe that $$B_{\alpha}=\begin{pmatrix} 0& B^-_\alpha\\
                    B^+_\alpha & 0 \end{pmatrix}
  $$ is a residual perturbation of $D_{\L}+\C_W$ in the $b$-calculus.
We shall also use ${\rm pr}:= \text{orthogonal projection
onto}\,\, \im P_-^* $.
\end{remark}

\begin{lemma}\label{lem:independece_of_w_alpha for boundary}
$\STR(e^{-tB^2_{\alpha}})$ is independent of $\alpha$.
\end{lemma}
\begin{proof}
  The result follows from Duhamel's formula, the formula
  for
  the $b$-supertrace of a commutator and the fact that the
  perturbation involved in the definition of $B_\alpha$ is
  residual.
\end{proof}

Proceeding exactly as in the closed case one  proves the analog
of  lemma \ref{lem:pr_smoothing} and, in particular, that $\pr$ is
a residual operator; the analog of lemma
\ref{lem:inverse_of_one_plus_smoothing} is established using the
semi-ideal property of the residual operators (see \cite{Melrose},
formula (5.23) and Proposition 5.38). Then, as in the closed case,
\begin{equation*}
  \bSTR (e^{-t B^2_{1}}) \xrightarrow{t\to\infty} {\rm bTR} (P_+) - {\rm bTR}(\pr)=
  \TR(P_+)-\TR(\pr) \in A_{{\rm ab}},
  \end{equation*}
with the last equality following once again from the fact that
$P_+$ and $\pr$ are residual.
  Finally, from the trace property we get $\TR(\pr)=\TR(P_-)$,
  exactly as in the closed case.
  Using lemma \ref{lem:independece_of_w_alpha for boundary}
we finally get

\begin{corollary}\label{cor:alg=int:b-case}
  The $A_{{\rm ab}}$-valued $b$-index $\Ind_{b,[0]} (D_{\L},\C)$,
i.e.~the image under the induced map
$\tr^{{\rm alg}}\colon K_0(A)\to A_{{\rm ab}}$
of the index class $\Ind_b (D_{\L},\C)$, can be calculated as
  \begin{equation*}
 \int_M \tr^{\rm alg} P_+(x,x) - \int_M \tr^{\rm alg}
  P_-(x,x) \in
  A_{{\rm ab}},
  \end{equation*}
  where $P_+$ and $P_-$ are the projections onto $\mathcal{I}_+$ and $\mathcal{I}_-$ as
  given by the $b$-Mishchenko-Fomenko decomposition.
\end{corollary}
\begin{proof}
  Both expressions are limits for $t\to \infty$ of
$    \bSTR (e^{-t(D_{\L}+\C_W)^2})$ by Lemma
\ref{lem:independece_of_w_alpha for boundary} and Corollary
\ref{cor:large-b}.
\end{proof}

\subsubsection{Local expansion}
\label{sec:local-expansion_APS}

\begin{lemma}
  The local supertrace  $\str^{\rm alg}(e^{-t(D_{\L}+\C_W)^2}(x,x)) {\rm vol}_b(x)$
  has a limit for $t\to 0$
  which is exactly the differential form
  \begin{equation*}
(x\mapsto AS(D)(x)\wedge \ch(E)(x)\wedge \ch{\L}(x)_{[\dim
    w]}) \in \Omega^{\dim W} (W,A_{{\rm ab}}),
\end{equation*}
  where the Chern forms are
  defined as usual using Chern-Weyl theory and the curvature of the
  connections.
\end{lemma}
\begin{proof}

The formula holds for $\str(e^{-tD_{\L}^2}(x,x)) {\rm vol}_b(x)$,
with proof employing the rescaled $b$-heat calculus. See
\cite{Melrose}, Chapter VIII. For the perturbed  operator we
simply observe that $$ e^{-t (D_{\L}+\C_W)^2}= e^{-t D_{\L}^2}+ t
C^{\infty}([0,\infty),\Psi^{-\infty}_{b,A})\,,$$ See \cite{MPI},
Proposition 8.
\end{proof}

\subsubsection{The index formula}
\label{sec:defect}

The index formula, as stated in subsection
\ref{sec:aps-index-theory}, follows from the large and short time
limits for $\bSTR (e^{-t(D_{\L}+\C_W)^2})$ together with the
commutator formula \ref{prop:commutator for b trace}, with details
as in \cite{MPI}, \cite{LPGAFA}:

\begin{proposition}
  \begin{multline*}
  \lim_{t\to\infty} \bSTR(e^{-t(D_{\L}+\C_W)^2}) -   \lim_{t\to 0}
  \bSTR(e^{-t(D_{\L}+\C_W)^2}) \\
 = \int_0^\infty \frac{d}{dt}
  \bSTR(e^{-t(D_{\L}+\C_W)^2}) =
   -\frac{1}{2}\eta_{[0]}(D_{M,\L}+\C) \in A_{{\rm ab}}.
 \end{multline*}
\end{proposition}

Notice that the derivative with respect to $t$ of $
\bSTR(e^{-t(D_{\L}+\C_W)^2})$, as given by the commutator formula,
is not precisely the eta integrand $$
-\frac{1}{2}\frac{1}{\sqrt{\pi}} t^{-1/2} \TR \left (
(D_{M,\L}+\C)e^{-t (D_{M,\L}+\C)^2}\right).$$ However, as explained in
\cite{MPI}, its integral from $0$ to $+\infty$ is indeed equal to
$$-\int_0^\infty \frac{1}{2}\frac{1}{\sqrt{\pi}} t^{-1/2} \TR
\left ( (D_{M,\L}+\C)e^{-t (D_{M,\L}+\C)^2}\right)dt\,\equiv
-\frac{1}{2}\eta_{[0]}(D_{M,\L}+\C) \in A_{{\rm ab}}$$
 The
$A_{{\rm ab}}$-valued APS index theorem, as stated in
\ref{theo:general_APS}, now follows immediately from the above
results.


\section{Graded Hermitian complexes and the signature operator}
\label{sec:appendix-odd-sign}

In this Appendix we shall make precise our conventions for the
signature operator
We 
follow 
\cite[Section 3.1]{HiSka} and also  \cite{LLP}, giving the general definition of
graded regular $n$-dimensional Hermitian complex and its associated 
signature operator.
The only difference between our conventions and those of \cite{HiSka}
is that we deal with left modules, whereas \cite{HiSka} deals with
right
modules.  The following material is taken directly from \cite{LLP}.

Let $A$ be a $C^*$-algebra with unit.

\begin{definition}\label{Hcomplex}
A graded regular $n$-dimensional Hermitian complex consists of the
following data
\begin{itemize}
\item  A $\ZZ$-graded cochain complex $({\cal E}^*, d)$ of finitely-generated
projective left $A$-modules,
\item  A nondegenerate quadratic form
$Q : {\cal E}^* \times {\cal E}^{n-*} \rightarrow A$ and
\item  An operator $\tau \in Hom_{A} \left( {\cal E}^*,
{\cal E}^{n-*} \right)$ 
such that  $Q(bx, y) = b Q(x,y)$,  $Q(x,y)^* = Q(y,x)$,  $Q(dx, y) +
Q(x, dy) = 0$, $\tau^2 = I$, and
 $<x, y> := Q(x, \tau y)$ defines a Hermitian metric on ${\cal
  E}$.
\end{itemize}
\end{definition}

Let $M$ be a closed oriented $n$-dimensional Riemannian
manifold. Let ${\cal V}$ be a flat
$A$-vector bundle on
$M$.
We assume that the fibers of
${\cal V}$ have
$A$-valued Hermitian inner products which are compatible with
the flat structure.

Let $\Omega^*(M; {\cal V})$ denote the vector space of smooth
differential forms with coefficients in ${\cal V}$. If $n=\dim(M) > 0$ then
$\Omega^*(M; {\cal V})$ is not finitely-generated over
$A$, but we wish to show that it still has all of
the formal properties of a graded regular $n$-dimensional Hermitian complex.
If $\alpha \in
\Omega^*(M; {\cal V})$ is homogeneous, denote its degree by $|\alpha|$.
In what follows, $\alpha$ and $\beta$ will sometimes implicitly denote
homogeneous elements of $\Omega^*(M; {\cal V})$.
Given $m\in M$ and $(\lambda_1 \otimes e_1), (\lambda_2 \otimes e_2)
\in \Lambda^*(T^*_mM) \otimes
{\cal V}_m$, we define
$(\lambda_1 \otimes e_1) \wedge (\lambda_2 \otimes e_2)^*
\in \Lambda^*(T^*_mM) \otimes
A$ by
$$(\lambda_1 \otimes e_1) \wedge (\lambda_2 \otimes e_2)^* =
(\lambda_1 \wedge \overline{\lambda_2}) \otimes <e_1, e_2>.$$
Extending by linearity (and antilinearity), given
$\omega_1, \omega_2 \in \Lambda^*(T^*_mM) \otimes
{\cal V}_m$,
we can define
$\omega_1 \wedge \omega_2^* \in
\Lambda^*(T^*_mM) \otimes A$.

Define an $A$-valued
quadratic form $Q$ on $\Omega^*(M; {\cal V})$ by
$$
Q(\alpha, \beta) = i^{- |\alpha| (n - |\alpha|)} \int_M \alpha(m)
\wedge {\beta}(m)^*.
$$
It satisfies $Q(\beta, \alpha) = {Q(\alpha, \beta)}^*$.
Using the Hodge duality operator $*$, define
$\tau : \Omega^p(M; {\cal V}) \rightarrow
\Omega^{n-p}(M; {\cal V})$ by
$$
\tau (\alpha) = i^{- |\alpha| (n-|\alpha|)} * \alpha.
$$
Then $\tau^2 = 1$, and the inner product $<\cdot,\cdot>$ on
$\Omega^*(M; {\cal V})$ is
given by
$< \alpha, \beta> = Q(\alpha, \tau \beta)$.
Let $d_{{\cal V}}$ be the de Rham differential with values in the flat
bundle ${\cal V}$; define $d : \Omega^*(M; {\cal V}) \rightarrow
\Omega^{*+1}(M; {\cal V})$ by
\begin{equation}\label{signed-diff}
d \alpha = i^{|\alpha|} d_{{\cal V}} \alpha.
\end{equation}
It satisfies $d^2 = 0$.
Its dual $d^\dagger$ with respect to $Q$, i.e. the operator $d^\dagger$
such that $Q(\alpha, d\beta) = Q(d^\dagger\alpha, \beta)$, is given by
$D^\dagger= -D$. The formal adjoint of $D$ with respect to $<\cdot, \cdot>$ is
$D^* = \tau D^\dagger \tau = - \tau D \tau$.
\begin{definition}\label{sigop}
If $n$ is even, the signature operator is
\begin{equation}\label{evensign}
\D^{{\rm sign}} := d + d^* = d - \tau d \tau.
\end{equation}
It is formally self-adjoint and
anticommutes with the $\ZZ_2$-grading
operator $\tau$.  If we denote by $\Omega^\pm (M,\V)$ the $\pm
1$-eigenspaces
of $\tau$ then $$\D^{{\rm sign}}=\begin{pmatrix}0&\D^{{\rm sign}}_-\\
\D^{{\rm sign}}_+ & 0 \end{pmatrix} $$
The index class of the signature operator in $K_0 (A)$ is, by definition, the index
class of the elliptic operator $\D^{{\rm sign}}_+ $.\\
If $n$ is odd, the signature operator is
\begin{equation}\label{oddsign}
\D^{{\rm sign}} = -i (d \tau + \tau d).
\end{equation}
It is formally self-adjoint and defines an index class in $K_1 (A)$.
\end{definition}

\smallskip
Now suppose that $M$ is a compact oriented
manifold of dimension $n=2m$ with boundary $\partial M$.
We fix a non-negative boundary defining function
$x\in C^\infty(M)$ for $\pa M$ and a Riemannian metric
on $M$ which is isometrically a product in an (open) collar neighborhood
${\cal{U}}\equiv (0,2)_x\times\pa M$ of $\partial M$.
The signature operator $\Si$ is a well-defined differential
operator; it is  associated to the graded regular n-dimensional
Hermitian complex defined by $\Omega^*_c({\rm int}(M); {\cal V})$
and the Riemannian metric on $M$.
Let ${\cal V}_0$
denote the pullback of
${\cal V}$ from $M$ to $\partial M$; there is a natural
isomorphism
$${\cal V}|_{{\cal{U}}}\cong (0,2) \times {\cal V}_0 \,.$$
 Let $Q_{\partial M}$, $\tau_{\partial M}$,
$d_{\partial M}$
and $\Si (\partial M)$
denote the expressions defined above on
$\Omega^*(\partial M; {\cal V}^\infty_0)$.
One can decompose $Q$, $\tau$, $d$ and $\Si$, when restricted to
compactly-supported forms on $(0,2) \times \partial M$, in terms of
$Q_{\partial M}$, $\tau_{\partial M}$, $D_{\partial M}$
and $\Si (\partial M)$.
This computation is given in great detail in \cite{LLP}; one proves in
particular that with the definitions given above {\it the 
operator  $\Si(\pa M)$ is the boundary operator of
$\Si_+$ in the sense of Atiyah-Patodi-Singer [APS, (3.1)]}.

\smallskip
\begin{remark}
Let $(X,g)$ be closed oriented and of  dimension $2m-1$.  
Assume that $\V=X\times\complexs$ is the trivial complex line bundle.
We wish to compare the eta-invariant of the signature operator
defined with our sign-conventions, denoted $D^{{\rm sign}}$, and the eta-invariant 
of the operator $B$ appearing in the work of Atiyah-Patodi-Singer, see
\cite{MR53:1655b}.
In order to simplify the notation we set $D:=D^{{\rm sign}}$

Recall that $B\phi=(\sqrt{-1})^{m}(-1)^{p+1}(\epsilon \star d - d \star)\phi$
with $\epsilon=1$ if $\phi\in\Omega^{2p}(X)$
and $\epsilon=-1$ if $\phi\in\Omega^{2p-1}(X)$.
The operator $D$ is given instead by $D\phi=(-1)^m ((\sqrt{-1} \star d
-d\star)) \phi$ if $\phi \in\Omega^{2p}(X)$ and $D\phi=- \star d\phi
-\sqrt{-1}d\star \phi$ if $\phi \in\Omega^{2p-1}(X)$.
Each operator preserves the parity of the form-degree. Moreover, both commute with the 
following self-adjoint odd involution: $\Xi:=\eta\Theta$ with
$\Theta=(-1)^p\star \; {\rm on}\; \; {\rm both} \;
\Omega^{2p} \;{\rm and} \; \Omega^{2p-1}$ and 
 $\eta=1$ if
$m=2k$; $\eta=\sqrt{-1}$ if $m=2k-1$.
Thus $B=B^{{\rm even}}\oplus B^{{\rm odd}}$, $D=D^{{\rm even}}\oplus
D^{{\rm odd}}$,
with $B^{{\rm even}}=(\Xi)^{-1} B^{{\rm odd}}\Xi$ and similarly for
$D$ \footnote{Notice that there is a sign-mistake in \cite{LPAGAG}:
the $\Theta$ given there in Section 1 must be replaced with the $\Xi$ given here}.
The Hodge theorem implies the following orthogonal
decomposition of the space of differential forms on $X$, where $H$
stands for the space of harmonic forms:
$$\Omega^* =\Omega^0\oplus {\rm \Omega}^1\oplus\cdots
d{\rm \Omega}^{m-2}\oplus  {\bf  H^{m-1}}\oplus
{\bf  d^* \Omega^{m}}\oplus {\bf d\Omega^{m-1}}\oplus{\bf  H^{m}}
\oplus d^*\Omega^{m+1}\oplus\cdots
\oplus {\rm \Omega}^{2m-1}
$$
Consider now the following two subspaces of $\Omega^*$: 
$$ V=H^{m-1}\oplus d^*\Omega^{m}\oplus d\Omega^{m-1}\oplus H^m\,;$$
$$W= \Omega^0\oplus {\rm \Omega}^1\oplus\cdots
{\rm d\Omega}^{m-2}\oplus d^*\Omega^{m+1}\oplus\cdots
\oplus {\rm \Omega}^{2m-1}$$
Both operators $B$ and $D$ are block diagonal with respect to the
orthogonal decomposition $\Omega^*=V\oplus W$; 
we denote by $D_V, B_V$ and $D_W, B_W$ their restrictions.
Notice that 
both $V$ and $W$ are invariant under the action of $\Xi$
so that, as above,  $B_V=B_V^{{\rm even}}\oplus B^{{\rm odd}}_V$,
$B_W=B^{{\rm even}}_W\oplus B^{{\rm odd}}_W$ and similarly
for $D_V$ and $D_W$.
It's easy to check from the explicit expression given above
that $B^{{\rm even}}_V=-D^{{\rm even}}_V$.
We will show that 
\begin{equation}\label{aps-ours-vanish}
\eta(B_W)=0=\eta(D_W)
\end{equation}
which implies immediately that
\begin{equation}\label{aps-ours}
\eta(B)=-\eta(D)
\end{equation}
In order to establish (\ref{aps-ours-vanish})
we argue as follows.
Define
$$\Omega^{<}=\Omega^0\oplus {\rm \Omega}^1\oplus\cdots
\Omega^{m-2}\oplus{\rm d\Omega}^{m-2}\;;\;\;\quad
\Omega^{>}=d^*\Omega^{m+1}\oplus
\Omega^{m+1}\cdots\oplus {\rm \Omega}^{2m-1}$$
so that $W=\Omega^{<}\oplus \Omega^{>}$.
There is a natural  involution $\alpha$ on $W$ defined as follows:
$$\alpha={\rm Id} \;\;  {\rm on} \;\; \Omega^{<} \quad\quad
\quad \alpha=-{\rm Id} \;\; {\rm on} \; \; \Omega^{>}$$
It is immediate from the structure of $D$ and $B$ that
$$D_W \circ \alpha+\alpha\circ D_W=0\,,\quad B_W \circ \alpha+\alpha\circ B_W=0$$
In other words $\alpha$ gives a grading to $W$,
$W^+= \Omega^{<}$, $W^- = \Omega^{>}$
and both $D_W$ and $B_W$ are {\it odd} with respect to
such a grading; thus $\eta(B_W)=0=\eta(D_W)$ as required.\\
{\it Conclusion:} the odd-signature operator considered here and the
odd signature operator considered in the work of Atiyah-Patodi-Singer
are different. However, their eta-invariants are equal up to a sign.

\end{remark}

\section{Pseudodifferential operators and morphisms of $C^*$-algebras}
\label{sec:appendix1}

Let $A$ and $B$ be unital $C^*$-algebras. We assume that there
exists a morphism of unital $C^*$-algebras $\lambda:A\to B$. Let
$\mathcal{E}_A\to M $ be a bundle of finitely generated projective
Hilbert $A$-modules, with fibers isomorphic to a fixed finitely
generated projective Hilbert $A$-module $V$. In particular, there
exists $N\in \NN$ and $W$, a Hilbert $A$-module, such that
$V\oplus W= A^N$ as Hilbert $A$-modules; moreover $V$ is the image
of a self-adjoint projection $p\in M(N\times N,A)$, $p=(p(ij)),
i,j=1,\dots,N$. The morphism $\lambda$ induces in a natural way a
morphism of matrix algebras $M(N\times N,A)\to M(N\times N, B)$
that will be still denoted by $\lambda$. If we define $p_\lambda:=
(p_\lambda (ij))$ with $p_\lambda (ij) :=\lambda (p(ij))$, briefly
$p_\lambda=\lambda\circ p$,  then $p_\lambda$ is still a
self-adjoint projection in $M(N\times N,B)$ so that $\lambda
(V):=\Im p_\lambda$ is a finitely generated projective Hilbert
$B$-module. Applying this reasoning  to a more global situation,
we see that $\lambda$ induces  a bundle of finitely generated
projective Hilbert $B$-modules, $\mathcal{E}_B^\lambda$, with
fibers diffeomorphic to $\lambda (V)$: in fact, if $\mathcal{E}_A$
is defined by $p\in C^\infty (M, M_{N\times N} (A))$ with
$p=p^*=p^2$, then $\mathcal{E}_B^\lambda$ is simply defined by
$p_\lambda:=\lambda\circ p\in C^\infty (M, M_{N\times N} (B))$.
Alternatively, if $\{U_\alpha \}$ is a trivializing covering for
$\mathcal{E}_A\to M $, with transition function
$g^A_{\alpha,\beta}:U\alpha \cap U_\beta \rightarrow {\rm Iso}_A
(V)$ then $\mathcal{E}_B^\lambda\to M $ is defined by
$g^B_{\alpha,\beta}:=\lambda\circ g^A_{\alpha,\beta}: U_\alpha
\cap U_\beta \rightarrow {\rm Iso}_B (\lambda(V))$.

If $\mathcal{E}_A\to M $ is endowed with a connection $\nabla^A$,
then $\mathcal{E}_B\to M $ inherits in a natural way a connection
$\nabla^{B,\lambda}$ defined as follows.  If $\nabla^A$ is equal
to $p\circ d \circ p$ then $\nabla^{B,\lambda}$ is simply equal to
$p_\lambda \circ d \circ p_\lambda$; if, on the other hand,
$\nabla^A$ is arbitrary, then $\nabla^A=p\circ d \circ p +\omega$,
with $\omega\in C^\infty (M,{\rm End}_A (\mathcal{E}_A))$ and we
then define $$\nabla^{B,\lambda}= p_\lambda \circ d \circ
p_\lambda+\omega_\lambda$$ with $\omega_\lambda\in C^\infty
(M,{\rm End}_B (\mathcal{E}_B^\lambda))$ the 1-form induced by
$\omega$ and $\lambda$.
 Alternatively, if
$\nabla^A$ is defined by a collection of local 1-forms
$\{\omega_\alpha\}$ associated to the trivializing cover
$\{U_\alpha \}$, with $$\omega_\alpha\in \Omega^1 (U_\alpha,{\rm
End}_A (V))\,,\quad \omega_\beta= g_{\alpha,\beta}^{-1}
dg_{\alpha,\beta}+g_{\alpha,\beta}^{-1}\omega_\beta
g_{\alpha,\beta}$$ then we define the connection
$\nabla^{B,\lambda}$ through the local 1-forms
$\{\omega^\lambda_\alpha\}$ where, once again,
$\omega^\lambda_\alpha\in \Omega^1 (U, {\rm End}_B (\lambda (V)))
$ is defined in a natural way by $\omega_\alpha$ and the extension
of $\lambda$ to a morphism $M(N\times N,A)\to M(N\times N, B)$.

Consider now the
  graded vector space of pseudodifferential Hilbert $A$-module bundle operators
 $\Psi^*_A (M,\mathcal{E}_A, \mathcal{F}_A)$, with $\mathcal{F}_A$ a second bundle
 of finitely generated projective Hilbert $A$-modules. This space, defined
for the first time in \cite{Mish-Fom}, is nothing but
\begin{equation}
\Psi^* (M)\otimes_{C^\infty (M\times M)} C^\infty (M\times M, {\rm
HOM}_A (\mathcal{F}_A,\mathcal{E}_A))
\end{equation}
where we are considering the bundle $ {\rm HOM}_A
(\mathcal{F}_A,\mathcal{E}_A))\longrightarrow M\times M$ with
fiber at $(x,y)$ equal to ${\rm Hom}_A ((\mathcal{E}_A)_y,
(\mathcal{F}_A)_x)$. Here we are using the $C^\infty (M\times
M)$-module structure of both $\Psi^* (M)$ and $C^\infty (M\times
M, {\rm HOM}_A (\mathcal{F}_A,\mathcal{E}_A))$.
 The morphism $\lambda$ induces in a natural way a morphism of
graded vector spaces $\lambda_\natural: \Psi^*_A
(M,\mathcal{E}_A,\mathcal{F}_A)\to \Psi^*_B
(M,\mathcal{E}_B^\lambda,\mathcal{F}_B^\lambda )$ which is simply
induced by the natural map $$C^\infty (M\times M, {\rm HOM}_A
(\mathcal{F}_A,\mathcal{E}_A))\longrightarrow C^\infty (M\times M,
{\rm HOM}_B (\mathcal{F}_B^\lambda,\mathcal{E}_B^\lambda))\,.$$ If
$ \mathcal{F}_A =  \mathcal{E}_A$ then $\lambda_\natural$ is in
fact a morphism of graded algebras; we shall often denote
$\Psi^*_A (M,\mathcal{E}_A,\mathcal{E}_A)$ by $\Psi^*_A
(M,\mathcal{E}_A)$ or, more simply, by $\Psi^*_A$.  Thus the following
simple conclusion is true.

\begin{corollary}\label{corol:invertible_preserved}
  If $\mathcal{Q}\in \Psi^*_A$ is
  invertible as an element in $\Psi^*_A$, then $\lambda_\natural
  (\mathcal{Q})$ is invertible in $\Psi^*_B$.
\end{corollary}

As an important example we shall consider twisted Dirac operators.
Thus $\mathcal{E}_A= E\otimes \mathcal{V}_A$ with $E$ a vector
bundle over $M$ and $\mathcal{V}_A$ a line bundle with typical
fiber $A$ and endowed with a flat connection $\nabla^A$. We can
then choose trivializations of $\mathcal{V}_A$ with locally
constant transition functions $\{a_{\alpha,\beta}\in U_1
(A)\subset A\}$. Let $D\in \operatorname{Diff}^1 (M,E)$ be a
Dirac-type operator acting on the sections of $E$ and consider
$D_{\mathcal{V}_A}\in \operatorname{Diff}^1_A(M,\mathcal{E}_A)$,
the twisted Dirac operator defined by the connection $\nabla^A$.
Then $\lambda_\natural (D_{\mathcal{V}_A})=
D_{\mathcal{V}_B^\lambda}$ with $\mathcal{V}_B^\lambda$ the flat
bundle induced by $\lambda$.
 Notice that
the transition functions of the flat bundle
$\mathcal{V}_B^\lambda$ will only involve the image $\lambda (A)$.

Let now $\chi\in C^\infty_c (\RR,\RR)$. The operator
$D_{\mathcal{V}_A}$ defines a self-adjoint unbounded regular
operator on the (full) Hilbert A-module $L^2_A (M,\mathcal{E}_A)$.
It is well-known that there is a  continuous functional calculus
for regular operators on (full) Hilbert A-modules; in particular
$\chi (D_{\mathcal{V}_A})$ is a well defined bounded $A$-linear
endomorphism on $L^2_A (M,\mathcal{E}_A)$. As observed in
\cite{LPGAFA}, following the arguments in \cite[page 300]{Taylor},
one can show that $\chi (D_{\mathcal{V}_A})$ is in fact a
smoothing operator: $\chi( D_{\mathcal{V}_A})\in \Psi^{-\infty}_A
(M,\mathcal{E}_A)$. Since $\chi (D_{\mathcal{V}_A})$ is defined
through a functional integral, we see immediately that
\begin{equation}\label{chi-lambda=lambda-chi}
\lambda_\natural (\chi (D_{\mathcal{V}_A})= \chi (\lambda_\natural
(D_{\mathcal{V}_A}))= \chi ( D_{\mathcal{V}_B^\lambda})\,.
\end{equation}
The same is true for other functions and other operators, for the same
reason. E.g.~if $\mathcal{C}$ is a smoothing operator, then
\begin{equation}
  \label{eq:lambda_and_exp}
  \lambda_\sharp\left (
    (D_{\mathcal{V}_A}+\mathcal{C})
    \exp(t(D_{\mathcal{V}_A}+\mathcal{C})^2)\right) =
  (D_{\mathcal{V}_B^\lambda}+\lambda_\sharp\mathcal{C})\exp(t(D_{\mathcal{V}^\sharp_B}+\lambda_\sharp\mathcal{C})^2).
\end{equation}

\begin{remark}
Let $\phi\in C^\infty_c (\RR,[0,1])$ be any real function equal to 1
on $[-\epsilon, \epsilon]$ and equal to 0 on
$(-\infty,-2\epsilon]\cup [2\epsilon,\infty)$. Let
$\mathcal{C}_A\in \Psi^{-\infty}_A$ and let
$\mathcal{C}_B^\lambda:= \lambda_\natural (\mathcal{C}_A)$. Then
the following implication holds:
\begin{equation}\label{localize-at-0}
\text{if}\;\;(\Id-\phi(D_{\mathcal{V}_A}))\circ
\mathcal{C}_A=0\;\;\text{then}\;\;(\Id-\phi(D_{\mathcal{V}_B^\lambda}))\circ
\mathcal{C}_B^\lambda=0\,.
\end{equation}
 Notice that if
$(\Id-\phi(D_{\mathcal{V}_A}))\circ \mathcal{C}_A=0$ then, taking
adjoints, $\mathcal{C}_A \circ (\Id-\phi(D_{\mathcal{V}_A}))=0\,.$
\end{remark}

Let now $K\in \Psi^{-\infty}_A(M,\mathcal{E}_A)$ be a smoothing
operator, and let us consider ${\rm TR}(K)\in
A_{{\rm ab}}$. As an immediate consequence of the definitions,
\begin{equation}\label{TR-commutes-with-lambda}
\lambda ({\rm TR}(K))={\rm TR}(\lambda_\natural K)\in
B_{\rm ab}.
\end{equation}
If now $\mathcal{C}_B^\lambda=\lambda_\natural
(\mathcal{C}_A)$, then
\begin{lemma}\label{lem:stable_eta_preserved}
\begin{equation}\label{lambda-eta=eta-lambda}
\lambda(\eta_{[0]}(D_{\mathcal{V}_A}+\mathcal{C}_A))=\eta_{[0]}(D_{\mathcal{V}_B^\lambda}+\mathcal{C}_B^\lambda))
\;\;\text{in}\;\;B_{\rm ab}
\end{equation}
\end{lemma}
\begin{proof}
  This simply follows from the definition of $\eta_{[0]}$ as a
  convergent integral. Then, using for the second equality
  \eqref{TR-commutes-with-lambda},
  \begin{equation*}
    \begin{split}
      \lambda \eta_{[0]}(D_{\mathcal{V}_A}+\mathcal{C}_A) &=
    \frac{1}{\sqrt{\pi}} \lambda \int_0^\infty \mathrm{TR}\left(
    (D_{\mathcal{V}_A}+\mathcal{C}_A)\exp(t
    (D_{\mathcal{V}_A}+\mathcal{C}_A)^2)\right)\,\frac{dt}{\sqrt{t}} \\
    &=     \frac{1}{\sqrt{\pi}}  \int_0^\infty \mathrm{TR}\left(
    \lambda_\sharp\left((D_{\mathcal{V}_A}+\mathcal{C}_A)\exp(t
      (D_{\mathcal{V}_A}+\mathcal{C}_A)^2)\right)\right)\,\frac{dt}{\sqrt{t}} \\
  &=     \frac{1}{\sqrt{\pi}}  \int_0^\infty \mathrm{TR}\left(
    (D_{\mathcal{V}^\lambda_B}+\mathcal{C}_B^\lambda)\exp(t
      (D_{\mathcal{V}_B^\lambda}+\mathcal{C}^\lambda_B)^2)\right)\,\frac{dt}{\sqrt{t}} \\
  &= \eta_{[0]}(D_{\mathcal{V}^\lambda_B}+\mathcal{C}^\lambda_B).    \end{split}
  \end{equation*}
\end{proof}

\section{Twists with finite dimensional representations}
\label{sec:twists-with-finite}

As a particular case of the principles explained in the previous
subsection, we consider a discrete group $\Gamma$ and a finite
dimensional unitary representation $\lambda\colon\Gamma\to
U(d)\subset M(d\times d, \complexs)$ . This induces a morphism of
$C^*$-algebras (we also call it $\lambda$) $\lambda\colon
C^*\Gamma\to M(d\times d, \complexs)$. Here, $C^*\Gamma$ is the
maximal $C^*$-algebra of $\Gamma$, and the extension of $\lambda$
follows from the universal property of $C^*\Gamma$.

Given, as in Section \ref{sec:class-mishch-fomenko},  a
classifying map $u\colon M\to B\Gamma$ defined on a closed
manifold $M$, we obtain the corresponding Mishchenko-Fomenko
bundle $\mathcal{L}$ and the twisted Dirac operator $D_{\L}$; thus
in this case $$A=C^* \Gamma\,,\quad \mathcal{V}_A\equiv \L :=
\tilde{M}\times_\Gamma C^* \Gamma\,, \quad B=M(d\times d,
\complexs)$$ Let $V_\lambda$ be the flat {\it vector bundle}
$\tilde{M}\times_\lambda \CC^d$ and let $D_\lambda$ be the Dirac
operator twisted by $V_\lambda$. Then it is easy to see that the
bundle $\mathcal{V}_B^\lambda$ introduced in the previous
subsection is equal to the direct sum of $d$-copies of $V_\lambda$
and that $D_{\mathcal{V}^\lambda_B}$ is the diagonal operator
$D_\lambda\otimes I_d$, with $I_d$ the $d\times d$-identity
matrix. If $\C\in \Psi^{-\infty}_{C^* \Gamma}$ is a trivializing
perturbation, then $\lambda_\natural \C= C_\lambda\otimes I_d$,
with $C_\lambda\in \Psi^{-\infty}$ and $D_\lambda+C_\lambda$
invertible. Let $\tr\colon M(d\times d,\CC)\to\complexs$ be the usual trace and
consider $\tr_d:= d^{-1}\tr$; let $\tr_\lambda:= \tr_d\circ
\lambda\colon  C^* \Gamma_{\rm ab}\to \CC$.
Then, by Lemma \ref{lem:stable_eta_preserved},
\begin{equation}\label{eq:aps_equals_tr}
\tr_\lambda (\eta_{[0]}(D_{\L} + \C)=\eta(D_\lambda+C_\lambda)
\end{equation}


Let us go back to closed manifolds and let $X=M\sqcup (-M^\prime)$
with $M$ and $M^\prime$ homotopy equivalent. Let us fix a unitary
representation $\lambda\colon \Gamma\to U(d)$ and let
$D_{X,\lambda}^{{\rm sign}}$ be  signature operator twisted by the
flat vector bundle associated to $\lambda$. Let $0<\epsilon$ such
that ${\rm spec} (D^{{\rm sign}}_{\lambda,X})\cap
(-\epsilon,\epsilon)\subseteq\{0\}$. Let $\phi\in C^\infty_c (\RR,[0,1])$
be a real function equal to 1 on $[-\epsilon, \epsilon]$ and equal to
0 on $(-\infty,-2\epsilon]\cup [2\epsilon,\infty)$. Let $\D^{{\rm
sign}}_X$ be the Mishchenko-Fomenko signature operator on $X$. We
know from Section \ref{sec:stable-=-unstable} that in this situation
we can  construct a trivializing perturbation $\C$
of $\D^{{\rm sign}}_X$ such that $(\Id-\phi(\D^{{\rm
sign}}_X))\circ \mathcal{C}=0$. Thus, using (\ref{localize-at-0}),
we see that $\mathcal{C}_\lambda$ is a smoothing operator such
that $(\Id-\phi(D^{{\rm sign}}_{\lambda,X}))\circ C_\lambda=0$; by
Remark \ref{rem:discrete_case_much easier} we can now conclude as in
Section \ref{subsec:limits}, but without making use of the limiting procedure of
Section \ref{sec:proof-lemma-}, that
$\eta(D_\lambda)=
(\eta(D_\lambda+C_\lambda)+\eta(D_\lambda-C_\lambda))/2$. Using Remark
\ref{rem:exact_perturbed_eta_calc} we can even conclude that the
stable and unstable APS-rho
invariants coincide, without making use of the limiting procedure of
Section \ref{sec:proof-lemma-}.

\section{Classical $L^2$-invariants versus $\NeumannN\Gamma$-invariants}
\label{sec:class-l2-invar-1}

\subsection{Classical $L^2$-invariants}
\label{sec:class-l2-invar}

In Section \ref{sec:intro} we introduced (delocalized)
$L^2$-invariants by working directly on a normal covering $\tilde M$
of a closed
manifold $M$.
However, all the main arguments of this
paper are given  in terms of $A/\overline{\Commutator{A}{A}}$-valued
 invariants for Dirac
operators twisted by Hilbert $A$-module bundles.
We briefly called these invariants {\it degree zero invariants}.
 In this appendix,
we will explain how the classical
$L^2$-invariants can be derived from these degree zero invariants,
where the $C^*$-algebra $A$ in question is
$\NeumannN\Gamma$, and the bundle to twist with is the
Mishchenko-Fomenko line bundle $\NeumannN=\tilde
M\times_{\Gamma}\NeumannN\Gamma$.
This procedure is well known to the experts. In this appendix, we will
give a detailed description of it, since we are not aware of a place
in the
literature where this would be covered, in particular when dealing
with eta invariants (the focus of interest here).

\subsection{Twisting with $l^2\Gamma$}
\label{sec:twist-with-l2gamma}

The passage from the covering situation to the twist with the
$\NeumannN\Gamma$ bundle is achieved by using an intermediate step:
instead of twisting with $\N:=\tilde M\times_\Gamma\NeumannN\Gamma$ we
twist with $\H:=\tilde M\times_\Gamma l^2\Gamma$.

The typical fiber of this bundle is isomorphic to $l^2\Gamma$, which
we consider as an $\NeumannN\Gamma$-Hilbert space (in the notation of
\cite{Schick03}).

We start with giving the basic definitions, repeated from
\cite[Section 8]{Schick03}.
\begin{definition}\label{def:of_A_Hilbert_space}
  A \emph{finitely generated projective $\NeumannN\Gamma$-Hilbert space} $V$ is a Hilbert
  space together with a right action of $\NeumannN\Gamma$ such that $V$ embeds
  isometrically preserving the $\NeumannN\Gamma$-module structure as a direct
  summand into $l^2(\Gamma)^n$ for some $n$.

  A (general) \emph{$\NeumannN\Gamma$-Hilbert space} $V$ satisfies the same
  conditions a finitely generated projective $\NeumannN\Gamma$-Hilbert space does,
  with the exception that $l^2(\Gamma)^n$ is replaced by $H\tensor l^2(\Gamma)$
  for some Hilbert space $H$ with trivial $\NeumannN\Gamma$-action (the tensor product has to be
  completed).
\end{definition}

\begin{definition}
 An $\NeumannN\Gamma$-Hilbert space morphism is a bounded $\NeumannN\Gamma$-linear map
    between two $\NeumannN\Gamma$-Hilbert spaces. If it is an isometry for the
    Hilbert space structure, it is called $\NeumannN\Gamma$-Hilbert space isometry.
 A $\NeumannN\Gamma$-Hilbert space bundle $\H$ on a space $X$
    is a locally
    trivial bundle of $\NeumannN\Gamma$-Hilbert spaces, the transition functions
    being $\NeumannN\Gamma$-Hilbert space isometries. A smooth structure is given by
    a trivializing atlas where all the transition functions are
    smooth.

    If the fibers are finitely generated projective $\NeumannN\Gamma$-Hilbert space,
    the bundle is called a \emph{finitely generated projective
      $\NeumannN\Gamma$-Hilbert space bundle}.
\end{definition}

\begin{lemma}\label{lem:sections are A Hilbert space}
  The $L^2$-sections of a $\NeumannN\Gamma$-Hilbert space bundle $\mathcal{W}$ on a
  Riemannian manifold $X$ form themselves a $\NeumannN\Gamma$-Hilbert space.
\end{lemma}
\begin{proof}
  Compare \cite[Lemma 8.10]{Schick03}
\end{proof}

\begin{example}
  Essentially the only example we will be dealing with is the
  $\NeumannN\Gamma$-Hilbert space $l^2(\Gamma)$. This is obviously
  finitely generated projective. Since the left
  regular representation acts by $\NeumannN\Gamma$-Hilbert space
  isometries, $\H:=\tilde M\times_\Gamma l^2(\Gamma)$ is a finitely
  generated projective $\NeumannN\Gamma$-Hilbert space bundle on $M$.

  The trivial connection on $\tilde M\times l^2(\Gamma)$ descends to a
  canonical   flat connection on $\tilde M\times_\Gamma l^2(\Gamma)$.
\end{example}

\begin{remark}\label{rem:bundle_completion}
    Note that $\NeumannN\Gamma$ is canonically a subset of $l^2(\Gamma)$
  (the inclusion given by $b\mapsto 1\cdot b$), and $l^2(\Gamma)$ is
  the completion of $\NeumannN\Gamma$ with respect to the inner
  product $\innerprod{b_1,b_2}:= \tr_\Gamma(b_1^*b_2)$. In the same
  way, $\tilde M\times_\Gamma l^2(\Gamma)$ is the fiberwise completion
  of $\tilde M\times_\Gamma\NeumannN\Gamma$.  For all of this, compare
 \cite[Section 8.6]{Schick03}
\end{remark}

Given a Dirac type operator $D\colon C^\infty(M,E)\to
C^\infty(M,E)$ on $M$, we can now also form the
twisted Dirac operator $D_\H$ as usual.



  If $M$ is closed, this is an elliptic differential operator of order $1$ on
 $\NeumannN\Gamma$-Hilbert space bundles in the
  sense of \cite{BFKM(1996)}. In any case, it extends to an
  unbounded operator on $L^2(E\tensor \H)$.
 The main point, as observed in \cite[Section 8]{Schick03} is now the
  following:
\begin{proposition}
The operator $D_\H$ maps the sections of
    the subbundle
    $E\tensor \NeumannN$ of $E\tensor \H$ to sections of
    $E\tensor \NeumannN$, and the restriction is exactly the
    operator $D_\NeumannN$. This holds for smooth sections and all
    kinds of completions. Moreover, all the functions of the
    operator $D_\H$ we are considering (like $D_\H\exp(-tD_\H^2)$) restrict
    to the corresponding functions of $D_\NeumannN$.

Vice versa, we get $D_\H$ by
    applying the procedure of completion
    of \cite[Section 8]{Schick03} to
    $D_\NeumannN$.

The above properties together imply  that
  there is a canonical bijection of algebras
  \begin{equation*}
    \{\psi(D_\N)\mid \psi \text{ is a Schwartz function}\}\text{ and }
    \{\psi(D_\H)\mid \psi\text{ is a Schwartz function}\}.
  \end{equation*}
From now on, making use of this identification, we will write $D_\H$
and $D_\N$ interchangeably. Put it differently, we now have defined the
invariants $\eta_\tau(D_\H):=\eta_\tau(D_\N)$,
$\eta_{\tau_{<g>}}(D_\H):=\eta_{\tau_{<g>}}(D_\N)$.
\end{proposition}

\subsection{Translation between $\tilde M$ and $M\times_\Gamma
  l^2\Gamma$}
\label{sec:transl-betw-tilde}

In this subsection, we are going to explain why it is equivalent
to consider
$\tilde D$ acting on $\tilde E \to \tilde M$ and
$D_\H$, $\H=\tilde M\times_\Gamma l^2(\Gamma)$, acting
on $E\otimes\H\to M$.
In fact, there is a dictionary that allows to pass back and forth between objects on
  the covering and twisted objects.
For the sake of completeness we indicate the
  constructions. Other accounts (where certain aspects are explained
  in more detail) can be found e.g.~in \cite[Section
  III]{LottI}, \cite[Section
  3.1]{SchickICTP} and \cite[Example 3.39]{SchickBologna}.

The translation is summarized in the following table.

  \begin{tabular}[t]{c|c|c}\\
no.&    $\tilde M$ & $\cdot \tensor \H$\\ \hline
 1&   $L^2(\tilde M, \tilde E)$ & $L^2(M,E\tensor \H)$\\
  2&  $\{s\in C^\infty(\tilde M,\tilde E)\mid \sum_{\gamma\in\Gamma} \abs{s(\gamma
      x)}^2<\infty\;\forall x\in \tilde M\}$ &  $C^\infty(M,E\tensor
    \H)$ \\
  3&  $\tilde D$ & $D_\H \;(=D_\N)$\\
4&   $\tilde D\exp(-t\tilde D^2)$ & 
   $D_\H\exp\left(-t D^2_\H\right)$\\ 
  5& $\Tr_{<g>}$ & $\tau_{<g>}\circ \TR$\\
  6& $\Tr_{(2)}$ & $\tau_\Gamma\circ \TR$\\
7& $\eta_{<g>}(\tilde D)$ & $\eta_{\tau_{<g>}}(D_\H)$.
\end{tabular}

Recall that for an operator $\tilde A$ given by the smooth integral
kernel $\tilde k$, by definition
\begin{equation}
\Tr_{<g>}(\tilde A)=\sum_{\gamma\in [g]}\int_{\tilde M/\Gamma} \tr_x k(x,\gamma
   x)\;dx,
\label{eq:trace_def}
 \end{equation}
and in particular $\Tr_{(2)}(\tilde A)= \int_{\tilde M/\Gamma} \tr_x
\tilde k(x,x)\;dx$.

\medskip
 We give some explanation how to translate between the two sides.

\noindent (1) A section $s$ of $\tilde E$ corresponds to the section $\hat
    s$ of $E\tensor \H$ with $\hat s(x) = \sum_{\gamma\in\Gamma}
    s(\gamma\tilde x)\tensor [\tilde x,\gamma]$, where $\tilde
    x\in\tilde M$ is an arbitrary
    lift of $x\in M$ along the covering projection. We identify the fibers $E_x$ and $\tilde
    E_{\gamma\tilde x}$. Moreover, by definition $\H_x= \Gamma\tilde x
    \times_\Gamma l^2(\Gamma)$. This construction is well defined by
    the very
    definition of the twisted bundle $\H$, with fiber identified with
    $l^2(\Gamma)$ using the chosen lift $\tilde x$.

\noindent (2) The above identification defines an isometry of the spaces of
    $L^2$-sections. Moreover, it is compatible with the
    $\Gamma$-actions. $L^2(\tilde M,\tilde E)$ is well known to be an
    $\NeumannN\Gamma$-Hilbert space, and the identification is an
    identification of $\NeumannN\Gamma$-Hilbert spaces. In
    addition, it preserves smoothness and continuity,
    where the condition as given in the table is used to really get a
    section of $E\tensor \H$.

\noindent (3) The operators $\tilde D$ and $D_\H$ are conjugated to each
    other under the isomorphism of the section spaces. This follows
    from their local definition as follows. For a small
    connected neighborhood $U$ of $x\in M$, we can choose a lift
    $\tilde U$, a connected neighborhood of a lift $\tilde x$, such
    that there is a unique section $U\to \tilde U$ of the restriction
    of the covering $\tilde M\to M$ to $U$, and then $y\mapsto [\tilde
    y,\gamma]$ is a by definition of the connection a \emph{flat}
    section of $\H|_U$ for each $\gamma\in\Gamma$. Consequently, using
    the identification of (2) and this flatness, we see that
    $\tilde D(\tilde s)$  corresponds on the set $U$ to
    \begin{equation*}
    \sum_{\gamma\in\Gamma}\tilde D\tilde s(\gamma \tilde x) \tensor
    [\tilde x,\gamma] = D_\H (\sum_{\gamma\in\Gamma} s(\gamma\tilde
    x)\tensor [x,\gamma]),
  \end{equation*}
  i.e.~to $D_\H$ applied to the section corresponding to $\tilde s$.

\noindent (4) Since the self-adjoint unbounded operators $\tilde D$ and $D_\H$ are
    unitarily equivalent, the same is true for all bounded measurable
    functions of them, using functional calculus. In particular, this
    is the case for $\tilde D\exp(-t\tilde D^2)$, but also for any
    other bounded measurable function $\phi\colon\reals\to\reals$.

\noindent
  (5) Choose a subset $U\subset M$ such that $M\setminus U$ has
    measure zero and such that the restriction of the covering $\tilde
    M\to M$ to $U$ is trivial. If we choose an appropriate lift
    of $U$
    then $\tilde M|U \iso U\times \Gamma$. This induces a trivialization
    $\H|_U\iso U\times l^2(\Gamma)$. Using this, we identify
$L^2(U,E\tensor \H|_U) =
    L^2(U,E|_U)\tensor l^2(\Gamma)$. 

    On the other hand, using the corresponding trivialization of the
    covering $\tilde M|U\iso U\times \Gamma$ we get the identification
    $L^2(\tilde U, \tilde E|U)\iso
    L^2(U,E|U)\tensor l^2(\Gamma)$, and our unitary isomorphism defined
    above becomes the identity under these identifications.

    So, in this picture, also the operators $\tilde D$ and $D_\H$ are
    identical. It remains to show that the two ways to define the trace are
    equal. We take traces of operators like $S=\tilde De^{-t\tilde
      D^2}$, which have integral kernels $\tilde k(\tilde x,\tilde y)$
    on $\tilde M\times\tilde
    M$. Note that it is well known that these operators are of
    $\Gamma$-trace class, by
    the results of \cite{Atiyah-covering}. They are also well known
    to be smoothing operators in the
    $\NeumannN\Gamma$-Mishchenko-Fomenko calculus, and therefore
    $\TR(S)$ is defined.

Using the identification $\tilde U=U\times\Gamma$, $\tilde
    k(u,g,v,h)$ has the property that for a section $s$ of $E|_U$ and
    $u\in U$, $g,\gamma\in \Gamma$
    \begin{equation*}
      S(s(\cdot)\tensor \gamma)(u,g) =\int_U \tilde k(u,g,x,\gamma)
      s(x) \,dx.
    \end{equation*}

  On the other hand, if we look at the integral kernel $k_\H$ on
  $M\times M$ of $S$, now considered as obtained from $D_\H$,
  using again the identification $L^2(U,E|_U\tensor \H)= L^2(U,
  E|_U)\tensor l^2(\Gamma)$, we get for a section $s$ of $E|_U$
  \begin{equation*}
    S(s(\cdot)\tensor \gamma)(u) = \int_U k_\H(u,x)
    (s(u)\tensor \gamma) \,dx = \sum_{g\in\Gamma} \int_U k_\H(u,x)(\gamma,g)
    (s(u)\tensor \gamma)\tensor g \,dx.
  \end{equation*}
  where we write the homomorphism $k_\H(u,x)$ from $E_x\tensor
  l^2(\Gamma)$ to $E_u\tensor l^2(\Gamma)$ as a matrix
  $\left(k_\H(u,x)(\gamma,g)\right)_{\gamma,g\in \Gamma}$ of
  homomorphisms
 from $E_x$ to $E_u$, using
  the orthonormal basis $\Gamma$ of $l^2(\Gamma)$.

  Since these are two representations of the same operator $S$, it
  follows that
  \begin{equation}
k_\H(u,x)(\gamma,g)=\tilde
k(u,\gamma,x,g).
\label{eq:equality_of_matrices}
\end{equation}

  Now, by definition of $\Tr_{<g>}$, we get
  \begin{equation}\label{eq:formula_for_tr}
    \Tr_{<g>}(S) = \sum_{h\in<g>} \int_U \tr_x \tilde k(x,1,x,h) \,dx.
  \end{equation}

  On the other hand, by the definition of $\tau_{<g>}$ in terms of a
  matrix decomposition of elements of
  $\NeumannN\Gamma\subset\boundedops(l^2\Gamma)$, using the basis
  $\Gamma$ of $l^2(\Gamma)$,
  \begin{equation}\label{eq:formula_for_tau}
    \tau_{<g>}\TR(S) = \int_U\tau_{<g>}\tr^{alg}_x(k_\H(x,x))\,dx =
    \int_U \sum_{h\in<g>}\tr_x (k_\H(x,x)(1,h))\,dx.
  \end{equation}
  From Equation \eqref{eq:equality_of_matrices}, the right hand sides
  of Equations \eqref{eq:formula_for_tr} and
  \eqref{eq:formula_for_tau} are equal. Consequently,
  \begin{equation*}
      \Tr_{<g>}(S) = \tau_{<g>}\circ \TR(S),
    \end{equation*}
    as claimed.

  \noindent (6)Since $\Tr_{(2)}=\Tr_{<1>}$, this is a special case of (5).



  \noindent (7) This is a direct consequence of (4) and
    of (5), using the definition of the $\eta$-invariants.


    For more details, in particular with respect to the delocalized
    trace, compare also \cite[Section 4]{MR2000k:58039} and
    \cite[Section 3]{LottI}

We want to single out the relevant proposition which is used in the
body of this paper, using the results of the above table and of
Section :
\begin{proposition}\label{prop:compute_L2_eta}
   Let $M$ be a compact manifold with $\Gamma$-covering $\tilde M$, and
   $D$ a Dirac type operator on $M$. Let $<g>$ be a finite conjugacy
   class in $\Gamma$. If $\tilde D$ is the lift of $D$ to $\tilde M$,
   and $D_\NeumannN$ is the twist of $D$ with the Mishchenko-Fomenko
   line bundle $\NeumannN=\tilde M\times_\Gamma \NeumannN\Gamma$, then
   \begin{equation*}
       \eta_{<g>}(\tilde D) = \eta_{\tau_{<g>}}(D_\NeumannN).
     \end{equation*}
     In particular, for $g=e$ we compute the classical $L^2$-eta
     invariant
     \begin{equation*}
       \eta_{(2)}(\tilde D) = \eta_{\tau_\Gamma}(D_\NeumannN).
     \end{equation*}
\end{proposition}

\begin{remark}\label{rem:a-la-Lott}
There is an alternative route to Proposition
\ref{prop:compute_L2_eta}, due to John Lott, and based
on particular properties of the heat-kernel. We briefly explain the
main ideas, referring to \cite{LottI} for the details.

There is a sequence of inclusions of algebras
\begin{equation}\label{sequence-inclusion}
\B^\omega_\Gamma\subset \B^\infty_\Gamma\subset C^*_r \Gamma\subset
\NeumannN\Gamma\,,
\end{equation}
where $\B^\omega_\Gamma$ is the algebra of functions on $\Gamma$
which are exponentially rapidly decreasing and where
$\B^\infty_\Gamma$ is the
Connes-Moscovici algebra.
A generalized Dirac operator can be twisted with the corresponding
flat bundles
\begin{equation}\label{sequence-bundles}
\V^\omega:=\tilde M\times _\Gamma \B^\omega_\Gamma\;,\quad  \V^\infty:=\tilde M\times
_\Gamma\B^\infty_\Gamma
\;,\quad
\V:= \tilde M\times _\Gamma C^*_r \Gamma\;,\quad
\NeumannN:=\tilde M\times _\Gamma \NeumannN\Gamma.
\end{equation}
producing $\D^\omega$, $\D^\infty$, $\D:=D_{\V}$, $D_{\NeumannN}$
There are obvious compatibility conditions for these operators, coming
from the inclusions of vector spaces
\begin{equation*}
C^\infty(M,E\otimes\V^\omega)\subset C^\infty (M,E\otimes
\V^\infty)\subset C^\infty (M,E\otimes\V) \subset
C^\infty (M,E\otimes \NeumannN)
\end{equation*}
We also have natural  inclusions for the corresponding
Mishchenko-Fomenko calculi
\begin{equation*}
\Psi^*_{\B^\omega_\Gamma}\xrightarrow{j_\omega}
\Psi^*_{\B^\infty_\Gamma}\xrightarrow{j_\infty} \Psi^*_{C^*_r \Gamma}\xrightarrow{j_r}\Psi^*_{\NeumannN\Gamma}
\end{equation*}
For the smoothing operators we have natural traces and a
  commutative diagram
\begin{equation*}
\begin{CD}
\Psi^{-\infty}_{\B^\omega_\Gamma} @>{j_\omega}>> \Psi^{-\infty}_{\B^\infty_\Gamma}
@>{j_\infty}>> \Psi^{-\infty}_{C^*_r  \Gamma} @>{j_r}>> \Psi^{-\infty}_{\NeumannN\Gamma}\\
@V{\TR_\omega}VV                               @V{\TR_\infty}VV    @V{\TR}VV @V{\TR_{\NeumannN\Gamma}}VV\\
(\B^\omega_\Gamma)_{{\rm ab}} @>>> (\B^\infty_\Gamma)_{{\rm ab}}
@>>> (C^*_r  \Gamma)_{{\rm ab}}  @>>> (\NeumannN\Gamma)_{{\rm ab}}
\end{CD}
\end{equation*}
Using Proposition 6 and Proposition 7 in \cite{LottI} we have that
\begin{equation*}
\tau_{<g>} \TR_\omega (\D^\omega e^{-t(\D^\omega)^2})=
\Tr_{<g>}
(\tilde D e^{-t(\tilde D)^2})\,.
\end{equation*}
Since from the commutative diagram and the definition of
$\tau_{<g>}$ we clearly have $$\tau_{<g>} \TR_\omega  ( \D^\omega e^{-t(\D^\omega)^2}) =
\tau_{<g>} \TR_{\NeumannN \Gamma}  ( D_\NeumannN e^{-t D_\NeumannN ^2}),$$
we see that
Proposition \ref{prop:compute_L2_eta} is proved once again.
Notice that we have also established that
$$
\tau_{<g>} \TR_\infty (\D^\infty e^{-t(\D^\infty)^2})
=\Tr_{<g>}
(\tilde D e^{-t(\tilde D)^2})\,.$$

\end{remark}

\end{appendix}

{\small
\bibliographystyle{plain}
\bibliography{rho.bib}
}

\end{document}